\documentclass[DIV16,twocolumn,10pt,headsepline,abstract=true]{scrartcl}

\usepackage[english]{babel}
\usepackage[utf8x]{inputenc}
\usepackage[intlimits]{amsmath}
\usepackage{graphicx}
\usepackage{float}
\usepackage{booktabs}
\usepackage{xspace}
\usepackage{url}
\usepackage{subfig}
\usepackage{color}
\usepackage[colorlinks,%
pagebackref=false, % bibliography -> text
linktocpage=true, % toc etc: make page number active (not name)
plainpages=false, % distinguish roman and arabic pagenumbers
bookmarksopen=true,%
bookmarksnumbered=true,%
pdfauthor={Adrian Q.T. Ngo},%
pdftitle={Numerical solution of steady-state groundwater flow and solute
  transport problems: Discontinuous Galerkin based methods compared to the Streamline Diffusion approach},%
pdfsubject={},%
pdfkeywords={reactive modelling, capillary fringe, bacterial growth, adhesion, transport},%
]{hyperref}        % clickabe references
\usepackage[expproduct=cdot,alsoload={hep,binary}]{siunitx}
\usepackage{nicefrac}
\usepackage{ctable}
\usepackage{listings}
\usepackage[square]{natbib}
\usepackage{amssymb}
\usepackage{stmaryrd}
\usepackage{algorithm}
\usepackage{algorithmicx}
\usepackage{algpseudocode}
\usepackage{tikz}
\usetikzlibrary{chains, shapes.misc}
\usepackage{pgfplots}
\usepackage{pgfplotstable}
\usetikzlibrary{backgrounds}
\usetikzlibrary{pgfplots.groupplots}

\newcommand{\ra}[1]{\renewcommand{\arraystretch}{#1}}
\DeclareMathOperator{\erfc}{erfc} 

\usepackage{geometry}

\addtokomafont{pagehead}{\small}
\addtokomafont{section}{\normalsize}
\addtokomafont{subsection}{\normalsize\normalfont}

\newunit{\years}{\text{year}}

\title{Numerical solution of steady-state groundwater flow and solute
  transport problems: Discontinuous Galerkin based methods compared to the Streamline Diffusion approach}

\author{A.Q.T.~Ngo$^1$, P.~Bastian$^1$, O.~Ippisch$^{1,2}$}
\date{}
\publishers{\small \begin{flushleft}\noindent $^1$University of Heidelberg, Interdisciplinary Center for Scientific Computing, Im Neuenheimer Feld 368, 69120 Heidelberg, Germany\\\vspace{2mm}
$^4$Clausthal University of Technology, Department of  Mathematics, Erzstra{\ss}e 1, 38678 Clausthal-Zellerfeld, Germany \\\vspace{2mm}
  Email: adrian.ngo@iwr.uni-heidelberg.de
\end{flushleft}
}

\renewcommand{\hat}{\widehat}

\newcommand{\lvecx}{\vec{x}} % location vector x
\newcommand{\lvecv}{\vec{v}} % v
\newcommand{\lvecq}{\vec{q}} % q
\newcommand{\matId}{\mathbf{Id}} % Id matrix
\newcommand{\matD}{\mathcal{D}}   % Scheidegger Dispersion tensor

\newcommand{\grad}{\nabla}  % the gradient is a vector
\newcommand{\divergence}{\textnormal{div}}

\newcommand{\injectiontime}{T_{\mbox{\tiny{inj}}}}
\newcommand{\injectionRate}{\tilde{w}_{\hbox{\tiny{inj}}}}
\newcommand{\extractionRate}{\tilde{w}_{\hbox{\tiny{ext}}}}
\newcommand{\injectionWell}{W_{\hbox{\tiny{inj}}}}
\newcommand{\extractionWell}{W_{\hbox{\tiny{ext}}}}

\newcommand{\GammaInflow}{\Gamma_{\mbox{\footnotesize{--}}}}
\newcommand{\GammaOutflow}{\Gamma_{\mbox{\footnotesize{+}}}}
\newcommand{\GammaCharacteristic}{\Gamma_{\mbox{\tiny{0}}}}

\definecolor{farbe1}{RGB}{255,0,0}
\definecolor{farbe2}{RGB}{0,0,255}
\definecolor{farbe3}{RGB}{0,255,0}
\definecolor{farbe4}{RGB}{192,192,192}

\begin{document}

\maketitle

\begin{abstract}\noindent {
    In this study, we consider the simulation of subsurface flow and solute
    transport processes in the stationary limit. In the convection-dominant case,
    the numerical solution of the transport problem may exhibit non-physical
    diffusion and 
    under- and overshoots.
    For an interior penalty discontinuous Galerkin (DG) discretization, we present 
    a $h$-adaptive refinement strategy and, alternatively, a new efficient
    approach for reducing numerical under- and overshoots using a diffusive
    $L^2$-projection. Furthermore, we illustrate an efficient way of solving the
    linear system arising from the DG discretization.
    In $2$-D and $3$-D examples, we compare the DG-based methods to the streamline
    diffusion approach with respect to computing time and their ability to resolve
    steep fronts.
  }
\end{abstract}

\section{Introduction}

Natural flowing conditions that are nearly at steady-state over a long period of time 
(i.e. several days during which a tracer experiment is being conducted) 
can be expected, if at all possible, in a confined aquifer.
Due to very small dispersivities and very small molecular diffusion, 
(non-reactive) solute transport in groundwater is convection-dominated.
Considering temporal moments of its concentration
% in the coupled forward model 
leads to the \textit{steady-state singularly perturbed convection-diffusion equation}.
%The transport of a conservative tracer can be described by 
%a singularly perturbed convection-diffusion problem.
% the stationary linear convection-diffusion equation which is formally 
% elliptic, but nearly first order hyperbolic in the convection-dominated case.
The numerical solution of this kind of equations %boundary value problems 
has a long tradition. 
Due to the fact that a linear monotonicity preserving scheme can
be at most first-order accurate (Godunov's Theorem),
all existing schemes suffer from 
a trade-off between numerical diffusion (too much smearing) and
spurious oscillations (under- and overshoots) near internal or boundary layers
where the gradient of the solution is very large (steep fronts).
From a practical point of view, the decision whether the one
or the other deficit is tolerated has to be made, leading to the appropriate
choice of a numerical scheme.

Amongst the vast literature on the subject, the books of
\citet{Roos:2008} and \citep{Kuzmin:2010} provide excellent overviews of
state-of-the art classes of schemes.
We give a short overview 
of the most prominent methods, list their main advantages and
disadvantages, hereby drawing on the results presented by \citet{John:2011},
who have worked on a special $2$-D problem (Hemker problem).
%(with a constant velocity field and a constant diffusion coefficient and  
%where the solution contains a boundary layer that transitions into a
%characteristic layer at some critical points)
They compared the numerical solutions 
at specific cut lines with respect to the
size of maximal under- and overshoots, the
width of smeared internal layers and the performance in computing
time, revealing important properties of the different schemes.

\begin{itemize}
\item The Scharfetter-Gummel scheme is a first-order finite volume scheme. It
is efficient and oscillation-free, but the solution is strongly smeared at
layers. A higher order extension is not available.

\item The Streamline Diffusion finite element method (SDFEM), also known as
Streamline Upwind Petrov-Galerkin (SUPG) method, adds a residual-based
stabilization term to the standard Galerkin
method \citep{Brooks:1982}. 
SDFEM belongs to the less time-consuming methods that are capable of 
resolving the steep fronts well.
Due to its simplicity, it 
has been the mainstream approach for decades and a standard method
in the hydrogeologists' community \citep{Cirpka:2001,Nowak:2006,Gordon:2000,Couto:2008,Bear:2010} where our
practical application originates from.
However, the optimal choice of a user-defined 
stabilization parameter is an open question. 

\item The Continuous Interior Penalty (CIP) method \citep{Roos:2008} adds a symmetric stabilization
term to the standard Galerkin method that penalizes jumps of the gradient
across faces (edge stabilization technique).
It introduces connections between unknowns of neighboring mesh cells and leads
to a discretization with a wider matrix stencil.
Compared to the SDFEM method it is in general less performant.

\item Spurious Oscillations at Layers Diminishing (SOLD) methods, originally
developed in \citep{Hughes:1986} and further investigated
in \citep{John:2007,John:20072,John:2008}, suppress oscillations caused by
SDFEM by adding a further stabilization term introducing diffusion orthogonal to
streamlines (crosswind diffusion). This term is in general non-linear.
Therefore, a non-linear equation has to be solved for a linear problem.
Furthermore, the stabilization term contains another user-defined parameter whose
optimal choice might become difficult for complicated problems.
SOLD methods are capable of reducing numerical oscillations 
at a higher computational cost. The larger the stabilization
parameter, the better the reduction. However, non-linearity also increases and 
the iterative non-linear solver might not converge \citep{John:2011}.

\item Algebraic Flux Correction (AFC) is a general approach to design high
resolution schemes for the solution of time-dependend transport problems
that ensure the validity of the discrete maximum principle \citep{Kuzmin:2006,Kuzmin:2010}.
Whereas the aforementioned stabilization methods modify the bilinear form of a
finite element method (FEM), AFC methods modify the linear system arising from a
FEM discretization by adding discrete diffusion to the system
matrix and appropriate anti-diffusive fluxes to the right hand
side.
The anti-diffusive fluxes are non-linearly dependent on the computed solution.
Depending on whether the algebraic constraints are being imposed on the semi-discrete or
the fully discrete level, flux limiters of TVD-type (total variation
diminishing) or FCT-type (flux corrected transport) can be
constructed.
Only the FEM-TVD schemes can be used to solve the steady-state 
convection-diffusion equation directly. Using FEM-FCT schemes, 
a pseudo time stepping to the stationary limit of the associated time-dependent problem
would deliver the steady-state solution \citep{Kuzmin:2006}.
A suitable linearization technique for the anti-diffusive fluxes exists only
for the FEM-FCT scheme \citep{Kuzmin:2009}.
According to the studies in \citep{John:2008_2,John:2009}, FEM-FCT schemes yield
qualitatively the best solution and, beyond that, the linear FEM-FCT scheme is
efficient. The authors recommend the linear FEM-FCT for the
solution of instationary problems.
Linearized AFC-based methods for the solution of the stationary transport
problem are not available.
% In summary, AFC-based methods for the solution of the stationary transport problem 
% can be very accurate, but time consuming.

\item
Discontinuous Galerkin (DG) methods use piecewise polynomials, that are not
required to be continuous across faces, to approximate the solution.
%On a cuboid mesh, 
The total number of degrees of freedom on a structured mesh with cuboidal cells is 
% equal to the
% number of mesh cells multiplied by a power function of the polynomial degree.
$O(n\cdot (k+1)^d)$ where $n$ is the number of mesh cells, $d$ is the
dimension of the domain and $k$ is the polynomial degree.
Compared to a continuous Galerkin FEM method, a DG method using the same
polynomial space on the same structured mesh requires more unknowns.
This disadvantage is balanced by a long list of advantages that has made DG
increasingly attractive in computational fluid dynamics in the last decade:
%Due to the localization of test functions to single mesh cells, 
DG methods are readily parallelizable, lead to discretizations with 
compact stencils (i.e. the unknowns in one mesh cell are only connected to the unknowns 
in the immediate neighboring cells),
a higher flexibility in mesh design (non-conforming meshes are possible in adaptive
h-refinement) and the availability of different polynomial degrees on
different mesh cells (adaptive p-refinement).
Furthermore, DG schemes satisfy the local, cell-wise mass balance which is a crucial
property for transport processes in a porous medium.
They are particularly well-suited for problems with discontinuous
coefficients and effectively capture discontinuities in the solution.
In the comparative study by \citep{John:2011}, the DG method gives the best
results regarding
sharpness of the steep fronts and produces small errors with respect to
reference cut lines, whereas under- and overshoots are larger than those
produced by the SDFEM method.
For the discretization of first-order hyperbolic problems, upwinding is
incorporated into the formulation of DG schemes, evading the need for 
user-chosen artificial diffusion parameters.
The books of \citet{KanschatBook:2008}, \citet{Riviere:2008} and \citet{Pietro:2012} 
offer a comprehensive introduction to this class of methods.
\end{itemize}
For time-dependent problems, where explicit time stepping schemes combined with 
finite volume or DG discretizations can be used, 
slope limiters may be constructed from the solution of one time-step to 
preserve monotonicity in the following time-step.
To the best of our knowledge, for the immediate solution of stationary problems, 
a post-processing technique of this type is not available.

In the simulation of many applications (e.g. biochemical reactions or combustion),
the concentration of a species
must attain physical values (numerical under- and overshoots are not accepted) 
although the position of the plume may be allowed to be inaccurate.
%
%But in the context of a parameter estimation scheme, 
%which allows for measurement errors, 
%a small amount ($\approx 5\%$) of spurious oscillations in the solution is tolerable 
%whereas the correct localization of steep fronts is of primary interest.
%
Our practical application stems from the field of geostatistical
inversion \citep{Cirpka:2001} in which the hydraulic conductivity is estimated on
the basis of indirect measurements of related quantities such as the tracer
concentration or arrival time.
The computed solution of the transport problem is used for a pointwise comparison
with real measurements. Hereby, the parameter estimation scheme
allows for measurement errors, i.e. a small amount ($\approx 5\%$) of spurious
oscillations in the solution is tolerable. 
By contrast, the correct localization of steep fronts is of primary interest.
% Furthermore, we are interested in the best possible  solution that can
% be resolved on the same structured mesh on which the hydraulic conductivity is
% resolved.

%The comparative study made in \citep{John:2011} has shown 
%that 
%SDFEM belongs to the less time-consuming methods that are capable of 
%resolving the steep fronts well.

%Another observation was that the DG method gives the best results regarding
%sharpness of the
%steep fronts and produces small errors with respect to reference cut lines,
%whereas spurious oscillations are larger than those produced by the
%SDFEM method.

For high resolution $3$-D simulations of real-world applications, 
direct sparse solvers are limited by their memory
consumption. 
The main purpose of the stabilization term in the SDFEM method is not only to 
provide a solution with bounded under- and overshoots
but also to improve the iterative solvability of the 
linear system arising from the SDFEM discretization.
For an upwind scheme applied to a first-order hyperbolic problem, it is
well-known that numbering the unknowns in a fashion that follows 
approximately the direction in which information is propagated
will improve the performance and stability of iterative linear solvers
of ILU or Gauss-Seidel type \citep{Bey:1997,Hackbusch:1997,Reed:1973}.

Motivated by the requirements for our practical application,
we have chosen a DG-based method to solve the %groundwater flow and 
solute transport problem.
The remaining part of this paper is structured as follows:
Section \ref{paragraph:model} introduces the steady-state groundwater flow and
solute transport equations,
in Section \ref{paragraph:disc} we consider two combinations of
discretizations for their numerical solution: FEM/SDFEM and CCFV/DG.
\begin{itemize}
\item We present two approaches to reduce the under- and overshoots
      of the DG-solution of the transport problem:
    \begin{enumerate}
    \item The first approach (Section \ref{paragraph:adaptive}) uses 
    $h$-adaptive hanging-nodes $1$-irregular refinement on a cuboidal 
    axis-parallel
    mesh based on the residual error estimator
    by \citet{Schoetzau:2009} combined with an error-fraction marking
    strategy.
    \item The second approach (Section \ref{paragraph:diffusiveL2projection}) is a
    diffusive $L^2$-projection of the DG
    solution into the continuous Galerkin finite element subspace. It works directly
    on the structured coarse mesh.
    \end{enumerate}
    \item In addition, we have implemented an efficient way to solve the 
      linear system for the DG discretization
      iteratively by exploiting a downwind cell-wise 
      numbering of unknowns {\it before} the stiffness matrix is
      assembled (Section \ref{paragraph:renumbering}).
\end{itemize}

\noindent Numerical studies are presented in
      Section \ref{paragraph:numerical_tests}.
      The DG method is compared to a first order SDFEM
      implementation within the same code, i.e. the performance of the two
      different discretizations can be compared on the same computational 
      grid using the same linear solver. 
      The results are summarized in Conclusion \& Outlook section.

The numerical software used to perform the simulations 
is written in C++ and 
based on the libraries of the {\it Distributed and Unified Numerics Environment \tt DUNE}
\citep{dunegridpaperI:08,dunegridpaperII:08,dune-web-page}
and the finite element discretization module DUNE-PDELab
(\url{www.dune-project.org/pdelab}).
DUNE offers a structured parallel grid (YASP) and interfaces to
the unstructured grids UG \citep{ug} and {\tt DUNE-ALUGrid} \citep{DednerKN14}.

\section{Model equations}\label{paragraph:model}
%---------------------------------------
%Steady-state flow
%---------------------------------------
The physical models describing flow and transport processes in a confined
aquifer are well developed and can be
found in the textbooks \citep{Bear:2010,DeMarsily:1986} or in the lecture note
\citep{Roth:2012}.
\subsection{Groundwater flow}\label{s:GWE}
For convenience, we assume $\Omega\subset\mathbb{R}^d,~ d\in\{2,3\}$ to be a 
rectangular cuboid in which the boundary is
subdivided into a Dirichlet boundary
($\Gamma_{D}$) and a Neumann boundary ($\Gamma_{N}$) section.
%with polygonal or polyhedral boundary.
We consider the {\it steady-state groundwater flow equation} 
\begin{equation}\label{eqn:GWE}
  \nabla \cdot\big( -K \nabla \phi \big) = \injectionRate - \extractionRate \quad \mbox{in} \quad \Omega
\end{equation}
subject to the boundary conditions:
\begin{equation}\label{eqn:GWE_BC}
\begin{array}{rclcl}
  \phi & = & \hat{\phi}_{D}   &\mbox{on}& \Gamma_{D}\\
  \vec{n}\cdot (-K\nabla \phi) & = & 0  &\mbox{on}& \Gamma_{N}
\end{array}
\end{equation}

\noindent in which $\vec{n}(\lvecx)$ is the unit outer normal vector, 
$K(\lvecx)>0$ is the spatially variable, but locally isotropic 
hydraulic conductivity [$m/s$],
$\phi(\lvecx)$ is the hydraulic head $[m]$ and
the source terms on the right hand side % $\tilde{w}(\lvecx)$ %is the volumetric flux per unit volume and 
prescribe the rates
(volumetric flux per unit volume [$1/s$])
of injection and extraction wells 
$\injectionWell, \extractionWell \subsetneq \Omega$:
\begin{equation}\label{eqn:wellRates}
  \begin{split}
    & \injectionRate(\lvecx) \left \{
      \begin{array}{cl}
        > 0  & \lvecx \in \injectionWell \\
        \\
        = 0  & \lvecx \in \Omega\backslash \injectionWell\\
      \end{array}
    \right.
    \\
    & \extractionRate(\lvecx) \left \{
      \begin{array}{cl}
        > 0  & \lvecx \in \extractionWell \\
        \\
        = 0  & \lvecx \in \Omega\backslash \extractionWell\\
      \end{array}
    \right.
  \end{split}
\end{equation}
This is an elliptic equation for which
the coefficient $K$ may be highly variable.
% 
% 
% The domain boundary can be subdivided into two types:
% Dirichlet boundary sections $\Gamma_{\mbox{\footnotesize{D}}} \subset \partial
% \Omega$ on which fixed values $\phi_0$
% of the hydraulic head {\it must} be prescribed and 
% Neumann boundary sections $\Gamma_{\mbox{\footnotesize{N}}} = \partial \Omega
% \backslash \Gamma_{\mbox{\footnotesize{D}}}$ on which the normal flux $q_0$ may be prescribed:
%
%
%
The fluid motion is induced by the head distribution.
The volumetric flux
% \footnote{The volumetric flux $\vec{q}$, sometimes called
%  the Darcy velocity, is related to the pore water velocity $\vec{v}$ via
%  $\vec{q}=\theta\vec{v}$ where $\theta$ is the porosity. \citep{Roth:2012}}
%[$LT^{-1}$]
is given by Darcy's law:
\begin{equation}\label{eqn:Darcy}
  \vec{q} = - K\nabla \phi \quad \mbox{in} \quad \Omega.
\end{equation}
$\vec{q}$ is sometimes called the Darcy velocity. It is related to the pore water velocity $\vec{v}$ via
$\vec{q}=\theta\vec{v}$
where $\theta$ is the porosity.
We work with $\lvecq$ and the porosity 
$\theta$ is supposed to be constant.

\subsection{Subsurface solute transport}\label{s:GWE}
A conservative tracer %like a nonreactive dye tracer 
used to track flow motion has no influence on the flow itself.
Its concentration $c(t,\lvecx)$ [$kg/m^{3}$] is described by the transient
convection-diffusion-reaction equation \citep{Bear:2010}
\begin{equation}\label{eqn:TPE_transient}
  \begin{split}
    \dfrac{\partial(\theta c)}{\partial t} 
    + \divergence ( - \matD \nabla c ) + \lvecq~ \grad c 
    & = 
    \injectionRate ~\tilde{c}_{\hbox{\tiny{inj}}}
    - \extractionRate ~c
    \\
    & \quad \mbox{in~} (0,T] \times \Omega
  \end{split}
\end{equation}
in which $\tilde{c}_{\hbox{\tiny{inj}}}(t,\lvecx)$ is the concentration 
at the injection well and
\begin{equation}\label{eqn:DispersionTensor}
  \matD = \theta \matD_{\text{S}}
\end{equation}
is the dispersion tensor [$m^{2}/s$] 
given by \citep{Scheidegger:1961}:
\begin{equation}\label{eqn:Scheidegger}
  \matD_{\text{S}} = \big(\alpha_{\ell} -\alpha_t \big) ~\dfrac{\vec{v}\cdot\vec{v}^T}{\|~\vec{v}~\|_2} ~+~  \bigg( \alpha_t~\|\vec{v}\|_2 + D_m \bigg) ~\matId, 
\end{equation}
where $\alpha_{\ell}$ and $\alpha_t$ are the longitudinal and transversal
dispersivities [$m$],
$D_m$ is the mo\-le\-cu\-lar diffusion coefficient [$m^{2}/s$]
and $\matId$ is
the identity matrix.
For the transport equation, we distinguish three types of boundaries,
inflow, outflow and characteristic boundary:
\begin{equation}\label{eqn:TPE_BT}
  \begin{array}{rcl}
    \GammaInflow &=& \{ \vec{x}\in\partial\Omega :
    \lvecq(\vec{x})\cdot \vec{n}(\vec{x}) < 0 \} \\
    \GammaOutflow &=& \{ \vec{x}\in\partial\Omega :
    \lvecq(\vec{x})\cdot \vec{n}(\vec{x}) > 0 \} \\
    \GammaCharacteristic &=& \{ \vec{x}\in\partial\Omega :
    \lvecq(\vec{x})\cdot \vec{n}(\vec{x}) = 0 \} \\
  \end{array}
\end{equation}

Tracer in a constant concentration may enter over a fixed time
period $\injectiontime>0$ through the
injection well or somewhere on the inflow boundary:
\begin{equation}\label{eqn:TPE_BC1_transient}
  \begin{array}{rcll}
    c(t,\lvecx) & = & \hat{c}_D(\lvecx) & \text{on}\quad \GammaInflow\\
  \end{array}
\end{equation}

On the whole boundary $\partial\Omega=\GammaInflow \cup \GammaOutflow \cup \GammaCharacteristic$,
we assume that the flux is non-diffusive for convection dominant transport:
\begin{equation}\label{eqn:TPE_BC2_transient}
  \begin{array}{rcll}
    \vec{n}\cdot(-\matD\nabla c) & = & 0 & \text{on}\quad \GammaInflow \cup \GammaOutflow \cup \GammaCharacteristic
  \end{array}
\end{equation}
This implies the no-flux condition for impermeable boundaries:
\begin{equation}\label{eqn:TPE_BC3_transient}
  \begin{array}{rcll}
    \vec{n}\cdot(-\matD\nabla c + \lvecq c) & = & 0 & \text{on}\quad \GammaCharacteristic
  \end{array}
\end{equation}

\noindent We are interested in the steady-state solution ($\partial c/\partial t = 0$), or
more adequately, in the zeroth or first order temporal moments
of the resident concentration \citep{Harvey:1995}.
In all these cases, the differential equation takes the form of the 
{\it steady-state convection-diffusion equation}:

\begin{equation}\label{eqn:TPE}
  \divergence ( - \matD \nabla u + \tau \lvecq~ u) + \mu u = \tilde{s}
  \qquad \mbox{in} \quad \Omega
  ~.
\end{equation}
The reaction coefficient $\mu\in\mathbb{R}$ may be used as a sink term within extraction wells.
The source term $\tilde{s}$ may be used to describe the behavior at the injection well.
The boundary conditions read:
\begin{equation}\label{eqn:TPE_BC}
  \begin{array}{rcll}
    u & = & \hat{u}_D & \text{on}\quad \Gamma_{-}\\
    \vec{n}\cdot(-\mathbf{D}\nabla u) & = & 0 & \text{on}\quad \partial\Omega\\
  \end{array}
\end{equation}
For $\mu' = \mu - \divergence(\lvecq)$, equation (\ref{eqn:TPE})
can be written in non-conservative form as 
\begin{equation}\label{eqn:TPE_convective}
  \divergence(- \mathbf{D} \nabla u) + \lvecq~ \nabla u + \mu' u = \tilde{s}
  \qquad \mbox{in} \quad \Omega
  ~.
\end{equation}

\section{Discretizations}\label{paragraph:disc}

\subsection{Preliminary definitions}\label{paragraph:PreDefs}
Let $\{\mathcal{T}_{h_{\nu}}\}_{\nu\in\mathbb{N}}$ be a family of structured
or adaptively refined 
meshes (comprised of axis parallel cuboidal cells) %as in (\ref{eqn:smesh})
that we get from a successive refinement of an initially structured mesh.
Each $\mathcal{T}_{h_{\nu}}$
forms a partitioning of $\Omega$ into $n_{\nu}$ disjoint cells (mesh elements), this means
\begin{equation}\label{eqn:Partitioning}
\begin{split}
  \mathcal{T}_{h_{\nu}} &= \{t_j^\nu\}_{j=0,...,n_\nu-1} \\
  \overline{\Omega} &= \bigcup\limits_{j=0}^{n_{\nu}-1} \overline{t_j^{\nu}}, 
  \quad t_i^{\nu} \cap t_j^{\nu} = \emptyset \quad \forall i \ne j.
\end{split}
\end{equation}
The variable $\nu$ indicates the refinement level.
As refinement proceeds, the meshsize 
\begin{equation}\label{eqn:meshsize}
h_{\nu} = \max\limits_{t \in\mathcal{T}^*_{\nu}} \left\{
\max\limits_{\vec{x},\vec{y}\in \overline{t}} \|\vec{x}-\vec{y}\| 
\right\}
\end{equation}
tends to $0$.
$\mathcal{T}^*_{\nu}$ is understood to be the subset of
$\mathcal{T}_{h_{\nu}}$ that contains only
the finest level cells. 
To keep notation readable, we write $h$ instead of $h_\nu$,
$\mathcal{T}_{h}$ instead of $\mathcal{T}_{h_{\nu}}$ and $n$ instead of $n_\nu$
when it is clear or irrelevant which refinement level $\nu$ we are considering.

The hydraulic conductivity is resolved on the structured
mesh $\mathcal{T}_{h_0}$. It is described by a cell-wise constant function:
\begin{equation}\label{eqn:Kh}
  \begin{split}
    & K_h(\lvecx) = K(\lvecx_{t}) = \exp( Y(\lvecx_{t}) ) \\
    & \forall 
    \lvecx\in t,~ \lvecx_{t} \textnormal{~is the center the cell~} t \in \mathcal{T}_{h_0}
  ~.
  \end{split}
\end{equation}
Inside the wells, the hydraulic conductivity is supposed to be very high:
\begin{equation}\label{eqn:easyWellModel}
  \begin{split}
  K_h(\lvecx)=1.0 & \qquad  \forall~ \lvecx\in t, t \in \mathcal{T}_{h_0}
  \quad \textnormal{with} \\ 
  &
  ~t\cap \injectionWell\ne\emptyset \quad \textnormal{or} \quad 
  t\cap \extractionWell\ne \emptyset
  ~.
  \end{split}
\end{equation}
The hydraulic head distribution and the Darcy velocity (\ref{eqn:Darcy}) 
are computed on the same mesh.
The transport equation (\ref{eqn:TPE}), whose convective part is prescribed by
the Darcy velocity, may be solved either on the same mesh or on a hierarchy of
adaptively refined meshes based on $\mathcal{T}_{h_0}$, i.e. a cell of a subsequently refined mesh $\mathcal{T}_{h_{\nu+1}}$ is always a subset of a cell in $\mathcal{T}_{h_{\nu}}$.

In our practical application, the estimated hydraulic conductivity 
fields %(\ref{eqn:cokrigingSolution}) in the inversion algorithm
have the smoothness of a Gaussian variogram model. 
The meshsizes are chosen in such a way that they resolve the correlation
lengths well. % are many times greater than the meshsize.
Thus, the flow field on the mesh $\mathcal{T}_{h_0}$ can be regarded as
sufficiently accurate 
and we consider adaptive mesh refinement only for the solution of the transport problem.
Finite elements on cuboids $t\in\mathcal{T}_h$ are based on the polynomial space
\begin{equation}\label{eqn:Qkspace}
  \begin{split}
  \mathbb{Q}_k^d = \bigg\{ p(\lvecx) = \sum\limits_{0\le \alpha_i \le k\atop 1\le i \le d}
  \gamma_{\alpha_1,...,\alpha_d} & \cdot x_1^{\alpha_1} \cdots x_d^{\alpha_d} ,
  \\
  & \lvecx\in t \bigg\} 
  \end{split}
\end{equation}
with maximal degree $k$ in each coordinate direction.
Discontinuous Galerkin (DG) approximations are based on the 
broken polynomial space
\begin{equation}\label{eqn:DGspace}
  \begin{split}
    W_{h,k} = W_{h,k}(\Omega,\mathcal{T}_h) = \bigg \{ u \in & L^2(\Omega) :
    u|_{t} \in \mathbb{Q}_{k}^d
    \\
    & \forall ~ t\in \mathcal{T}_h \bigg \}
    ~.
  \end{split}
\end{equation}
We restrict ourselves to the case where the maximal polynomial degree $k$ is constant for all
cells:
\begin{equation}
\dim(W_{h,k}) = n \cdot \dim( \mathbb{Q}_k^d ) = n \cdot (k+1)^d
~.
\end{equation}
The continuous polynomial space
\begin{equation}\label{eqn:FEMspace}
  \begin{split}
    V_h 
    ~ = ~ 
    V_h(\Omega,\mathcal{T}_h) 
    ~ = ~ 
    \bigg \{ u \in & C^0(\overline{\Omega}) : u|_{t} \in \mathbb{Q}_1^d \\
    & \quad \forall ~ t\in \mathcal{T}_h \bigg \}
  \end{split}
\end{equation}
is used to describe the standard Galerkin FEM and the streamline diffusion method.
For a structured mesh with $n=\prod\limits_{i=1}^{d}n_i$ cells, its dimension is 
$\prod\limits_{i=1}^{d}(n_i+1)$.

\medskip
All cell-wise or face-wise defined integrals
that will occur in the following numerical schemes
are integrals over products of at least two polynomials of order $k$.
A Gaussian quadrature rule of order $k+1$ guarantees the exact evaluation of polynomials of order $2k+1$.

\subsection{FEM / SDFEM}
Let $V_h^0 = \big\{ u\in V_h ~:~ u_{|\partial\Omega} = 0 \big\}$ and assume $\hat{\phi}_{D}$ to be a piecewise linear
approximation of the Dirichlet b.c. in (\ref{eqn:GWE_BC}).
%
%
% Standard FEM:
%
%
%
The standard Galerkin \textbf{FEM method} (cf. \citep{Elman:2005}) for solving (\ref{eqn:GWE})\&(\ref{eqn:GWE_BC}) reads:
\begin{center}
  \framebox[0.49\textwidth]{
    \begin{minipage}[c]{0.45\textwidth}
      \smallskip
      Find $\phi_h\in V_h^{\text{D}} = \big\{ u\in V_h ~:~ u_{|\Gamma_D} = \hat{\phi}_{D} \big\}$ such that 
      \begin{equation}\label{eqn:FEM}
        \sum\limits_{t\in \mathcal{T}_h} \big(K_h\nabla \phi_h,\nabla v_h \big)_{0,t} = \tilde{w}_h(v_h)
        \quad \forall ~v_h \in V_h^0
        ~.
      \end{equation}
      % \smallskip
    \end{minipage}
  }
\end{center}
The discrete source term on the right-hand side is 
\begin{equation}\label{eqn:FEMsourceterms}
  \begin{split}
  \tilde{w}_h(v_h)
  = &
  \sum\limits_{\stackrel{t\in \mathcal{T}_h}{t\cap \injectionWell\ne \emptyset}} \big( \injectionRate, v_h \big)_{0,t} 
  \\
  & - 
  \sum\limits_{\stackrel{t\in \mathcal{T}_h}{t\cap \extractionWell\ne \emptyset}} \big( \extractionRate, v_h \big)_{0,t} 
  \end{split}
\end{equation}
where $K_h \in W_{h_0,0}$ as defined in (\ref{eqn:Kh}). 
The discrete Darcy velocity can be computed by direct pointwise
evaluation of gradients of the polynomial basis on each cell $t\in\mathcal{T}_{h_0}$:
\begin{equation}\label{eqn:FEM_Darcy}
  \vec{q}_h = -K_h \nabla\phi_h.
\end{equation}
%
%
% Streamline Diffusion:
%
%
%
\noindent The \textbf{SDFEM method} (cf. \citep{Brooks:1982}) for solving
(\ref{eqn:TPE_convective}) % \&(\ref{eqn:TPE_BC}) 
reads: 
\begin{center}
  \framebox[0.49\textwidth]{
    \begin{minipage}[c]{0.45\textwidth}
      \smallskip
      Find $u_h\in V_h^{\text{--}} = \big\{ u\in V_h ~:~ u_{|\Gamma_{-}} = \hat{u}_D \big\} $ such that 
      \begin{equation}\label{eqn:SDFEM}
        \begin{split}
          & \sum\limits_{t\in \mathcal{T}_h} 
          \bigg\{ \big(~\mathbf{D} \nabla u_h, ~ \nabla v_h ~\big)_{0,t} 
          \\
          & + \big(~ \vec{q}_h~\nabla u_h + \mu'_h u_h ,~ v_h + \delta^{\text{SD}}_t \cdot \vec{q}_h~\nabla v_h ~\big)_{0,t}
          \bigg\} \\
          & = 
          \sum\limits_{t\in \mathcal{T}_h} 
          \bigg\{ \big( ~ \tilde{s}_h, ~ v_h + \delta^{\text{SD}}_t \cdot \vec{q}_h~\nabla v_h~ \big)_{0,t}
          \bigg\} \\
          &
          \hspace{2cm} \forall ~v_h \in V_h^0.
        \end{split}
      \end{equation}
%      \smallskip
    \end{minipage}
  }
\end{center}
$\mu'_h$ and $\tilde{s}_h$ evaluate the corresponding terms of equation (\ref{eqn:TPE_convective}) at quadrature points.
The matrix $\mathbf{D}\in\mathbb{R}^{d\times d}$ is the discretized version of the dispersion tensor $\matD$ in (\ref{eqn:DispersionTensor}) which is depending on $\lvecv=\lvecq_h/\theta$.
%
%
%\noindent 
The \textbf{stabilization parameter} determined by 
\begin{equation}\label{SDFEM_parameter}
  \delta^{\text{SD}}_t = \dfrac{h}{2\|\vec{q}_h\|_2}\cdot\zeta( \mathcal{P}^t_h )
\end{equation}
is used to tune the amount of artificial
diffusion depending on the magnitude of the \textbf{mesh P\'eclet number} $\mathcal{P}^t_h$.
If $\mathbf{D}=\varepsilon\cdot\mathbf{Id}~(\varepsilon>0)$, % in (\ref{eqn:BVP_2}),
the general definition of the mesh P\'eclet number is
\begin{equation}\label{eqn:Peclet}
  \mathcal{P}^t_h = \frac{1}{2}\cdot\dfrac{\|\vec{q}_h\|_2\cdot h_t}{\varepsilon}.
\end{equation}
%where $h$ is a characteristic length of the mesh.
For the Scheidegger dispersion tensor (\ref{eqn:Scheidegger}), the effective mesh
P\'eclet number according to \citep{Cirpka:2001} is 
\begin{equation}\label{eqn:Peclet_longitudinal}
  \mathcal{P}^t_h = \frac{1}{2}\cdot\dfrac{\|\vec{q}_h\|_2\cdot
    h_t}{\alpha_{\ell}\cdot\|\vec{q}_h\|_2 + \theta D_m}.
\end{equation}
There is a large variety of definitions for the function $\zeta$ in
the literature. The original choice in \citep{Brooks:1982} is 
\begin{equation}
  \zeta(\mathcal{P}^t_h) = \coth(\mathcal{P}^t_h) - 1/\mathcal{P}^t_h.
\end{equation}
We go for the more efficient approximation (cf. \citep{Elman:2005})
%\begin{equation}
%  \zeta(\mathcal{P}^t_h) = \left\{
%  \begin{array}{cr}
%    1 - \dfrac{1}{\mathcal{P}^t_h} & \qquad\text{if}\quad \mathcal{P}^t_h > 1\\
%    0 & \qquad\text{if}\quad \mathcal{P}^t_h \le 1.\\
%  \end{array}
%  \right.
%\end{equation}
\begin{equation}
  \zeta(\mathcal{P}^t_h) = \max \big\{ ~ 0,  ~1 - 1/\mathcal{P}^t_h ~\big\}.
\end{equation}

\subsection{CCFV / DG}
\subsubsection{Preliminary definitions}\label{s:preDG}
\noindent The following notation is inspired by the presentation of \citet{Ern:2008}.
For a given mesh $\mathcal{T}_h$, % of the domain $\Omega$, 
each cell $t\in\mathcal{T}_h$ has a cell-center $\vec{x}_t$ and a 
$d$-dimensional cell-volume $|t|$.
Given two neighboring cells $t^-$ and $t^+$ in $\mathcal{T}_h$, 
an interior face or interface $f$ is defined as the intersection of their boundaries
$\partial t^- \cap \partial t^+$. To be more precise, we write $t^-_{f}=t^-$ and
$t^+_{f}=t^+$.
The \textbf{unit normal vector} $\vec{n}_{f}$ to $f$ is assumed to be oriented from $t^-_{f}$ to
$t^+_{f}$.
In the same manner, $f$ is called a boundary face if there exists a $t\in\mathcal{T}_h$ such that
$f=\partial t \cap \partial \Omega$ and we write $t^-_{f}=t$. 
In this case, $\vec{n}_{f}$ is chosen to be
the unit outer normal to $\partial \Omega$.
We denote by $\mathcal{E}_h$ the \textbf{set of interior faces} and by $\mathcal{B}_h$
the \textbf{set of boundary faces}.
%For $d=3$ a face is a 2-D boundary side of a cell, for $d=2$ a face is a 1-D edge.
%
For a face $f\in\mathcal{E}_h\cup\mathcal{B}_h$, we denote by $\vec{x}_{f}$ its 
midpoint (face center), whereas $\vec{x}_{f-}$ is the center of the cell
$t_{f}^-$. %``containing'' $f$. 
For $f\in\mathcal{E}_h$, we denote by
$\vec{x}_{f+}$ the center of 
% the ``neighboring'' cell 
$t_{f}^+$. 
$|f|$ is the $(d-1)$-dimensional volume of the face $f$.

A function $u\in W_{h,k}$ is in general double valued on internal faces $f\in\mathcal{E}_h$.
There, we set
\begin{equation}\label{eqn:doublevalued}
\begin{split}
%u^{\pm}_{f} = 
u_f^{\pm}(\vec{x}) = 
\lim\limits_{\varepsilon\rightarrow 0{\pm}} u(\vec{x}+\varepsilon\vec{n}_{f})
%\quad \text{and} \quad 
%u^+_{f}(\vec{x}) = 
%\lim\limits_{\varepsilon\rightarrow 0+} u(\vec{x}+\varepsilon\vec{n}_{f}),
\qquad 
\vec{x}\in f.
\end{split}
\end{equation}
The jump across $f$ and the arithmetic mean value on $f$ are given by 
%For any $\vec{x}\in e$ we define the {\bf jump of a function} $u\in W_{h,k}$ as 
\begin{equation}%\label{eqn:jump}
\begin{split}
\llbracket u \rrbracket_{f} = u^-_{f} - u^+_{f}
%\llbracket u \rrbracket_{f}(\vec{x}) = u^-_{f} - u^+_{f}
%\qquad \forall ~ \vec{x}\in f, \quad f\in \mathcal{E}_h
%$$ 
%and the mean value of $u$ on $f$ as 
%$$
\quad \text{and} \quad
\big< u \big>_{f}(\vec{x}) = \frac{1}{2} \biggl(
u^-_{f} + u^+_{f}
\biggr)
%\qquad \forall ~
%\vec{x}\in f, \quad f\in \mathcal{E}_h.
\end{split}
\end{equation}
respectively.
Following convention, the definition of these terms is extended to the boundary
$\partial \Omega$ by:
\begin{equation}
  \llbracket u \rrbracket_{f}(\vec{x}) = \big< u \big>_{f}(\vec{x}) = u (\vec{x}) 
  \qquad \forall ~
  \vec{x}\in f, \quad f\in \mathcal{B}_h.
\end{equation}
We will suppress the letter $f$ in subscripts if there is no ambiguity.

% In this paper, we seek for approximations of (\ref{eqn:GF}) and
% (\ref{eqn:TPE}) in the broken polynomial spaces $W_{h,p}$ which opens the
% door to a greater variety of numerical schemes.
%
%
%
% Two-point CCFV:
%
%
%
\subsubsection{Two-point flux cell-centered finite volume method (CCFV)}\label{paragraph:GWE_CCFV}
Using the function space $W_{h,0}$, the approximation of the boundary value
problem (\ref{eqn:GWE}) \& (\ref{eqn:GWE_BC})
is defined as follows:
\begin{center}
  \framebox[0.49\textwidth]{
    \begin{minipage}[c]{0.45\textwidth}
      \smallskip
      Find $\phi_h\in W_{h,0}$ such that 
      \begin{equation}\label{eqn:CCFV}
        \begin{split}
          a^{\text{FV}}( \phi_h, v_h ) & = \ell^{\text{FV}}(v_h) 
          \\
          & \quad \forall ~v_h \in W_{h,0}.
        \end{split}
      \end{equation}
    \end{minipage}
  }
\end{center}
The bilinear form $a^{\text{FV}}: W_{h,0}\times W_{h,0} \longrightarrow \mathbb{R}$ is defined by 
\begin{equation}\label{eqn:a_FV}
  \begin{split}
    a^{\text{FV}}( \phi, v ) ~ &= ~ 
    \sum_{f\in \mathcal{E}_h} 
    q_h^{\phi}(\vec{x}_{f}) \cdot \llbracket v \rrbracket_{f} \cdot |f| 
    ~ \\
    & + \sum_{f\in \mathcal{B}_h \cap \Gamma_D}
    q_h^{\phi}(\vec{x}_{f}) \cdot v(\vec{x}_{f-}) \cdot |f|
  \end{split}
\end{equation}
where 
\begin{equation}\label{eqn:CCFV_flux}
  q_h^{\phi}(\vec{x}_{f}) := \left\{ 
  \begin{array}{cl}
    - K_h^{\mbox{\tiny{eff}}}(\vec{x}_{f+},\vec{x}_{f-}) \cdot
    \dfrac{ \phi(\vec{x}_{f+}) - \phi(\vec{x}_{f-})}{\|\vec{x}_{f+} - \vec{x}_{f-}\|_2} 
    \\
    &\\
    \quad\text{for} \quad f \in \mathcal{E}_h, \\
    \\
    \\
    - K_h(\vec{x}_{f-}) \cdot \dfrac{ \hat{\phi}_{D}(\vec{x}_{f}) - \phi(\vec{x}_{f-})}{\|\vec{x}_{f} - \vec{x}_{f-}\|_2}
    \\
    &\\
    \quad \text{for}  \quad f \in \mathcal{B}_h \cap \Gamma_D .\\
  \end{array} \right.
\end{equation}
is a two-point finite difference approximation
of the normal component $\lvecq\cdot\vec{n}_{f}(\vec{x}_{f})$ of the
Darcy velocity $\lvecq=-K\nabla \phi$ through the
face $f$.
The discrete hydraulic conductivity $K_h \in W_{h_0,0}$ is defined as in (\ref{eqn:Kh}) and
\begin{equation}\label{eqn:harmonicAverage}
  K_h^{\mbox{\tiny{eff}}}(\vec{x}_{f+},\vec{x}_{f-}) 
  = 
  \dfrac{2 \cdot K_h(\vec{x}_{f+}) \cdot K_h(\vec{x}_{f-})}
        {K_h(\vec{x}_{f+}) + K_h(\vec{x}_{f-})}
\end{equation}
is the \textbf{harmonic average} of $K_h(\vec{x}_{f+})$ and $K_h(\vec{x}_{f-})$. 
In the definition (\ref{eqn:CCFV_flux}), we make use of the fact that the
face $f$ is perpendicular to the line connecting $\vec{x}_{f+}$ and $\vec{x}_{f-}$.

The linear functional $\ell^{\text{FV}}:W_{h,0}\longrightarrow\mathbb{R}$ is given by 
\begin{equation}\label{eqn:CCFVsourceterms}
  \begin{split}
    \ell^{\text{FV}}(v)
    = &
    \sum\limits_{\stackrel{t\in \mathcal{T}_h}{t\cap \injectionWell\ne \emptyset}} \injectionRate(\vec{x}_t) \cdot v(\vec{x}_t) \cdot |t| \\
    & - 
    \sum\limits_{\stackrel{t\in \mathcal{T}_h}{t\cap \extractionWell\ne \emptyset}} \extractionRate(\vec{x}_t) \cdot v(\vec{x}_t) \cdot |t| 
  \end{split}
\end{equation}

\noindent This approximation yields a cell-wise constant solution for which the
Dirichlet boundary values $\hat{\phi}_D$ are satisfied weakly.

%
%
%
% Flux reconstruction:
%
%
%
\subsubsection{Flux reconstruction}\label{paragraph:RT0}
So far, only the normal flux component $q_h^{\phi}(\vec{x}_e)$ from (\ref{eqn:CCFV_flux}) is available on the
midpoints of the faces of a cell $t\in\mathcal{T}_h$. But we need to 
evaluate the Darcy velocity on internal points %$\vec{x}\in \overline{t}$ 
of a cell. The simplest $H(\text{div})$-conforming flux reconstruction is
achieved using Raviart-Thomas elements of order $0$ \citep{BrezziFortin:1991,RaviartThomas:1977},
\begin{equation}\label{eqn:RT0_elements}
  \begin{split}
    \mathbb{RT}_0(t) = 
    \bigg\{
    \vec{\tau}(\vec{x}) & ~:~  \tau_i = 
    a_i + b_i x_i,
    \quad \lvecx \in t, \\
    & 1\le i\le d, \quad a_i, b_i \in \mathbb{R} \text{~fixed} 
    \bigg\},
  \end{split}
\end{equation}
for which the polynomial
space is defined by 
\begin{equation}\label{eqn:RT0_space}
  \begin{split}
    \mathbb{RT}_0(\Omega,\mathcal{T}_h) = 
    \bigg\{
    \vec{\tau} \in [L^2(\Omega)]^d &~:~
    \vec{\tau}_{|t} \in \mathbb{RT}_0(t) \\
    &
    ~\forall t\in\mathcal{T}_h
    \bigg\}.
  \end{split}
\end{equation}

\noindent
The discrete Darcy velocity 
\begin{equation}\label{eqn:RT0_Darcy}
\lvecq_h\in\mathbb{RT}_0(\Omega,\mathcal{T}_h)
\end{equation}
can be evaluated component-wise by a linear interpolation between the
normal fluxes on opposing face midpoints.
It can be shown that the reconstructed 
Darcy velocity $\lvecq_h\in\mathbb{RT}_0(\Omega,\mathcal{T}_h)$ is indeed in
$H(\text{div},\Omega)$ and satisfies the projection condition:
\begin{equation}\label{eqn:RT0_projection}
  \int\limits_{f}\big( ~\lvecq_h - \lvecq ~\big) \cdot \vec{n}_{f} %\cdot v
  ~ds = 0
  %\qquad \forall v \in \mathbb{P}_0(f)
\end{equation}
Compared to the direct evaluation (\ref{eqn:FEM_Darcy}), this
reconstructed Darcy velocity field is pointwise divergence-free if the
groundwater equation (\ref{eqn:GWE}) is free of any source or sink terms ($\injectionRate=\extractionRate=0$).
%
%
%
% DG method:
%
%
%
\subsubsection{Discontinuous Galerkin method (DG)}\label{paragraph:DG}
The symmetric weighted interior penalty (SWIP) method presented in \citep{Ern:2008} is
a robust discontinuous Galerkin method accounting for anisotropy and discontinuity in the
diffusion tensor of (\ref{eqn:TPE}).

For the discretization of the diffusive term, the authors have introduced a
scalar- and double-valued weighting function $\omega$ on internal
faces. The two values $\omega^-$ and $\omega^+$ are constructed
based on the double-valued diffusion tensor $\mathbf{D}$ with $\mathbf{D}^-$
and $\mathbf{D}^+$ defined element-wise following (\ref{eqn:doublevalued}).
%
%
%
%   weighted average by Ern et.al.:
%
%
%
Using the normal component of $\mathbf{D}^{\mp}$ across the face, namely
$\delta^{\mp} = \vec{n}_{f}\cdot \mathbf{D}^{\mp} \cdot
\vec{n}_{f}$, the weighting factors are defined as  
\begin{eqnarray}%\label{eqn:harmonic_coefficients}
  \omega^{-} = \dfrac{\delta^+}{\delta^- + \delta^+} \qquad \text{and} \qquad
  \omega^{+} = \dfrac{\delta^-}{\delta^- + \delta^+}.
\end{eqnarray}
Both factors are non-negative and add up to unity 
$\omega^- + \omega^+ = 1$.

\noindent For $v\in W_{h,k}$, the weighted average of the diffusive flux is defined as 
\begin{equation}%\label{eqn:weighted_avarage_diffusive_flux}
  \big< \mathbf{D}\nabla v \big>^{\omega} = \omega^- (\mathbf{D}\nabla
  v )^- +  \omega^+ (\mathbf{D}\nabla v )^+ .
\end{equation}
On the boundary face $f\in\mathcal{B}_h$, we set $\omega=1$ and $\big< \mathbf{D}\nabla v \big>^{\omega} = \mathbf{D}\nabla v$.

%
%
%
%   upwind term:
%
%
%

For the convective term, we choose 
% a conservative 
an \textbf{upwind flux}
formulation 
that is equivalent to the presentations in \citep{Georgoulis:2009}, 
in \S 4.2 of \citep{Riviere:2008} or in \S 4.6.2 of \citep{Pietro:2012}.
% where non-homogeneous Dirichlet boundary conditions are taken into account.
%
% Given a cell $t\in\mathcal{T}_h$, we denote by $\partial_{-}t$ the inflow part
% and by $\partial_{+}t$ the outflow part of the boundary $\partial t$.
%
For an internal face $f\in\mathcal{E}_h$ lying between two neighboring cells
$t_f^-$ and $t_f^+$, recall that the unit normal vector $\vec{n}_{f}$ is assumed to be
oriented from $t_f^-$ to $t_f^+$. For a boundary face $f\in\mathcal{B}_h$,
$\vec{n}_{f}$ is the unit outer normal.
By $u^{\text{upwind}}$ we denote the upwind value of a function $u\in W_{h,k}$.
For $\vec{x}\in f, f\subset \partial t\backslash\Gamma_-$, it is defined by 
$$
u^{\text{upwind}}(\vec{x}) = \left\{
  \begin{array}{cc}
    u^+(\vec{x}) & \qquad \text{if} \quad \lvecq_h(\vec{x})\cdot\vec{n}_{f}<0\\
    \\
    u^-(\vec{x}) & \qquad \text{if} \quad \lvecq_h(\vec{x})\cdot\vec{n}_{f}\ge 0
  \end{array}
  \right.
$$

%
%
%
%   DG scheme:
%
%
%
\noindent The higher order DG approximation\footnote{~DG($k$) indicates that the polynomial basis is from $\mathbb{Q}_k^d$.
~The colors in the following terms are 
{\color{farbe1}red for the convection terms}, 
{\color{farbe2}blue for the diffusion terms}, 
{\colorbox{farbe3}{green} for the reaction term} and 
{\colorbox{farbe4}{grey} for the source term}.} of the boundary value problem
(\ref{eqn:TPE})\&(\ref{eqn:TPE_BC}) reads:

\begin{center}
  \framebox[0.49\textwidth]{
    \begin{minipage}[c]{0.45\textwidth}
      \smallskip
      Find $u_h\in W_{h,k}$ such that 
      \begin{equation}\label{eqn:DG_variational}
        a^{\text{DG}}(u_h,v_h) = \ell^{\text{DG}}(v_h) \quad \forall ~v_h \in W_{h,k}.
      \end{equation}
    \end{minipage}
  }
\end{center}
The bilinear form is defined by 
\begin{equation}\label{eqn:DG_bilinear}
\footnotesize
  \begin{split}
    & a^{\text{DG}} (u,v) ~ = ~ \\
    & \sum\limits_{t\in\mathcal{T}_h} 
    \bigg\{~ 
    {\color{farbe2} % diffusion term in alpha_volume()
      (\mathbf{D}\nabla u,\nabla v)_{0,t}
    }
    {\color{farbe1} % convection term in alpha_volume()
      ~ - ~( u, \vec{q}_{h}\cdot\nabla v )_{0,t} 
    }
    { % reaction term in alpha_volume()
      ~ + ~ \colorbox{farbe3}{$( \mu_h u, v )_{0,t}$}
    }
    ~ \bigg\} \\
    & \hspace{3cm}{\footnotesize\textit{(cell terms)}} \\
    & ~ + ~ \sum\limits_{f\in\mathcal{E}_h} 
    \bigg\{~ 
    % 
    % diffusion and convection terms:
    {\color{farbe2}
      % 
      % standard IP term for the diffusion:
      \big( ~\gamma \llbracket u \rrbracket, ~ \llbracket v \rrbracket ~ \big)_{0,f} 
    }
    {\color{farbe1}
      % upwind flux:
      ~+~ \big(~ |\vec{n}_{f}\cdot \vec{q}_{h}|~ u^{\text{upwind}}, ~ \llbracket v \rrbracket ~ \big)_{0,f}
    }
    \\
    &
    {\color{farbe2}
      % diffusion term:
      - \big(~ \big< \vec{n}_{f} \cdot \mathbf{D}\nabla u \big>^{\omega} ,~ \llbracket v \rrbracket ~\big)_{0,f} 
      % 
      % (non-)symmetric IP term
      - ~ \big( ~ \big< \vec{n}_{f} \cdot \mathbf{D}\nabla v \big>^{\omega} ,~
      \llbracket u \rrbracket ~\big)_{0,f}
    }
    ~ \bigg\}
    \\
    & \hspace{3cm}{\footnotesize\textit{(interior face fluxes)}} \\
    & ~+~ \sum\limits_{f\in\mathcal{E}_h\cap\GammaOutflow} 
    {\color{farbe1}
      % upwind flux:
      \big(~ \vec{n}_{f}\cdot \vec{q}_{h}~ u, ~ v ~ \big)_{0,f}
    }
    \\
    & \hspace{3cm}{\footnotesize\textit{(convective outflow)}}\\ 
    & ~ + \sum\limits_{ f \in\mathcal{B}_h\cap \Gamma_{-}}
    \bigg\{~ 
    {\color{farbe2}
      \big(~\gamma (u - \hat{u}_D),~ v ~\big)_{0,f}
    }
    \\
    &
    \hspace{2cm}{\color{farbe2}
      -~ \big(~ 
      u - \hat{u}_D ,
      ~\vec{n}_{f} \cdot \mathbf{D}\nabla v
      ~\big)_{0,f}
    }
    ~\bigg\}\\
    & \hspace{3cm} {\footnotesize\textit{(Dirichlet B.C.)}}
  \end{split}
\end{equation}
The linear functional is given by
\begin{equation}
\footnotesize
  \begin{split}
    \ell^{\text{DG}}(v) = 
    \sum\limits_{t\in\mathcal{T}_h} & \colorbox{farbe4}{$ (\tilde{s}_h,v)_{0,t} $} \\
    &
    - \sum\limits_{ f \in\mathcal{B}_h\cap \Gamma_{-}}
    {\color{farbe1}
      \big(\vec{n}_{f}\cdot \vec{q}_{h}~ \hat{u}_D, ~ v \big)_{0,f}
    }
    \\
    & {\footnotesize\textit{(source term and Dirichlet B.C.)} }
  \end{split}
\end{equation}
$\mu_h$ and $\tilde{s}_h$ evaluate the corresponding terms of equation (\ref{eqn:TPE}) at quadrature points.
\noindent The parameter $\gamma$ penalizing discontinuity in the solution is given in
\citep{Bastian:2011} by
\begin{equation}\label{eqn:DGpenalty}
  \begin{split}
    \gamma = C_{\gamma} & \cdot\dfrac{D_{\text{eff}}\cdot k(k+d-1)}{h_{f}} \\
    & \forall~ f\in\mathcal{E}_h\cup(\mathcal{B}_h\cap\Gamma_{-})
  \end{split}
\end{equation}
where $C_{\gamma}>0$ is a constant to be chosen sufficiently large
($C_{\gamma}=10$ is usually enough). With this definition of $\gamma$ the 
user-chosen constant does not play a big role anymore in the
convection-dominated case.

\noindent Remember that for the discrete problem (\ref{eqn:DG_variational}),
$\mathcal{T}_h$ can be a non-conformingly refined cuboidal mesh. 
% For adaptive refinement, we will work with a $1$-irregular mesh.
For a face lying between two
possibly non-matching elements, we set
$h_{f}=\min\{h^-_{f},h^+_{f}\}$ 
where $h^{\mp}_{f} = \text{diam}(f\cap \partial t_{\mp})$ are face diameters.
On internal faces, the effective diffusivity is defined by the harmonic average
\begin{equation}\label{eqn:effective_diffusivity}
  D_{\text{eff}} = \dfrac{2 \delta^- \delta^+}{\delta^- + \delta^+}
\end{equation}
of the normal component of the diffusion tensor across the face.
On boundary faces, we set directly 
\begin{equation}
  D_{\text{eff}} = \vec{n}_{f} \cdot \mathbf{D}  \cdot \vec{n}_{f}.
\end{equation}
\section{Adaptive mesh refinement}\label{paragraph:adaptive}
We are interested in reducing numerical oscillations globally.
Adaptive mesh refinement is a particularly promising way of achieving this goal.

Residual-based a-posteriori error estimators offer the advantage of small
evaluation cost because they are based on local residual terms.
We consider two existing 
$h$-adaptive versions for the above mentioned discretization schemes of the transport
problem.
%Let $u_h\in V^-_{h}$ be the solution of (\ref{eqn:SDFEM}) or 
%$u_h\in W_{h,k}$ be the solution of (\ref{eqn:DG_variational})
%respectively.
Let $u_h$ be the SDFEM solution (\ref{eqn:SDFEM}) or 
the DG solution (\ref{eqn:DG_variational}) respectively.
The residual-based error estimator described by 
\citet{Verfuerth:1998,Verfuerth:2005},
developed for the finite element and the SDFEM discretization of the steady-state
convection-diffusion equation, is based on the following local error
indicator:
For each element $t\in\mathcal{T}_h$, the local error indicator $\eta_t^2$
is given by the sum of two terms, 
\begin{equation}
  \eta_t^2 = \eta_{R_t}^2 + \eta_{R_f}^2, 
\end{equation}
an element residual term\footnote{Note that $\Delta u_h$ can be omitted for
  the polynomial degree $k=1$.}
\begin{equation}\label{eqn:ree_1}
\eta_{R_t}^2 = \dfrac{h_t^2}{\varepsilon} ~\big\| s_h + \varepsilon \Delta u_h -
\vec{q}_{h} \cdot \nabla u_h - \mu_h' u_h \big\|^2_{L^2(t)}
\end{equation}
and a face residual term 
\begin{equation*}
\eta_{R_f}^2 = \frac{1}{2} \sum\limits_{f\in\partial t\backslash \Gamma}
\dfrac{h_f}{\varepsilon} ~\big\| ~\llbracket  \vec{n} \cdot (\varepsilon \nabla u_h) \rrbracket~ \big\|^2_{L^2(f)}.
\end{equation*}

\noindent To these two terms, \citet{Schoetzau:2009} added a third term 
measuring jumps of the solution on internal or inflow
boundary faces
% between two elements
\begin{equation*}
\begin{split}
\eta_{J_f}^2 & = \frac{1}{2} \sum\limits_{f\in\partial t\backslash \Gamma}
\bigg( \frac{\gamma\varepsilon}{h_f} + \frac{h_f}{\varepsilon} \bigg)
\big\|  \llbracket u_h \rrbracket \big\|^2_{L^2(f)}
\\ & + \sum\limits_{f\in\partial t \cap \Gamma_{-}}
\bigg( \frac{\gamma\varepsilon}{h_f} + \frac{h_f}{\varepsilon} \bigg)
\big\|  (u_h - \hat{u}_D) \big\|^2_{L^2(f)}
\end{split}
\end{equation*}

\noindent to construct a residual-based error estimator
for the interior penalty DG discretization scheme.
Hence, the local error indicator $\eta_t^2$ for each element
$t\in\mathcal{T}_h$
is given by the sum of three terms, 
\begin{equation}\label{eqn:eta}
  \eta_t^2 = \eta_{R_t}^2 + \eta_{R_f}^2 + \eta_{J_f}^2.
\end{equation}
In our concrete application, we choose $\varepsilon = \min\limits_{1\le i \le d}\{ \mathbf{D}_{ii} \}$
%in all the three indicator terms.
%This way, underestimation will be avoided.
to avoid underestimation.

\noindent The a-posteriori error estimator (in both cases) is defined by 
\begin{equation}\label{eqn:error_estimator}
\eta = \bigg(\sum\limits_{t\in\mathcal{T}_h} \eta_t^2 \bigg) ^{1/2}.
\end{equation}
%
%Set $\eta_{\mbox{\scriptsize{max}}} = \max\limits_{t\in\mathcal{T}_h}\{
%\eta_{t} \}$
Sch\"otzau and Zhu have shown %in numerical experiments
that their error estimator is 
robust in convection-dominated regimes, 
effective in locating characteristic and boundary
layers and that the error in the energy norm converges with optimal order as
soon as refinement reaches a state when the local mesh P\'eclet number is
of order $1$. An extension to $hp$-adaptivity can be found in
\citep{Schoetzau:2011}.
The performance of Verf\"urth's error estimator is assessed in a comparative
study by \citet{John:2000}.

Alternative error estimators for DG discretization schemes of the
convection-diffusion equation can be found in 
\citep{Ern:2010} and \citep{Georgoulis:2009}.

A robust error estimator is one aspect of adaptivity. Another important aspect
is a marking strategy that achieves an equitable distribution of error
contributions. % to all cells.
An error-fraction based refinement strategy has the following characteristics:
Given a fixed refinement fraction $p_r[\%]$
and the list of all cells sorted by the magnitude of the local error indicators
$$\eta_{t_{j_1}}^2 \le \eta_{t_{j_2}}^2 \le ... \le \eta_{t_{j_n}}^2,$$ 
the goal is to mark the cells with the largest local errors for refinement such that
their contribution to the total error is $p_r[\%]$. 
%
% To be more precise, we need to find the critical index 
% \begin{equation}\label{eqn:J_error_fraction}
%   J = \min\bigg\{ j: \sum\limits_{k=j}^{n} \eta_{t_{i_k}}^2 \ge
% %  \frac{p_r}{100} \sum\limits_{k=1}^{n} \eta_{t_{i_k}}^2 
%   \frac{p_r}{100} \cdot \eta^2 
% \bigg\}
% \end{equation}
%
To be more precise, we need to find the largest $\eta^{\star}$ such that 
\begin{equation}\label{eqn:J_error_fraction}
  \sum\limits_{t\in\mathcal{T}_h: \eta_t \ge \eta^{\star}} \eta_{t}^2 \ge
  \frac{p_r}{100} \cdot \eta^2 
\end{equation}
using the bisection method and mark the top contributors
$t\in\mathcal{T}_h$ with $\eta_t \ge \eta^{\star}$
for refinement.
This strategy is readily parallelizable:
the sum on the left hand side of (\ref{eqn:J_error_fraction}) and the total
error $\eta^2$ get their contributions from all processes.
The very same strategy can be used to mark the mesh cells with the lowest
error contribution for coarsening, given a fixed coarsening fraction
$p_c[\%]$:
Find the smallest $\rho^{\star}$ such that 
\begin{equation}\label{eqn:J_error_fraction_c}
  \sum\limits_{t\in\mathcal{T}_h: \eta_t \le \rho^{\star}} \eta_{t}^2 \le
  \frac{p_c}{100} \cdot \eta^2 
\end{equation}
and mark all cells 
$t\in\mathcal{T}_h$ for which $\eta_t \le \rho^{\star}$
for coarsening.
Note that the choice of the fixed parameters $p_r$ and $p_c$ are
dependent on the error distribution and have a strong influence on the 
efficiency of the scheme.
%
% And the error distribution depends on the base level and the problem.
%

For higher order polynomials, the second order derivative in the element
residual term (\ref{eqn:ree_1}) would be required.
In the small $2$-D test problems, this is neglected.
In the solute transport simulations, we apply adaptivity only to the DG($1$)
discretization.
% Adaptive SDFEM depends either on conforming
% mesh refinement or on a special treatment of hanging nodes 
% (this increases the matrix stencil). 
% For a parallel code, this would require extra interprocess communication cost.
% Adaptive SDFEM is therefore used only for a small $2$-D test problem 
% in a sequential code.
For a comparative study based on small $2$-D test problems, a sequential
version of the adaptive SDFEM code is sufficient.

\noindent Remember that in the practical application, $\mathcal{T}_{0}$ is the mesh on
which 
\begin{itemize}
\item the conductivity field $K$ is resolved (\ref{eqn:Kh}),
\item the flow equation (\ref{eqn:GWE}) is solved as in
\S \ref{paragraph:GWE_CCFV} and 
\item the Darcy flux (\ref{eqn:Darcy}) is evaluated as in
\S \ref{paragraph:RT0}.
% and step 2 is to be replaced by:
%\begin{enumerate}
%\item[ii')] Compute the solution of the solute transport equation (\ref{eqn:TPE}) on
%  $\mathcal{T}_{h_L}$ as described in \ref{paragraph:DG}.
%\end{enumerate}
\end{itemize}
The mesh $\mathcal{T}_{0}$ should be sufficiently fine such that 
the main features of the solution are visible
and a partitioning of the mesh among all available processes is possible.
The {\bf stopping criterion} for refinement is reached as soon as either of
the following conditions is fulfilled:
\begin{enumerate}
  \item[(i)] The estimated error is below some prescribed tolerance: 
    $\eta \le TOL$.
  \item[(ii)] Provided that the range of the true solution $u$ is given by 
    $[0,\hat{u}]$, we may choose a tolerance of $p_{\text{osc}} = 5\%$ for the 
    maximal under- and overshoots by
    \begin{equation}\label{eqn:stop}
      \max\bigg\{\dfrac{|u_{\min}|}{\hat{u}},\dfrac{u_{\max}-\hat{u}}{\hat{u}}\bigg\} <
      p_{\text{osc}}.
    \end{equation}
  \item[(iii)] The mesh refinement level $L$ or % has exceeded a maximal number of steps or the 
  %\item[(iii)]  
    the total number of unknowns has exceeded a certain limit.
\end{enumerate}

The {\bf $h$-adaptive mesh refinement algorithm} for the solution of the
convection-diffusion equation can be formulated as follows.

\begin{algorithm}[H]
%  \footnotesize
  \caption{$h$-adaptive refinement}
  \label{alg:Adaptive}
  \begin{algorithmic}
    \State\textbf{Input:} Appropriate values for $p_r$ and $p_c$. 
    \State (1) Start with mesh level $L=0$.
    \State (2) Compute the solution $u_{h_0}$ of (\ref{eqn:DG_variational}) on
    $\mathcal{T}_{h_0}$.
    \State (3) Compute the error estimator $\eta$ as in
    (\ref{eqn:error_estimator}) for $u_{0}$.
    \While{ $\eta > TOL$ }
    \State (4) Apply the marking strategy (\ref{eqn:J_error_fraction}).
    \State (5) Refine the mesh and set $L=L+1$.
    \If{$L>L_{\max}$}
    \State break; \Comment{// maximal number of refinement steps exceeded}
    \EndIf
    \State (2') Compute the solution $u_{h_L}$ of (\ref{eqn:DG_variational}) on
    $\mathcal{T}_{h_L}$.
    \If{(\ref{eqn:stop}) holds}
    \State break; \Comment{// overshoots and undershoots are small enough}
    \EndIf
    \State (3') Compute the error estimator $\eta$ as in
    (\ref{eqn:error_estimator}) for $u_{h_L}$.
    \EndWhile
    \State \textbf{Output: $u_{h_L}$ }
  \end{algorithmic}
\end{algorithm}

In each refinement step, the linear system for (\ref{eqn:DG_variational}) can be solved
independently of solutions from the previous refinement step since the problem
is linear and stationary.

\section{Diffusive \texorpdfstring{$\boldsymbol{L^2}$}{L2}-projection} \label{paragraph:diffusiveL2projection}

In this section, we present another method to reduce numerical oscillations.
Due to its simplicity we consider this a post-processing step for the DG solution.
Given the DG solution $u_{\hbox{\tiny{DG}}}\in W_{h,k}$ on the coarse level 
$h=h_0$, our goal is to find an approximation of $u_{\hbox{\tiny{DG}}}$ 
in the space $V_h$ of continuous Galerkin finite
elements %\footnote{for reasons of code reuse and comparability with
%  finite element solutions},
that preserves the profile of the DG solution, but with a significant
reduction of spurious oscillations.
The $L^2$-projection is a good candidate. It is well-known to give a good on
average approximation of a function and it does not require the approximated
function to be continuous. Furthermore, an extra term imitating a small amount
of diffusive flux can be added. 
This way, the $L^2$-projection can be interpreted as the solution of a
diffusion-reaction equation without boundary constraints.
This leads to the following variational problem: 
\begin{center}
  \framebox[0.49\textwidth]{
    \begin{minipage}[c]{0.48\textwidth}
      \smallskip
      \noindent Find $u_h\in V_h$ such that 
      \begin{equation}\label{eqn:L2Projection}
        \begin{split}
          & \big( \varepsilon_h \nabla u_h, \nabla v_h \big)_{0,\Omega}
          + \big( u_h, v_h \big)_{0,\Omega}
          = \big( u_{\hbox{\tiny{DG}}}, v_h \big)_{0,\Omega}
          \\
          & \hspace{5cm} \forall ~ v_h\in V_h.
        \end{split}
      \end{equation}
    \end{minipage}
  }
\end{center}
Hereby, we choose the extra diffusion $\varepsilon_h =
\frac{1}{2}h^2$ in such a way that the diffusivity of characteristic layers
are in the order of magnitude of the meshsize $\sim O(\sqrt{\varepsilon_h}) =
O(h)$.

\section{Efficient solution of the arising linear
  systems}\label{paragraph:renumbering}

\subsection{Flow equation and diffusive \texorpdfstring{$\boldsymbol{L^2}$}{L2}-projection}
The linear systems arising from the discrete elliptic problems 
(\ref{eqn:FEM}) or (\ref{eqn:CCFV}) 
and (\ref{eqn:L2Projection})
are all of the size $O(n^2)$, symmetric positive definite
%, fulfills the M-matrix properties 
and can be solved very efficiently using the combination 
CG with AMG.
The AMG preconditioner described by \citet{Blatt:2010}
is designed for the solution of problems of the type (\ref{eqn:GWE})
with a highly discontinuous coefficient $K$.

\subsection{Transport equation}
By contrast, the stiffness matrix of the discrete transport equation is
non-symmetric. BiCGSTAB or alternatively GMRES are used in our numerical tests.
For the SDFEM discretization (\ref{eqn:SDFEM}), the matrix size is also
$O(n^2)$.
In the more diffusive case of heat transport, parallel AMG may be used as
a preconditioner.
In the convection-dominated case, the SSOR or ILU(0) preconditioners are used.

As mentioned in the introduction, for the discretization of a first order
hyperbolic problem using an upwind scheme,
the order in which the unknowns are indexed,
plays an important role for the performance and stability of an iterative
solver. The main purpose of ordering unknowns in flow direction can already be
found in \citep{Reed:1973}. The downwind numbering algorithms described in the works
of \citet{Bey:1997} and \citet{Hackbusch:1997} handle arbitrary velocity fields.
Steady-state groundwater flow (with a scalar conductivity field) is a
potential flow and therefore always cycle-free.
Since the velocity field is induced by the hydraulic head, the latter can be
used directly as the sorting key for the unknowns.

For the DG discretization, it is advisable to collect the 
unknowns of the solution vector $\vec{u}_h$ block-wise where each vector block
$\vec{u}^{(t)}$ holds the unknowns of 
$\{u_1^{(t)},...,u_{n_{\text{local}}}^{(t)}\}$ of a single mesh cell $t$ with 
$n_{\text{local}}$ denoting the dimension of the local polynomial space.
The stiffness matrix $\mathbf{A}_h$ becomes a block
matrix with constant block-size $n_{\text{local}}\times n_{\text{local}}$. The arising
linear system
\begin{equation}\label{eqn:LGS}
  \mathbf{A}_h ~ \vec{u}_h = \vec{b}_h
\end{equation}
is of the size $O(n^2\cdot n_{\text{local}}^2)$ and 
can be solved efficiently using a block version of 
BiCGSTAB or GMRES combined with SSOR or ILU(0) preconditioning, after a
renumbering of mesh cells:
In the hyperbolic limit
the bilinear form of the DG discretization is reduced to the terms listed in
the first two lines of (\ref{eqn:DG_bilinear}). 
If the mesh cells are sorted 
according to the hydraulic head distribution $\phi$ 
the stiffness matrix $\mathbf{A}_h$ obtains the shape of a block-triangular matrix.
In this case, the symmetric block Gauss-Seidel method for (\ref{eqn:LGS}) becomes a
direct solver because it converges after one step.

If the groundwater flow and the solute transport equations are
solved on the same mesh, this procedure is straightforward.
Otherwise, if adaptive refinement is applied only to the solution of the
transport problem, as mentioned in subsection \ref{paragraph:PreDefs},
the hydraulic head $\phi$ must be reconstructed on the locally refined 
sub-cells. 
Given the discrete Darcy velocity $\vec{q}_h$ in the form
(\ref{eqn:RT0_Darcy}), % on the coarse mesh $\mathcal{T}_0$, 
the hydraulic head can be 
locally reconstructed as a quadratic function 
\begin{equation}\label{eqn:potential_function}
\tilde{\phi}_{|t} = \sum\limits_{j=1}^d \big( a_j x_j^2 + b_j x_j \big) + c_0 
\end{equation}
satisfying 
\begin{equation}\label{eqn:potential_function_center}
  \tilde{\phi}_{|t} = \phi_{|t}
\end{equation}
on the cell center and 
the discrete form of Darcy's law 
\begin{equation}\label{eqn:potential_function_Darcy}
-K \nabla \tilde{\phi}_{|t} = \vec{q}_h 
\end{equation}

\begin{figure}[H]
  \centering
  \begin{tikzpicture}
    %\node (title) at (0,5) {title};
    \begin{groupplot}[
        group style={
          group name=my plots,
          group size=2 by 1,
          horizontal sep=0.15cm,
          vertical sep=0.15cm
          %ylabels at=edge right
        },
        footnotesize,
        % tickpos=left,
        ytick align=inside,
        xtick align=inside,
        enlarge x limits=false,
        enlarge y limits=false,
        x tick label style={rotate=0},
        y tick label style={rotate=90}
      ]
      %Legend
      \nextgroupplot[title={},width=0.7\linewidth,height=0.7\linewidth,xticklabels=\empty,yticklabels=\empty]
      \addplot graphics[xmin=0,xmax=50,ymin=0,ymax=50,includegraphics={keepaspectratio}]{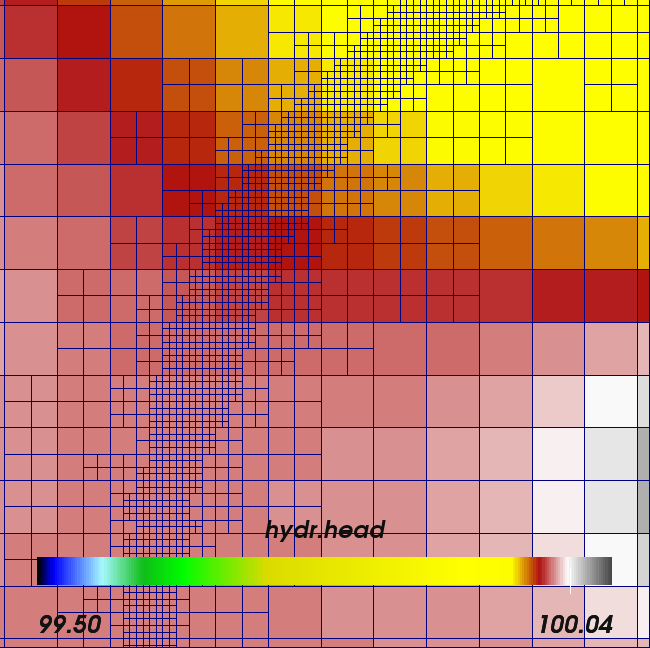};
      \node[fill=white,thick,draw=black,rounded corners] at (axis cs:30,11) {\tiny{\it hydraulic head $\phi_h$}};
      \nextgroupplot[title={},width=0.7\linewidth,height=0.7\linewidth,xticklabels=\empty,yticklabels=\empty]
      \addplot graphics[xmin=0,xmax=50,ymin=0,ymax=50,includegraphics={keepaspectratio}]{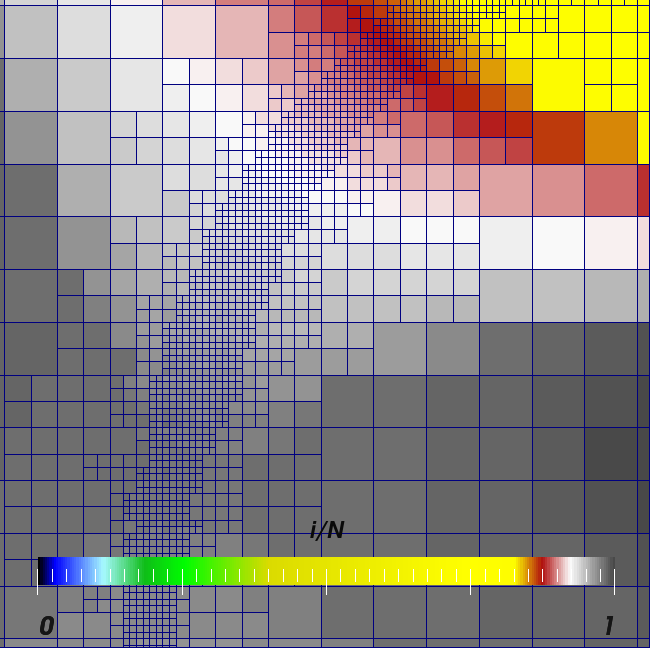};
      \node[fill=white,thick,draw=black,rounded corners] at (axis cs:25,11) {\tiny{$i/n$}};
    \end{groupplot}
  \end{tikzpicture}
  \caption[Downwind numbering of locally refined mesh cells]
          {
    Downwind numbering of locally refined mesh cells:
    The hydraulic head $\phi$ is constant on each coarse mesh cell
    (upper plot). For the renumbering of mesh cells according to the
    hydraulic head, it has to be resolved on the locally refined sub-cells.
    The cell index $i$ in the lower plot is sorted according
    to $\tilde{\phi}$. 
    $n$ is the total number of all cells of the locally refined mesh.
  }
  \label{fig:headcolor}
\end{figure}

on the $2d$ face centers 
of a coarse mesh cell $t\in\mathcal{T}_0$.
This yields $2d+1$
equations for the $2d+1$ coefficients of $\tilde{\phi}_{|t}$.
The locally refined sub-cells can then be sorted according to 
$\tilde{\phi}_{|t}$ evaluated at the centers of the sub-cells (Figure \ref{fig:headcolor}).

Due to the large problem size in $3$-D, parallelization 
is necessary for efficiency.
After each refinement step, the parallel partition-blocks\footnote{ In general, the shape of these blocks are not cuboidal.} 
of the coarse mesh $\mathcal{T}_0$ (and with it all locally refined sub-cells) 
may be altered to achieve a similar amount of refined sub-cells on every processor partition
(\textbf{dynamic load-balancing}).
To ensure that the renumbering procedure still works after each repartitioning
step, the two quantities $\phi$ and $K$ have to 

\noindent be made available on each
processor partition of $\mathcal{T}_0$ for a new reconstruction of
$\vec{q}_h$ and $\tilde{\phi}$.

\section{Numerical Studies}\label{paragraph:numerical_tests}

In \S \ref{s:John} and \S \ref{paragraph:Lopez}, we start with two singularly perturbed problems
on the unit square for which analytical solutions exist in the domain of
interest. Convergence tests can be performed using global
and adaptive refinement.
For the first problem, we demonstrate the influence of the
ordering of unknowns on the performance of the iterative solvers 
for linear systems arising from a DG($1$) discretization 
as described in section \ref{paragraph:renumbering}.
The second problem has a less regular solution. 
In \S \ref{paragragph:test3}, we take a closer look on the quality of the diffusive
$L^2$-projection compared to the SDFEM solution on a structured mesh.
In \S \ref{example:forward2D} and \S \ref{example:forward3D}, we show a $2$-D and a $3$-D example of the coupled
groundwater flow and transport problem. Since analytical solutions are not
available, numerical solutions computed on adaptively refined meshes are taken
as reference solutions for assessing the quality of different solution methods
computed on a coarse structured mesh.

\medskip

For computations on structured meshes in $2$-D and in $3$-D, we use the 
{\tt YASP} grid, an implementation of a structured parallel mesh available in
the {\tt dune-grid} module.
For sequential adaptive refinement in $2$-D, we use the 
{\tt UG} grid %\footnote{\scriptsize \url{http://www.iwr.uni-heidelberg.de/frame/iwrwikiequipment/software/ug}}
\citep{ug}, and for parallel adaptive refinement with dynamic load-balancing in $3$-D,
we use the
{\tt DUNE ALUGrid} module %\footnote{\url{http://aam.mathematik.uni-freiburg.de/IAM/Research/alugrid}}
\citep{DednerKN14}.
\medskip

In all computations, the best available linear solver / preconditioner
combination (in terms of robustness and speed)
is chosen for each linear system arising from a finite element discretization
of a stationary problem.

\medskip
All time measurements are based on the \textbf{wall-clock time}, 
i.e. the difference between the time at which a certain task finishes and the start time of that task.
It may include time that passes while waiting for resources to become available.

\medskip

All $2$-D computations are performed in sequential mode on a laptop with 
an {\it Intel\textregistered Core\texttrademark 2  Duo CPU (P9500, 2.53 GHz)}
and $4$ GB total memory.
All $3$-D computations are performed on a multi-core architectures % (a single-node machine or a cluster with several nodes)
with large memory and high-speed network communication links. % (InfiniBand).
Table \ref{t:IWRmachines} give an overview of the used hardware.

\subsection{An example with a regular solution}\label{s:John}
Let $(x,y)\in \overline{\Omega} = [0,1]^2$ and
consider the boundary value problem (from \citep{John:1997})
\begin{eqnarray}
  - \varepsilon \Delta u + \vec{q} \cdot \nabla u + \mu \cdot u &=& \tilde{s}_{\varepsilon} ~
  \hbox{~in~} \Omega\\
  u &=& 0 ~ \hbox{~on~} \partial \Omega
\end{eqnarray}
where $\vec{q}=( 2, 3 )^T$, $\mu=2$ and the source term $\tilde{s}_{\varepsilon}(x,y)$ is chosen such that 
\begin{equation}
\footnotesize
  u(x,y) = \frac{16}{\pi}x(1-x)y(1-y)\bigg(\frac{\pi}{2} + \arctan{ \bigg[
      \frac{2}{\sqrt{\varepsilon}} \xi(x,y) \bigg]} \bigg)
\end{equation}
with
\begin{equation}
  \xi(x,y) = 0.25^2-(x-0.5)^2-(y-0.5)^2.
\end{equation}

\noindent For our tests, we choose $\varepsilon=10^{-5}$. 
Note that in this example, the internal layer is generated by a source term
which itself depends on $\varepsilon$. For $\varepsilon\ll 1$,
an accurate representation of the source term requires a fine mesh,
because, in a finite element discretization of $\tilde{s}_{\varepsilon}$, the
error in the quadrature-rule might become dominating on a coarse mesh.
Since we investigate the convergence behavior for global and adaptive
refinement, this is not a severe problem. % for DG-based methods compared to the SDFEM method.

\begin{figure}[H]
\vspace{-2mm}
  \begin{center}
    \includegraphics[width=0.7\linewidth]{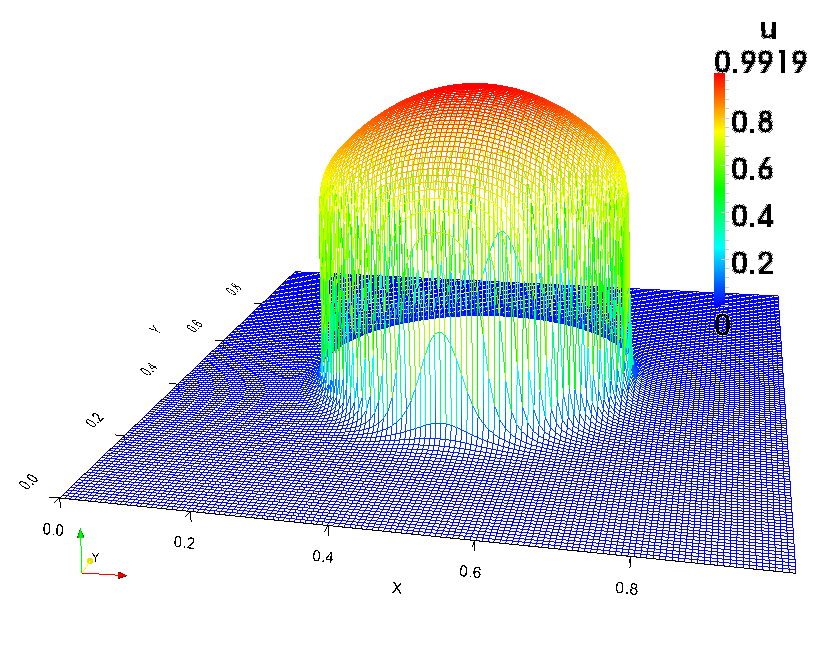}
    \vspace{-5mm}
    \caption[$2$-D example by John: Analytical solution]{Profile of the analytical solution $u$ for $\varepsilon=10^{-5}$.}
    \label{fig:Ex1_True_Solution}
  \end{center}
\end{figure}

\subsubsection*{Linear solver performance and accuracy}

\begin{table}[H]
  %\footnotesize
  \begin{center}
    \hspace{-5mm}
    \scalebox{0.69}{
      \ra{1.3}
      \begin{tabular}{@{}r|r|rlrlrlr@{}}
        \toprule
        &      & \multicolumn{2}{c}{random} &
        \multicolumn{2}{c}{horizontal} & \multicolumn{3}{c}{\color{blue}downstream} \\
        ~~$L$ &     $DOF$  &  $IT$ & $TIT[s]$ &  $IT$ & $TIT[s]$ &  $IT$ & $TIT[s]$ & $T_{\mbox{\tiny{sort}}}[s]$  \\
        \midrule
        ~~9  &    6,856  &   8 & 0.005  &    4 & 0.004  & \color{blue} 1  & \color{blue} 0.008  & 0.003 \\
        ~~10 &   10,528  &   9 & 0.008  &    5 & 0.007  & \color{blue} 1  & \color{blue} 0.008  & 0.006 \\
        ~~11 &   19,672  &  13 & 0.015  &    6 & 0.014  & \color{blue} 2  & \color{blue} 0.010  & 0.012 \\
        ~~12 &   34,540  &  17 & 0.028  &    7 & 0.026  & \color{blue} 2  & \color{blue} 0.028  & 0.022 \\
        13 &   85,468  &  25 & 0.081  &    9 & 0.079  & \color{blue} 3  & \color{blue} 0.068    & 0.059 \\
        14 &  150,016  &  34 & 0.152  &   10 & 0.138  & \color{blue} 4  & \color{blue} 0.127    & 0.123 \\
        15 &  278,836  &  48 & 0.325  &   12 & 0.262  & \color{blue} 5  & \color{blue} 0.246    & 0.284 \\
        16 &  544,300  &  73 & 0.770  &   14 & 0.516  & \color{blue} 7  & \color{blue} 0.496    & 0.661 \\
        \bottomrule
      \end{tabular}
    }
    \caption[$2$-D example by John: Comparing cell numbering strategies]
    { Performance of the linear solver (BiCGSTAB + SSOR with reduction $10^{-8}$) for different
      cell numbering strategies applied to the DG($1$) method:
      $L=$ refinement level,
      $DOF=$ degrees of freedom,
      $IT=$ number of iterations for the linear solver,
      $TIT=$ time per iteration,
      $T_{\mbox{\tiny{sort}}}=$ time for renumbering the grid cells.
    }
    \label{t:renumberDOF}
  \end{center}
\end{table}

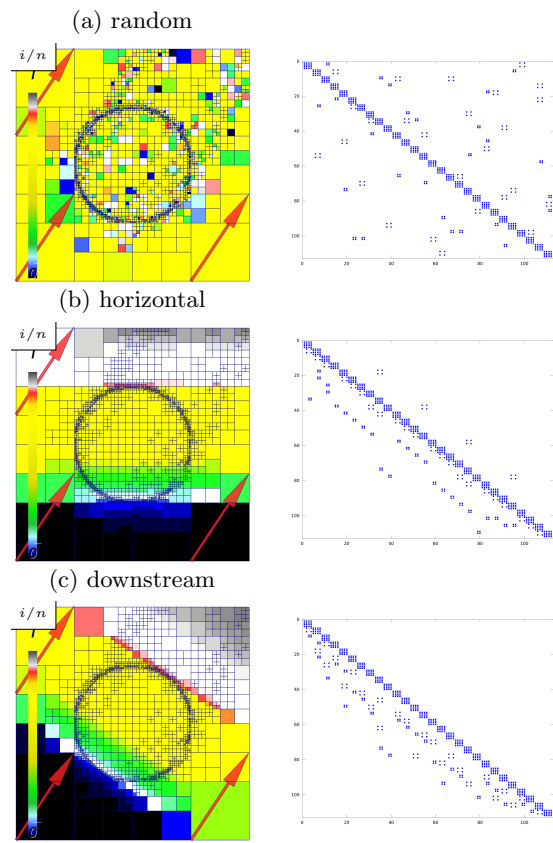
\begin{figure}[H]
  \centering
  \begin{tikzpicture}
    \begin{groupplot}[
        group style={
          group name=my plots,
          group size=2 by 3,
          horizontal sep=0.1cm,
          vertical sep=0.5cm
          %ylabels at=edge right
        },
        footnotesize,
        % tickpos=left,
        ytick align=inside,
        xtick align=inside,
        enlarge x limits=false,
        enlarge y limits=false,
        x tick label style={rotate=0},
        y tick label style={rotate=90}
      ]
      %Legend
      \nextgroupplot[title={(a) random},width=0.65\linewidth,height=0.65\linewidth,xticklabels=\empty,yticklabels=\empty,hide axis]
      \addplot graphics[xmin=0,xmax=50,ymin=0,ymax=50,includegraphics={keepaspectratio}]{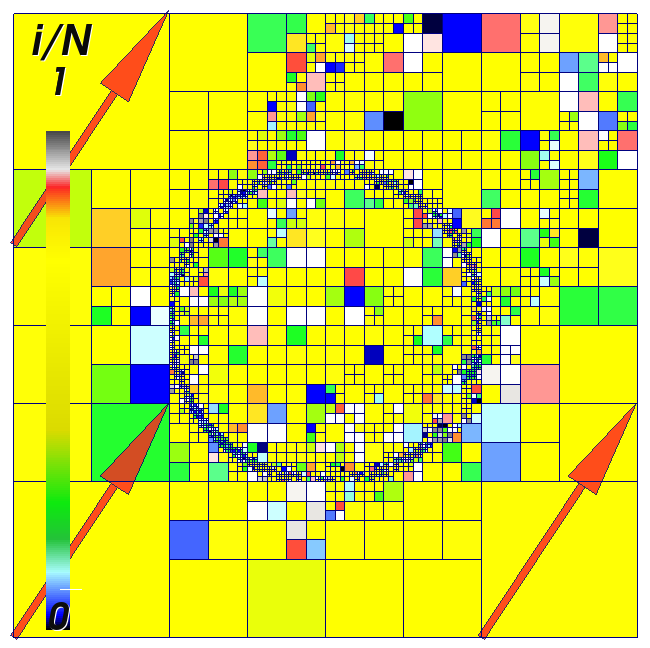};
      \node[fill=white,draw=black] at (axis cs:4.5,47.5) {\tiny{$i/n$}};
      \nextgroupplot[title={},width=0.80\linewidth,height=0.65\linewidth,xticklabels=\empty,yticklabels=\empty,hide axis]
      \addplot graphics[xmin=0,xmax=55,ymin=0,ymax=45,includegraphics={keepaspectratio}]{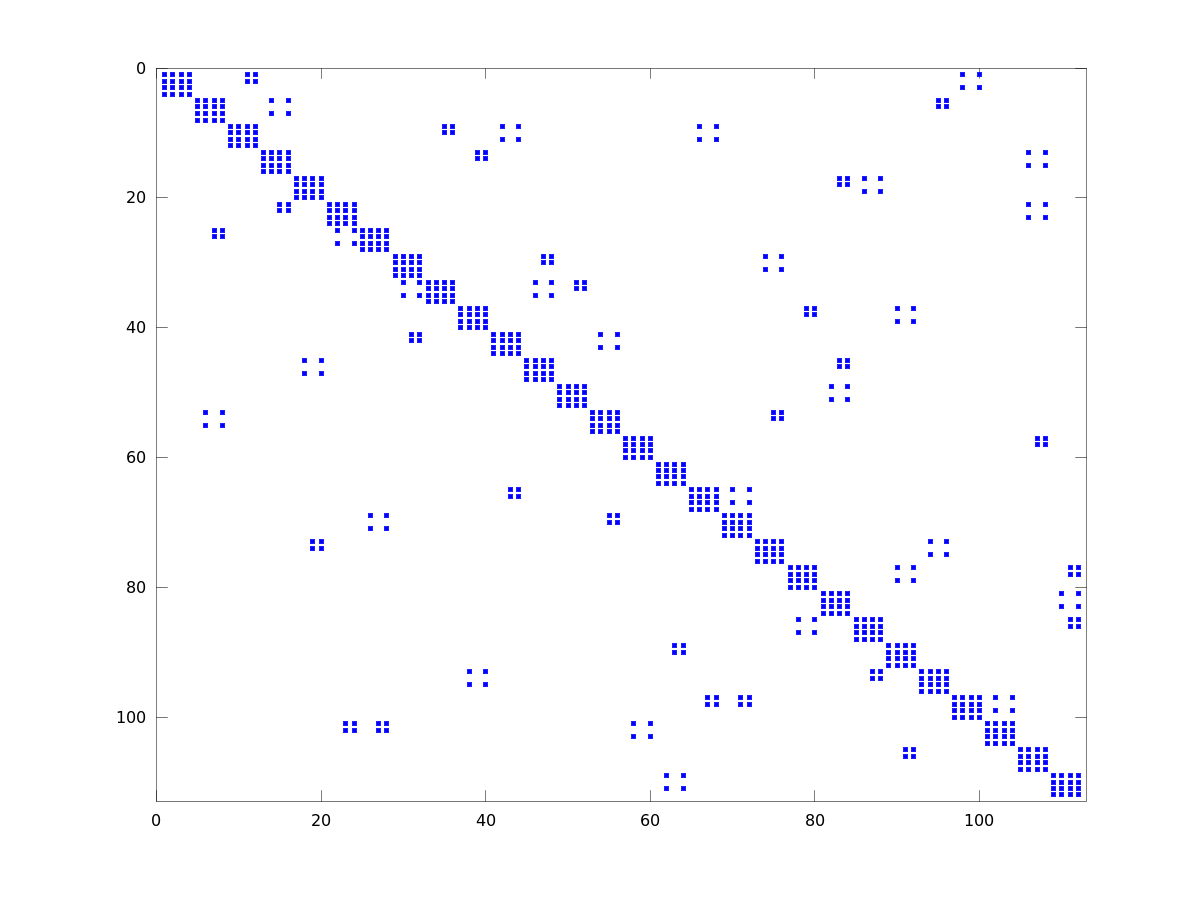};
      \nextgroupplot[title={(b) horizontal},width=0.65\linewidth,height=0.65\linewidth,xticklabels=\empty,yticklabels=\empty,hide axis]
      \addplot graphics[xmin=0,xmax=50,ymin=0,ymax=50,includegraphics={keepaspectratio}]{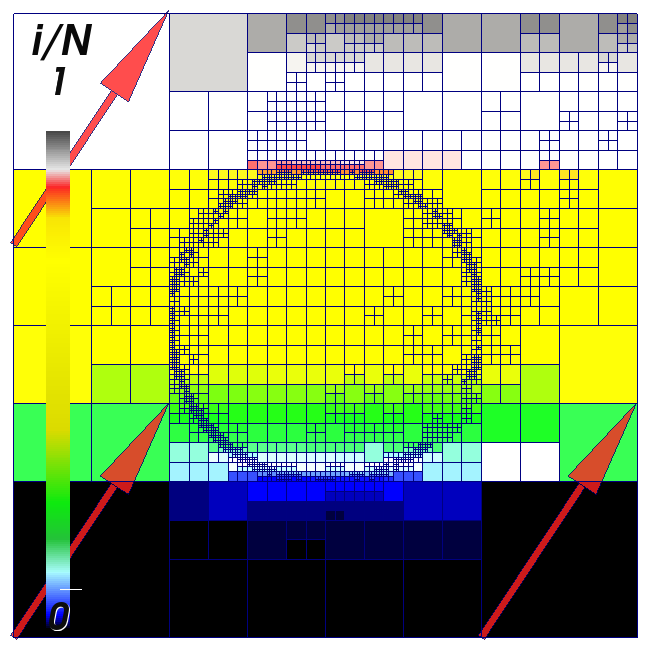};
      \node[fill=white,draw=black] at (axis cs:4.5,47.5) {\tiny{$i/n$}};
      \nextgroupplot[title={},width=0.80\linewidth,height=0.65\linewidth,xticklabels=\empty,yticklabels=\empty,hide axis]
      \addplot graphics[xmin=0,xmax=55,ymin=0,ymax=45,includegraphics={keepaspectratio}]{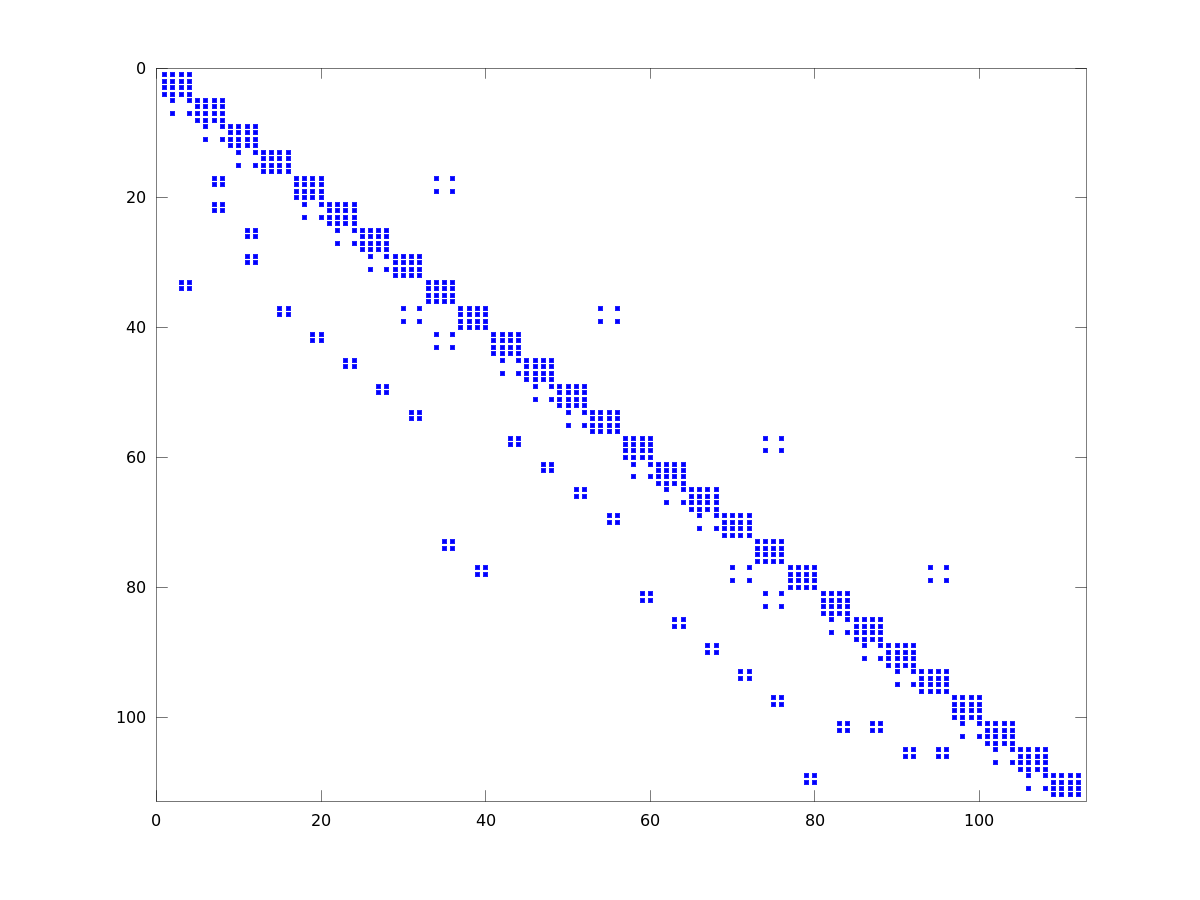};
      \nextgroupplot[title={(c) downstream},width=0.65\linewidth,height=0.65\linewidth,xticklabels=\empty,yticklabels=\empty,hide axis]
      \addplot graphics[xmin=0,xmax=50,ymin=0,ymax=50,includegraphics={keepaspectratio}]{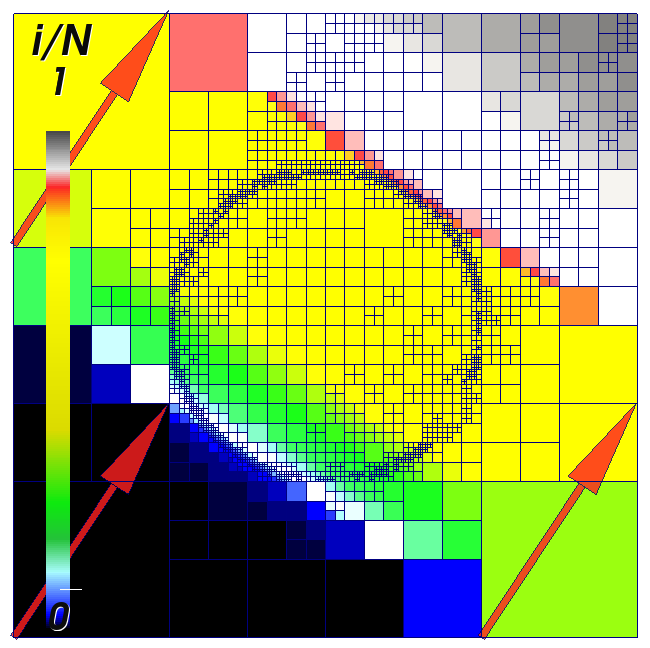};
      \node[fill=white,draw=black] at (axis cs:4.5,47.5) {\tiny{$i/n$}};
      \nextgroupplot[title={},width=0.80\linewidth,height=0.65\linewidth,xticklabels=\empty,yticklabels=\empty,hide axis]
      \addplot graphics[xmin=0,xmax=55,ymin=0,ymax=45,includegraphics={keepaspectratio}]{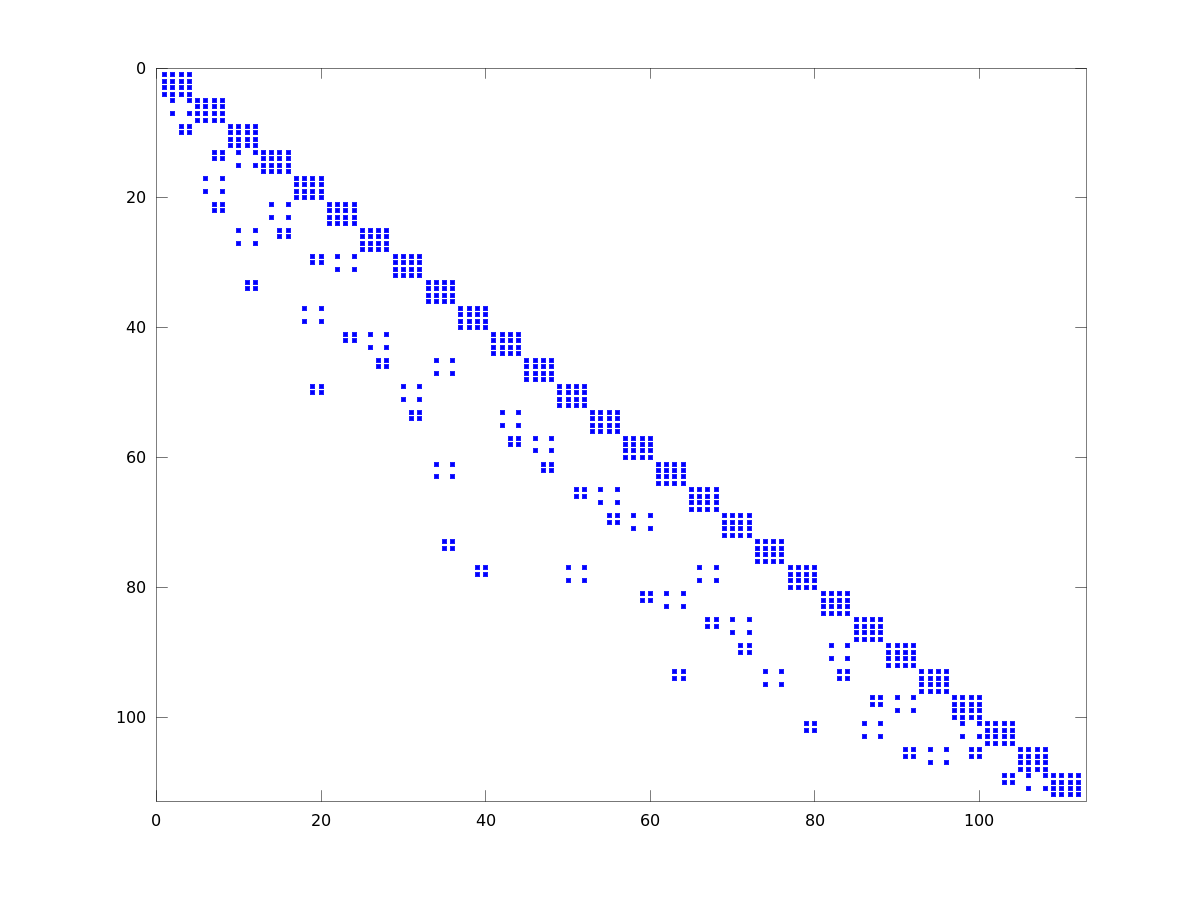};
%      \draw node[anchor=south west] (nodoFondo) at (axis cs:0,0) {\pgfuseimage{fondo}};
    \end{groupplot}
  \end{tikzpicture}
  \caption[$2$-D example by John: Comparing cell numbering strategies]
          {Different cell numbering strategies 
    %((a) random, (b) horizontal,(c) downstream) 
    and corresponding matrix patterns for the DG($1$) method.
    The block matrix for $\mathbb{Q}_1$ elements in $2$-D is made of 
    $4\times 4$ - blocks. % with the $\mathbb{Q}_1$-polynomial basis. 
    $i=$ cell index, $n=$ total number of mesh cells.
  }
  \label{fig:renumbering}
\end{figure}

Starting on a coarse structured mesh with $h_0=1/8$, we perform an adaptive
refinement loop three times, based on three different ways of cell numbering as
depicted in Figures \ref{fig:renumbering}(a)-(c).
For the SDFEM discretization, the different cell numbering has no influence on
the speed of the linear solver.
For the DG($1$) discretization, Table \ref{t:renumberDOF} confirms that an
optimal numbering of degrees of freedom (following the velocity field) results
in a faster solution of the arising linear system.
The time for renumbering the grid cells is comparable to the time for one step of the iterative linear solver.

\noindent On coarse meshes, the matrix pattern can assume a block-triangular form 
(Figure \ref{fig:renumbering}(c)). In these cases, the iteration number is
indeed $1$.
As refinement proceeds, the meshsizes and therefore
the mesh P\'eclet numbers decrease and we diverge from the hyperbolic
limit. Although this increases the iteration numbers, they stay at a low level for
the optimal numbering. Not only the number of iterations is reduced but also
the required time per iteration.
The linear solver used is BiCGSTAB with SSOR.
Similar results are obtained with the combinations BiCGSTAB + ILU(0),
GMRES + SSOR and GMRES + ILU(0).

\begin{figure}[H]
  \centering
  \begin{tikzpicture}
    \begin{groupplot}[
        group style={
          group name=my plots,
          group size=1 by 2,
          horizontal sep=0cm,
          vertical sep=0.8cm
          %ylabels at=edge right
        },
        footnotesize,
        % tickpos=left,
        ytick align=inside,
        xtick align=inside,
        enlarge x limits=false,
        enlarge y limits=false,
        x tick label style={rotate=0},
        y tick label style={rotate=90}
      ]
      %Legend
      \nextgroupplot[title={(a) global refinement},width=1.2\linewidth,height=0.9\linewidth,xticklabels=\empty,yticklabels=\empty,hide axis]
      \addplot graphics[xmin=0,xmax=50,ymin=0,ymax=50,includegraphics={keepaspectratio}]{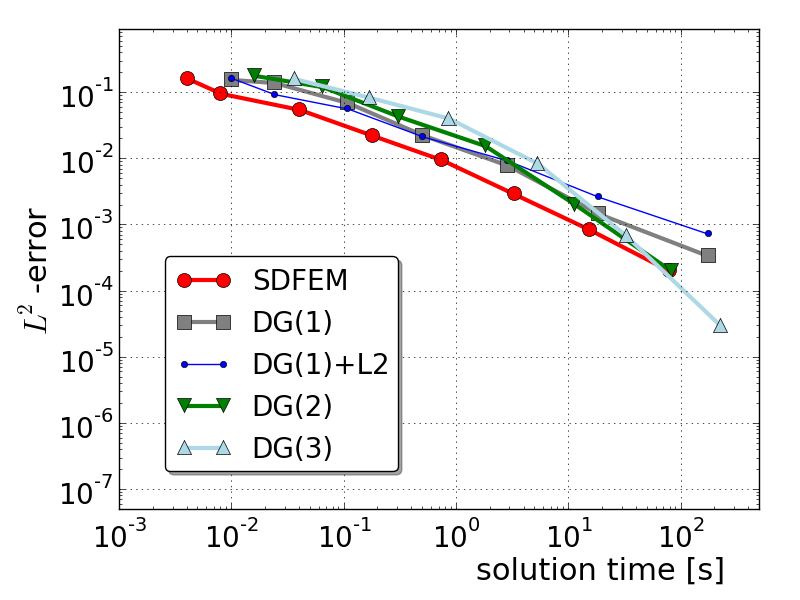};
      \nextgroupplot[title={(b) adaptive refinement},width=1.2\linewidth,height=0.9\linewidth,xticklabels=\empty,yticklabels=\empty,hide axis]
      \addplot graphics[xmin=0,xmax=50,ymin=0,ymax=50,includegraphics={keepaspectratio}]{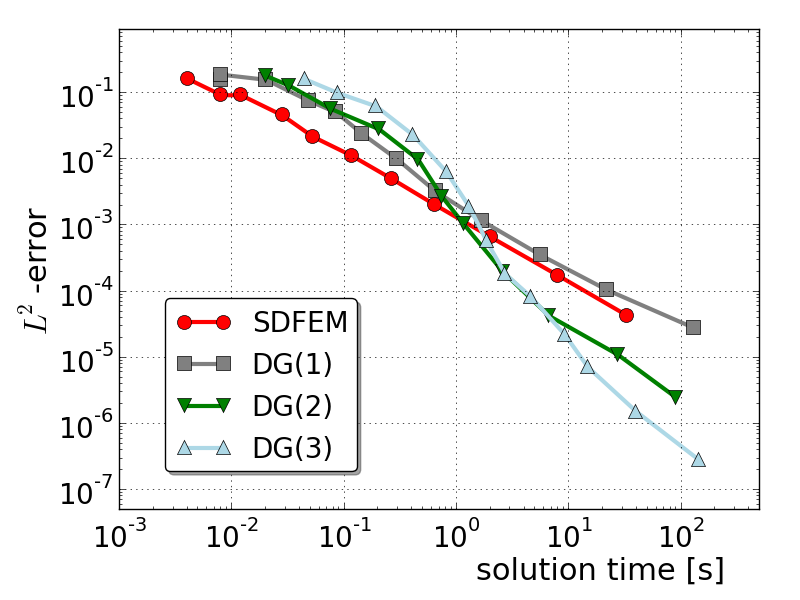};
    \end{groupplot}
  \end{tikzpicture}
  \caption{\footnotesize Example by John:
    (a) global refinement on a structured
    mesh and 
    (b) adaptive refinement on an unstructured mesh (UG),
    both starting on a coarse mesh with meshsize $h=1/8$.    
    The refinement and coarsening fractions for the adaptive refinement
    algorithm are $p_r=95[\%]$ and $p_c=0.5[\%]$.
    Linear solver used: BiCGSTAB + SSOR with reduction $10^{-8}$.
    Solution time = system assembly time + linear solver time.
  }
  \label{fig:john_convergence}
\end{figure}

We measure accuracy with respect to computing time.
The SDFEM method with bi\-li\-near elements is compared to the 
DG($k$) methods (with globally constant polynomial degree $k\in\{1,2,3\}$) in a
convergence test with uniform and adaptive refinement.
Accuracy is measured in the $L^2$-norm of the error taken
over the domain of interest: $\|u-u_h\|_{L^2(\Omega)}$.
The solution time is the sum of the system assembly time and the linear solver
time.
For the DG($k$) methods, we apply an optimal numbering of mesh cells to speed
up the linear solver. The time required to sort the mesh cells is negligible
compared to the solution time.

\noindent The solution time is linearly proportional to the number of unknowns.
For a comparable number of degrees of freedom (DOF), the iterative linear solvers
perform better for the DG-based methods.

From Figure \ref{fig:john_convergence} we can make the following observations:

\begin{enumerate}
\item For the same computing time, SDFEM is more accurate than DG($1$).
\item The accuracy of the higher order DG methods overtakes
  the accuracy of SDFEM at a certain refinement level. This happens as
  soon as the steep gradient is resolved and optimal convergence order is
  achieved.
\item With adaptive refinement, higher accuracy is achieved
  at an earlier stage.
\end{enumerate}
The blue curve (DG($1$)+L2) in the first plot displays the accuracy of the
post-processed DG($1$) solution on a structured mesh.
Although it is close to the DG($1$) curve, in this case, 
the diffusive $L^2$-projection adds an extra amount of error.
Since solution time for the post-processing step is a small fraction of the
solution time for the transport problem, it is neglected in this plot.

%\subsection{Example by L\'opez and Sinus\'ia}
\subsection{An example with a less regular solution}\label{paragraph:Lopez}
Let $(x,y)\in [0,1]^2$.
We consider the boundary value problem 
\begin{equation}
  -\varepsilon\Delta u + \vec{q} ~\nabla u = 0 \quad \text{~in~}[0,1]^2
\end{equation}
with constant velocity $\vec{q}=\frac{\sqrt{2}}{2}(1,1)^T$ and discontinuous boundary conditions
\begin{equation}
  u(x,0) = 1 \qquad \text{and} \qquad u(0,y) = 0.
\end{equation}
Obviously, the solution has a jump at the origin and is therefore not
$H^1$-regular. This jump causes a characteristic boundary layer along the
direction $\vec{q}$. This example is close to a real-world example in the
sense that the discontinuity may be used to describe a binary state and the
direction of the velocity is not aligned to the mesh.
% This is a well-suited example for
% the analysis of numerical schemes developed for the solution of convection
% dominated transport problems in a heterogeneous flow field.
Theorem 1 of \citep{Lopez:2004} provides an asymptotic expansion of the
solution $u$ on the subset $\Omega = [0,1]^2 \setminus \mathbf{U}_{0}$
where $\mathbf{U}_{0}=\{\vec{y}\in\mathbb{R}^2:\|\vec{y}\|_2 < r_0\}$ is a ball
with radius $r_0>0$ and center $(0,0)$. Introducing polar coordinates through
$x=r\sin\varphi$ and $y=r\cos\varphi$, we get
\begin{equation}
  u = u_{0}(r,\varphi) + \dfrac{e^{wr(\sin(\varphi+\beta)-1)}}{\pi\sqrt{2wr}}u_{1}(r,\varphi)
\end{equation}
where $\beta=\pi/4$, $w=\|\vec{q}\|_2/(2\varepsilon)$ and
%\twocolumn[
\begin{equation}
\footnotesize
  u_0(r,\varphi) = \frac{1}{2} \left\{
  \begin{array}{ll}
      \erfc\bigg( \sqrt{(1-\sin(\varphi+\beta))wr} \bigg) \\
      \quad \text{if~} \varphi <  \beta,  \\
    \\
      1   \\
      \quad \text{if~} \varphi = \beta, \\
    \\
      2 - \erfc\bigg( \sqrt{(1-\sin(\varphi+\beta))wr}  \bigg)   \\
      \quad \text{if~} \varphi >  \beta .
  \end{array}
  \right.
\end{equation}
The function $u_1(r,\varphi)$ has an asymptotic expansion from which we use only
the first term, 
\begin{equation}
\footnotesize
\begin{split}
  u_1(r,\varphi) = \left\{
  \begin{array}{ll}
      \Gamma(1/2) \bigg\{ \bigg(
      \frac{\cos(\varphi-\beta)}{\cos(\varphi+\beta)} - 
      \frac{\cos(\varphi+\beta)}{\cos(\varphi-\beta)} \bigg) 
      \\
      -\dfrac{1}{2\sin(\frac{1}{2}(\frac{\pi}{2}-\varphi-\beta))} \bigg\} 
      \\
      \hspace{3cm} \text{~if~} \varphi \ne \pi/4,
      \\
      \\
      0 \\ 
      \hspace{3cm} \text{~if~} \varphi = \pi/4,
  \end{array}
  \right.
      \\
      \\
\end{split}
\end{equation}
hereby neglecting higher order terms of $\varepsilon$.
%] % end of twocolumn

\noindent For our tests, we choose $\varepsilon=10^{-5}$ and $r_0=5\times
10^{-5}$.
% The $L^2$-norm of the error $u-u_h$ is taken over the 
% domain of interest $\Omega$.
The small area $\mathbf{U}_{0}(r_0)$ around the critical location $(0,0)$,
where the numerical errors are largest and a different asymptotic
expansion is necessary, is left out of consideration.
% This is legitimate for this simple example for which
%we want to observe the convergence behavior of the numerical methods.

\begin{figure}[H]
  \begin{center}
    \includegraphics[width=0.9\linewidth]{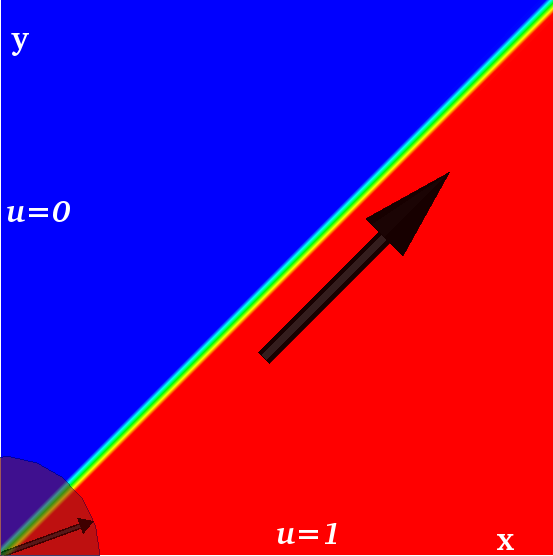}
    \caption[$2$-D example by L\'opez and Sinus\'ia: Reference solution]{Reference solution for $\varepsilon=10^{-5}$ on the domain of interest
      $\Omega = [0,1]^2 \setminus \mathbf{U}_{0}$ ($r_0=5\times 10^{-5}$).}
    \label{fig:Ex2_True_Solution}
  \end{center}
\end{figure}

%\subsubsection{Linear solver performance after renumbering}% (example 2)
\subsubsection*{Linear solver performance and accuracy}
The same numerical experiments as in \S \ref{s:John} are conducted on this example.
The influence of renumbering cells on the linear solver performance for the
DG($1$) discretization are very similar to the results presented in 
\S \ref{s:John} Table \ref{t:renumberDOF}.
%\S \ref{paragraph:LSperformance_John}.

\begin{figure}[H]
  \centering
  \begin{tikzpicture}
    \begin{groupplot}[
        group style={
          group name=my plots,
          group size=1 by 2,
          horizontal sep=0cm,
          vertical sep=0.8cm
          %ylabels at=edge right
        },
        footnotesize,
        % tickpos=left,
        ytick align=inside,
        xtick align=inside,
        enlarge x limits=false,
        enlarge y limits=false,
        x tick label style={rotate=0},
        y tick label style={rotate=90}
      ]
      %Legend
      \nextgroupplot[title={(a) global refinement},width=1.2\linewidth,height=0.9\linewidth,xticklabels=\empty,yticklabels=\empty,hide axis]
      \addplot graphics[xmin=0,xmax=50,ymin=0,ymax=50,includegraphics={keepaspectratio}]{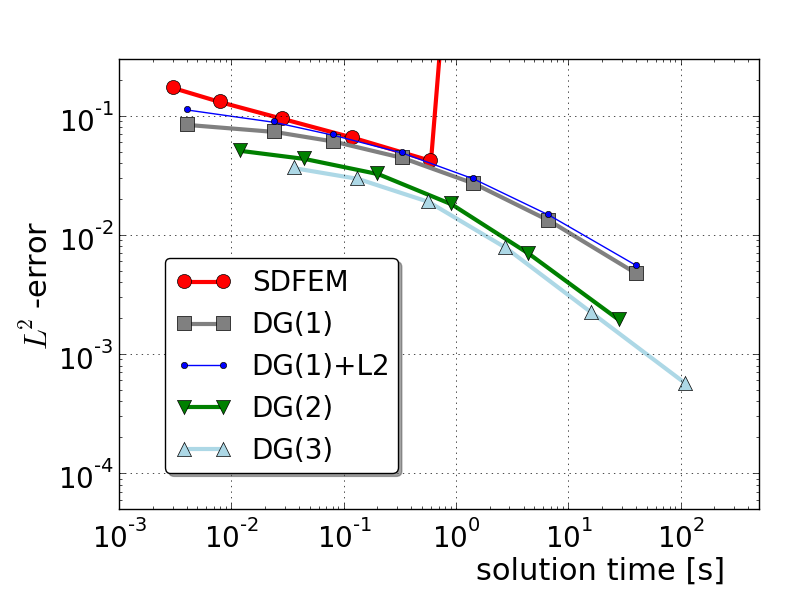};
      \nextgroupplot[title={(b) adaptive refinement},width=1.2\linewidth,height=0.9\linewidth,xticklabels=\empty,yticklabels=\empty,hide axis]
      \addplot graphics[xmin=0,xmax=50,ymin=0,ymax=50,includegraphics={keepaspectratio}]{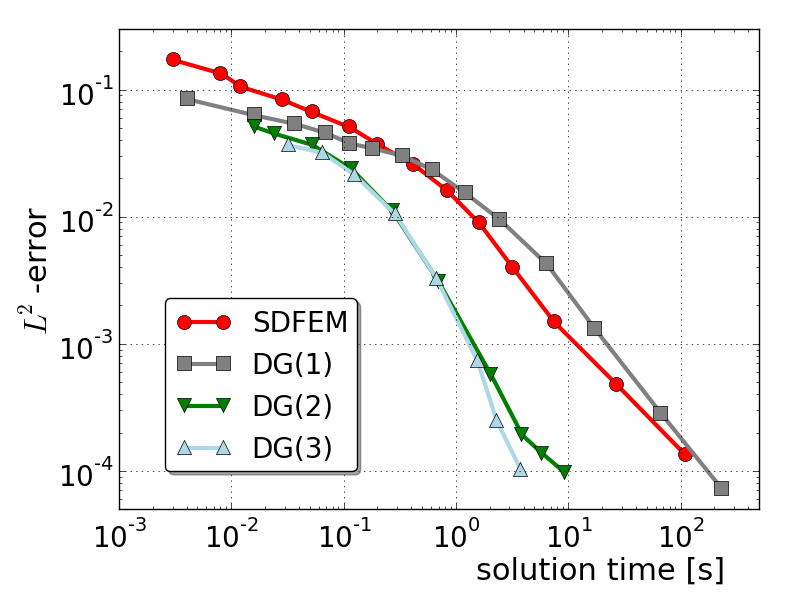};
    \end{groupplot}
  \end{tikzpicture}
  \caption[$2$-D example by L\'opez and Sinus\'ia: Numerical convergence behavior]{Example by L\'opez and Sinus\'ia:
    (a) global refinement on a structured
    mesh and 
    (b) adaptive refinement on an unstructured mesh (UG),
    both starting on a coarse mesh with meshsize $h_0=1/8$.
    The refinement fraction for adaptive refinement is $p_r=95[\%]$ (no coarsening).
    The solution time = matrix assembly time + linear solver time.
    Linear solver used: BiCGSTAB with SSOR with a reduction of $10^{-8}$.
  }
  \label{fig:lopez_convergence}
\end{figure}

From Figure \ref{fig:lopez_convergence} we can see
that the convergence behavior
for the approximation of this less regular solution
is different from the observations made in \S \ref{s:John}:
\begin{enumerate}
\item For the same solution time, 
  the accuracy of DG($1$) and SDFEM are comparable.
\item  Higher order DG methods are more accurate than SDFEM right
  from the beginning.
\item With adaptive refinement, higher accuracy is achieved
  at an earlier stage.
\end{enumerate}

%\subsubsection{Accuracy and computing time}

The blue curve (DG($1$)+L2) in the first plot is closer to the DG($1$) curve
than in the example of \S \ref{s:John} (Figure \ref{fig:john_convergence}).
Furthermore, the SDFEM plot (red curve in Figure
\ref{fig:lopez_convergence}) stops after $4$ refinement steps ($16,641$ unknowns).
In the next refinement step ($66,049$ unknowns), the iterative linear solver
BiCGSTAB with SSOR converges, but the solution is wrong. In the $6$-th step 
($263,169$ unknowns), the iterative linear solver does not converge, although
the linear system
can still be solved using the direct solver {\tt SuperLU}. This is most likely due
to the fact that the large sparse system has become very ill-conditioned.
However, direct solvers are not an option for large practical problems.
{\tt SuperLU} has reached the memory limit of the laptop already in the next
refinement level where the number of unknowns is $1,050,625$.

\subsection{Post-processed DG(1) versus SDFEM}\label{paragragph:test3}
\begin{figure}[H]
  \centering
  \begin{tikzpicture}
  \begin{axis}[width=\linewidth,height=\linewidth,hide axis,clip=false]
      \addplot+[only marks,mark=.,mark size=1pt,draw=white,fill=white] coordinates {
      (0.0,0.0)
      (0.0,100.0)
      (100.0,0.0)
      (100.0,100.0)
      };
      % links oben:
      \node[draw=black] at (axis cs:25,75) {\includegraphics[width=0.33\linewidth]{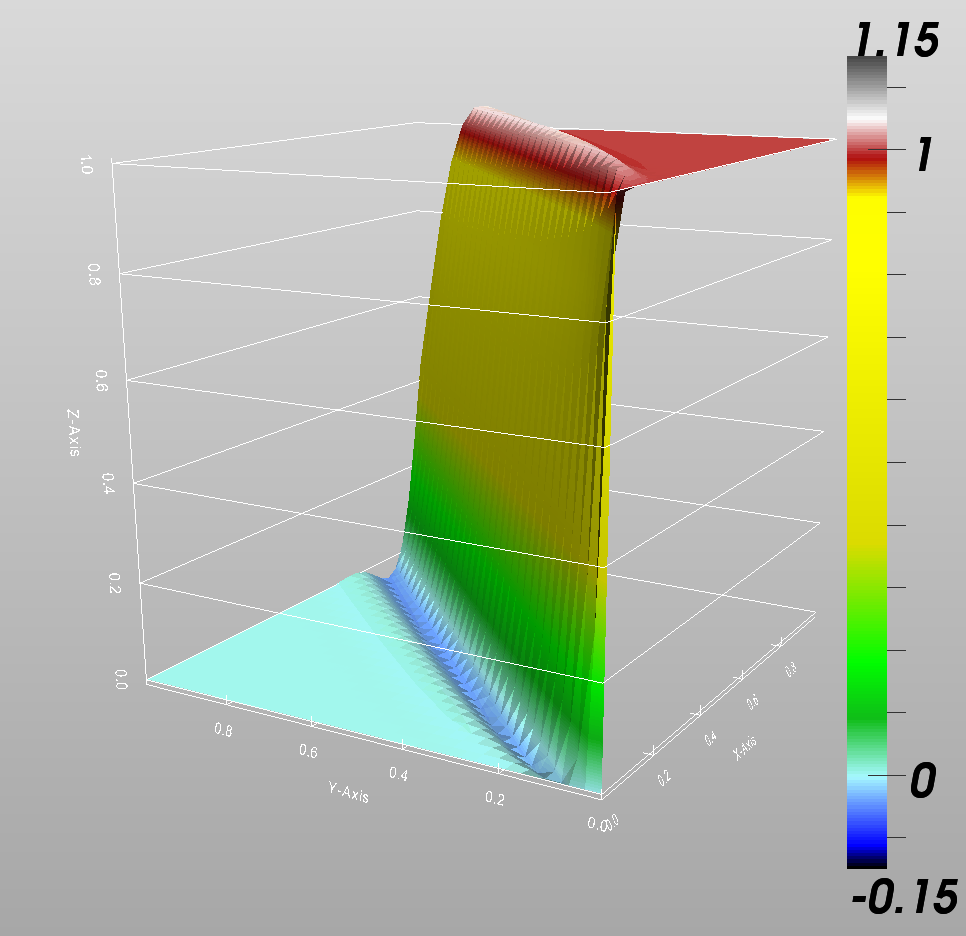}};
      \node at (axis cs:-8,95){(a)};
      \node at (axis cs:-8,90){\footnotesize $u_h = u_{sd}$};
      \node at (axis cs:-8,85){\footnotesize $h=1/32$};
      \node at (axis cs:-8,80){\footnotesize $N=1089$};
      \node at (axis cs:-8,75){\footnotesize $\epsilon = 0.096$};
      \node at (axis cs:-8,70){\footnotesize $u_{\max} = -1.045$};
      \node at (axis cs:-8,65){\footnotesize $u_{\min} = -0.047$};
      \node at (axis cs:-8,60){\footnotesize $T_{\text{sol}} = 0.028s$};
      % rechts oben:
      \node[draw=black] at (axis cs:75,75) {\includegraphics[width=0.33\linewidth]{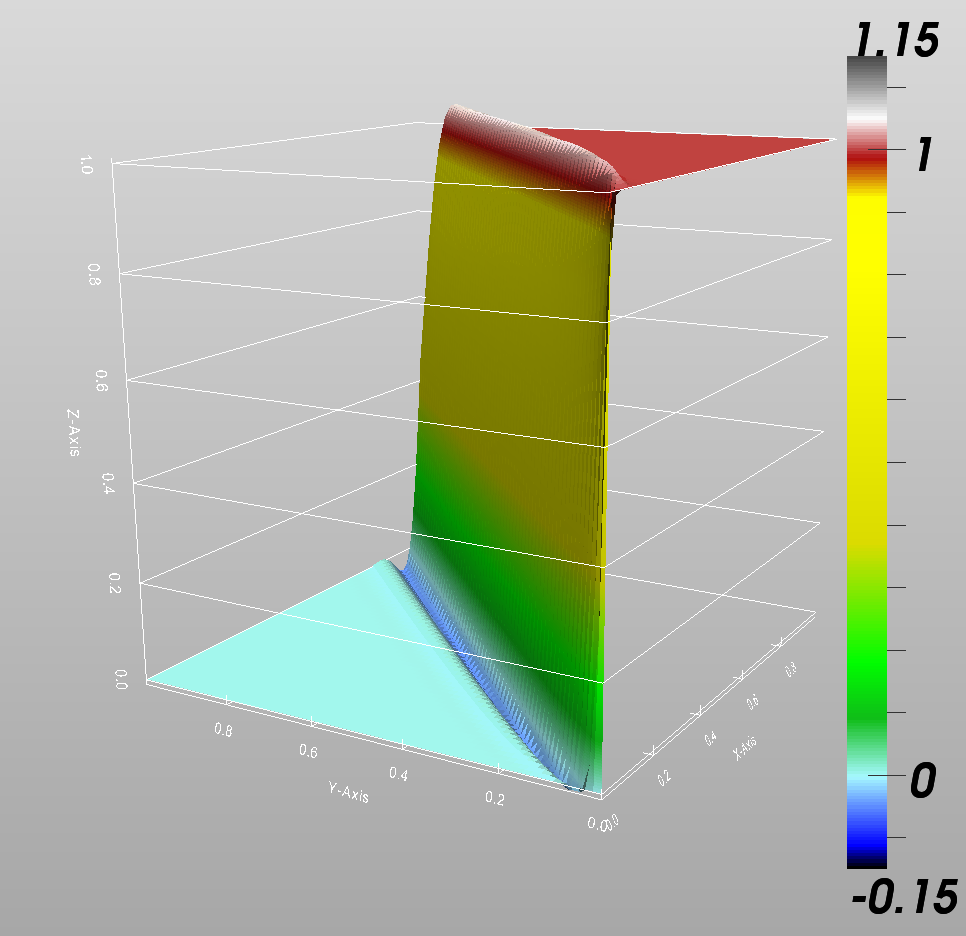}};
      \node at (axis cs:108,95){(b)};
      \node at (axis cs:108,90){\footnotesize $u_h = u_{sd}$};
      \node at (axis cs:108,85){\footnotesize $h=1/64$};
      \node at (axis cs:108,80){\footnotesize $N=4225$};
      \node at (axis cs:108,75){\footnotesize $\epsilon = 0.066$};
      \node at (axis cs:108,70){\footnotesize $u_{\max} = -1.046$};
      \node at (axis cs:108,65){\footnotesize $u_{\min} = -0.048$};
      \node at (axis cs:108,60){\footnotesize $T_{\text{sol}} = 0.12s$};
      % links unten:
      \node[draw=black] at (axis cs:25,25) {\includegraphics[width=0.33\linewidth]{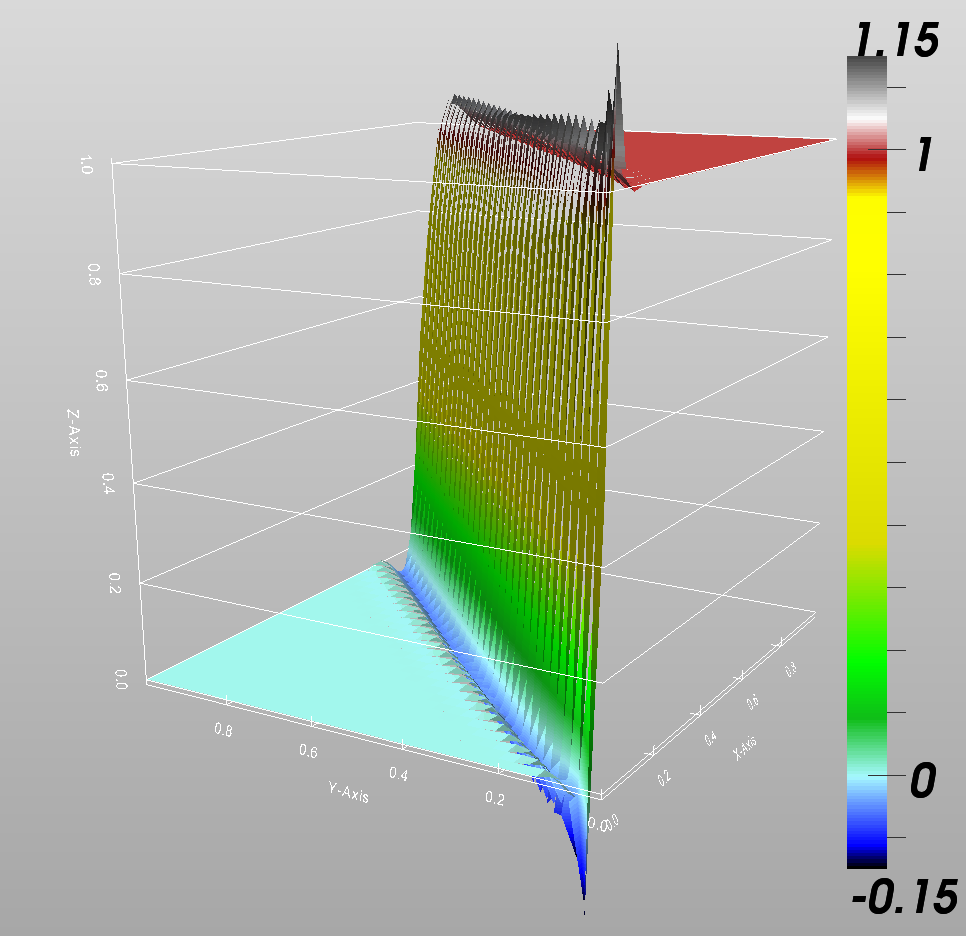}};
      \node at (axis cs:-8,45){(c)};
      \node at (axis cs:-8,40){\footnotesize $u_h = u_{dg}$};
      \node at (axis cs:-8,35){\footnotesize $h=1/32$};
      \node at (axis cs:-8,30){\footnotesize $N=4096$};
      \node at (axis cs:-8,25){\footnotesize $\epsilon = 0.062$};
      \node at (axis cs:-8,20){\footnotesize $u_{\max} = -1.249$};
      \node at (axis cs:-8,15){\footnotesize $u_{\min} = -0.249$};
      \node at (axis cs:-8,10){\footnotesize $T_{\text{sol}} = 0.08s$};
      % rechts unten:
      \node[draw=black] at (axis cs:75,25) {\includegraphics[width=0.33\linewidth]{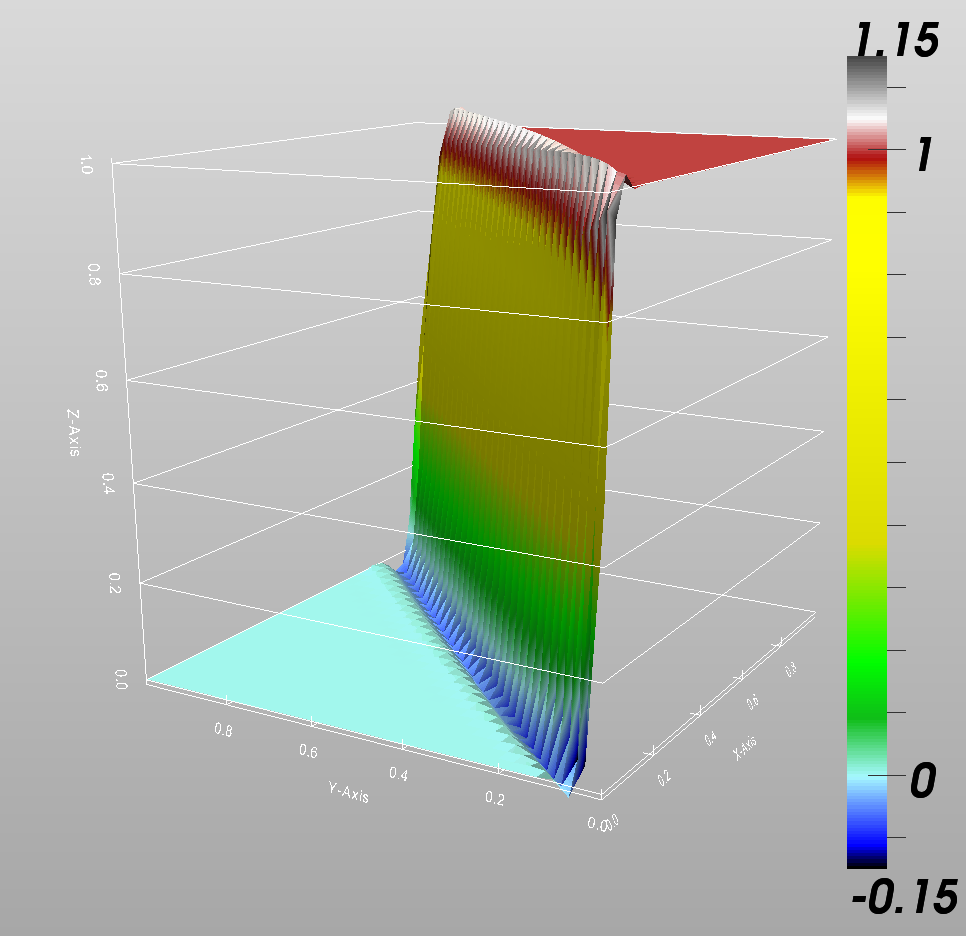}};
      \node at (axis cs:108,45){(d)};
      \node at (axis cs:108,40){\footnotesize $u_h = u_{cg}$};
      \node at (axis cs:108,35){\footnotesize $h=1/32$};             
      \node at (axis cs:108,30){\footnotesize $N=1089$};
      \node at (axis cs:108,25){\footnotesize $\epsilon = 0.069$};
      \node at (axis cs:108,20){\footnotesize $u_{\max} = -1.042$};
      \node at (axis cs:108,15){\footnotesize $u_{\min} = -0.042$};
      \node at (axis cs:108,10){\footnotesize $T_{\text{sol}} = 0.10s$}; 
  \end{axis}
  \end{tikzpicture}
  \caption[$2$-D example by L\'opez and Sinus\'ia: Warped plots]
          {Warped plots (all from the same perspective) of the numerical
    solution for different methods: (a) and (b): $u_{sd}$ is the SDFEM solution, (c): $u_{dg}$
    is the DG($1$) solution, (d): $u_{cg}$ is the diffusive $L^2$-projection of
    $u_{dg}$ (DG($1$)+L2). $h$ is the uniform meshsize, $N$ is the dimension of the
    solution space and $\epsilon=\|u-u_h\|_{L^2(\Omega)}$.
    $u_{\max}$ and $u_{\min}$ are the maximal and minimal values of the
    numerical solution $u_h$.
    $T_{\text{sol}}$ is the solution time (matrix assembly + linear solver).
  }
  \label{fig:Lopez_warp}
\end{figure}
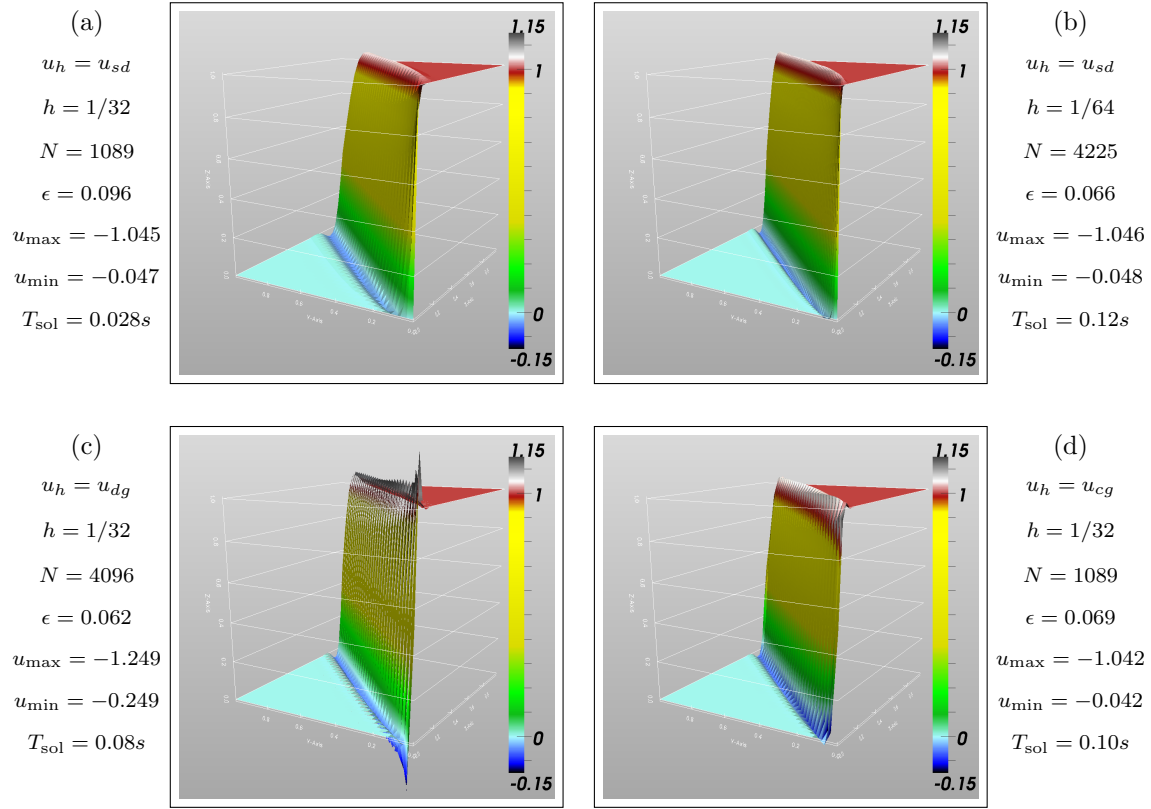

\begin{figure}[H]
  \centering
  \begin{tabular}{ccc}
  (a) $L=3 (h=1/32)$
  &
  (b) $L=4 (h=1/64)$
  &
  (c) $L=5 (h=1/128)$ \\
    \includegraphics[width=5cm]{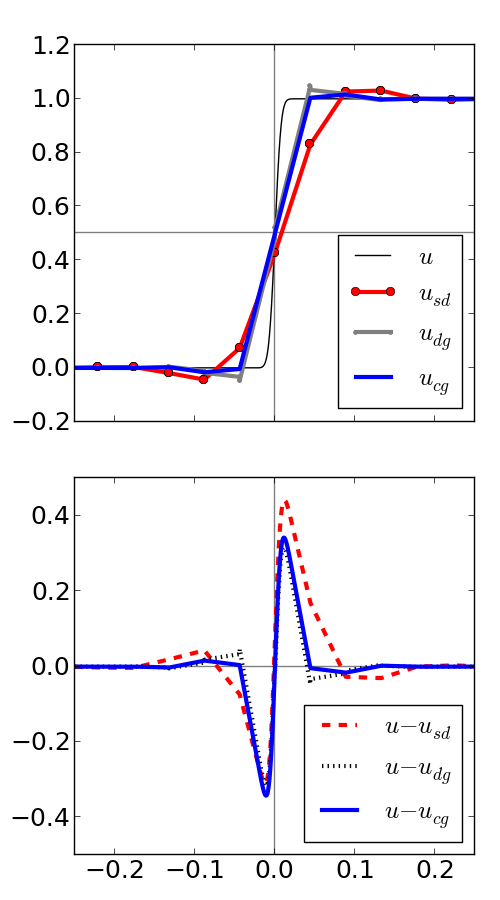}
    &
    \includegraphics[width=5cm]{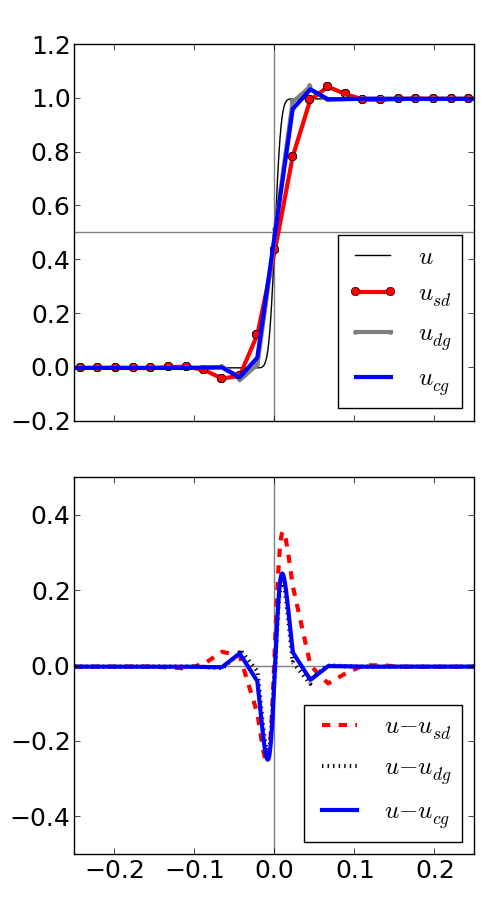}
    &
    \includegraphics[width=5cm]{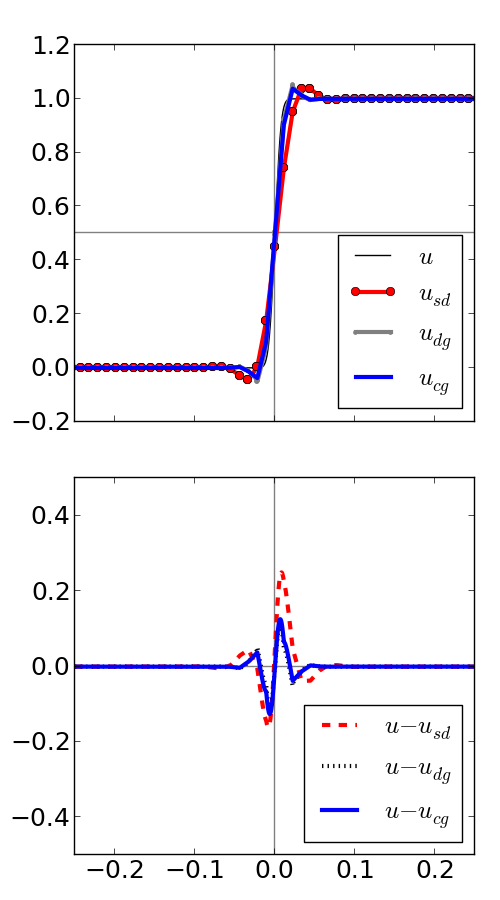}
  \end{tabular}
  \caption[$2$-D example by L\'opez and Sinus\'ia: Zoomed error plots]
          {Zoomed plots of the solution and the absolute error for different
    methods along the diagonal line connecting (1,0) and (0,1).
    $L$ is the global refinement level.
    $u$ is the true solution resolved on a very fine mesh ($h=10^{-5}$).
  }
  \label{fig:Lopez_diagonal_plots}
\end{figure}

Using the test problem from subsection \ref{paragraph:Lopez}, we take
a closer look
at the quality of the solution with respect to smearing effects and numerical
over- and undershoots around the characteristic layer.
We compare the post-processed DG($1$) method with the SDFEM method on
structured meshes.
Figure \ref{fig:Lopez_warp} shows the 3-D profile of four different numerical
solutions on the whole domain $[0,1]^2$.
Figure \ref{fig:Lopez_diagonal_plots} uses cross-sectional plots over
the diagonal line between $(0,1)$ and $(1,0)$ to zoom into the steep front.

\medskip
\noindent Observations from Figures \ref{fig:Lopez_warp} and \ref{fig:Lopez_diagonal_plots}:
\begin{itemize}
\item[] \hspace{-5mm} Near the discontinuity in the boundary condition:
\item[1.] On the same refinement level, DG($1$) exhibits larger over- and
  undershoots than SDFEM (see Figures \ref{fig:Lopez_warp} (a)+(c)).
\item[2.] The diffusive $L^2$-projection has a dampening effect on the DG($1$)
  solution, reducing the amount of large over- and undershoots significantly
  (see Figure \ref{fig:Lopez_warp}(c)+(d))
  without smearing out the steep front beyond a mesh cell
  (see Figure \ref{fig:Lopez_diagonal_plots} solution plots).
\end{itemize}

\begin{itemize}
\item[] \hspace{-5mm} Globally:
\item[3.] Small over- and undershoots in the DG($1$) solution are merely dampened by
  the diffusive $L^2$-Projection (see Figure \ref{fig:Lopez_diagonal_plots}).
\item[4.] On the same refinement level, both DG($1$) and DG($1$)+L2 capture
  the location of the steep front more accurately than SDFEM throughout the
  domain (see Figures \ref{fig:Lopez_warp}(a)+(c) and
  \ref{fig:Lopez_diagonal_plots}).
\item[5.] SDFEM on refinement level $L+1$ captures the steep front as well as
  DG($1$) or DG($1$)+L2 on level $L$ (comparable number of DOF)
  (see Figure \ref{fig:Lopez_warp}(b)+(c) and Figure \ref{fig:Lopez_diagonal_plots}:
  blue line in (a) vs. red line in (b))
\end{itemize}

\noindent
To achieve a comparable number of DOF as for the DG($1$) method, the SDFEM
method requires one extra level of global mesh refinement.
The resulting matrix assembly time for SDFEM on the refined mesh is higher
than for DG($1$) on the coarse mesh.

\subsection{Forward transport in \texorpdfstring{$\boldsymbol{2}$}{2}-D}\label{example:forward2D}

In the following, we demonstrate the applicability of the presented DG
methods to more realistic scenarios.

We solve the groundwater flow equation (\ref{eqn:GWE}) for the hydraulic head
distribution $\phi$ and evaluate the velocity field (\ref{eqn:Darcy})
on the structured mesh $\mathcal{T}_{h_0}$. 
The solute transport equation (\ref{eqn:TPE}) is then solved 
\begin{itemize}
\item using adaptive DG($1$) on a hierarchy of adaptively refined meshes
  $\{\mathcal{T}^{\mbox{\tiny{adapt}}}_{\nu}\}_{\nu\in\mathbb{N}}$, or
\item using DG($k$) with diffusive $L^2$-projection on the same mesh $\mathcal{T}_{h_0}$, % with an additional post-processing step using a , 
  or
\item using SDFEM on the same mesh $\mathcal{T}_{h_0}$ or on a globally refined
  mesh $\mathcal{T}_1^{\mbox{\tiny{global}}}$.
\end{itemize}

\noindent
The Gaussian field $Y$ depicted in Figure \ref{fig:2D_transport_plotsA}(a) has the
physical size of $100 \times 100 [m^2]$, the resolution of the structured mesh
$\mathcal{T}_0$ is $100 \times 100$ cells. $Y$ is described by its mean value
$\beta=-6.0$, its variance $\sigma^2=1.0$ and
the correlation lengths $(\ell_x,\ell_y)=(10,10)[m]$.
The hydraulic head is prescribed on the left ($\phi\big|_{x=0}=100[m]$) and on the
right ($\phi\big|_{x=100}=99.5[m]$) boundaries. The induced pressure gradient drives the main
flow.

%\newgeometry{left=2.5cm,bottom=1cm,top=1cm,right=2.5cm}
%\thispagestyle{empty}
\onecolumn

%\restoregeometry
\twocolumn

The injection well parameters are $\injectionRate=5\times 10^{-4} [m^3/s]$, 
$\tilde{c}_{\hbox{\tiny{inj}}}=1[g/m^3]$ and $T_{\hbox{\tiny{inj}}} = 100[s]$.
The stationary solution values range between $0$ and $100[gs/m^3]$.
The molecular diffusion is $D_m=2\times 10^{-9} [m^{2}/s]$, 
the longitudinal and transversal dispersivities
are $\alpha_l=10^{-3}[m]$ and $\alpha_t=10^{-4}[m]$.

Theoretically, $u_h$ tends to the true solution $u$ as the meshsize $h$ tends to
$0$. At least, all numerical oscillations should disappear if the meshsizes near the 
characteristic layer become so small that the effective mesh P\'eclet numbers\footnote{Replace
  $\alpha_{\ell}$ by $\alpha_t$ in (\ref{eqn:Peclet_longitudinal}) to consider
  longitudinal effects.}
are smaller than $1$.
Working with $\mathcal{O}(\alpha_t) \sim 10^{-4}$ and a base level meshsize
$\mathcal{O}(h)\sim 1$ we would get $\mathcal{O}(\mathcal{P}^t_h)\sim
10^{3}$. It would require $10$ levels of global refinement to achieve
$\mathcal{O}(\mathcal{P}^t_h)\sim1$, but already after
$7$ levels of global refinement, our mesh would have more than $10^8$ cells.
The problem size would become extremely large for a $2$-D simulation. 

Adaptive mesh refinement is the only chance to keep the problem size orders of
magnitude lower
%within an acceptable level. 
while reducing mesh P\'eclet numbers at the steep fronts.
We choose the stopping criterion (\ref{eqn:stop}) from section
\ref{paragraph:adaptive} with a tolerance of $p_{\text{osc}}=1\%$.
The result is achieved after $5$ steps of adaptive refinement (Table
\ref{t:2D_transport_dg1A} and Figure \ref{fig:2D_transport_plotsA}(b)).
This solution is taken as the reference
solution for assessing the quality of different methods on the
structured mesh (Table \ref{t:2D_transport_sdfem_dg_cg} and 
Figures \ref{fig:2D_transport_plotsA}) and \ref{fig:2D_transport_plotB}).

\begin{table}[H]
  %\footnotesize
  \centering
  \scalebox{0.69}{
    \ra{1.3}
    \hspace{-3mm}
    \begin{tabular}{@{}rrrrrrrrr@{}}
      \toprule
      \multicolumn{7}{c}{Adaptive DG($1$)}\\
      $L$ &   $DOF$  & max. $\mathcal{P}_h^t$  & $M[s]$  & $T[s]$ & $IT$  &  $TIT$  &   $u_{\min}$     &  $u_{\max}$ \\
      \midrule
      0  &  40,000  & 4675.89  &  0.91 &   0.04  &  4  & 0.011  & 134.47 &  -33.11  \\ % \\ %TIT=
      1  &  64,000  & 2232.44  &  1.49 &   0.08  &  4  & 0.018  & 132.26 &  -31.77  \\ % \\ %TIT=
      2  & 102,364  & 1105.47  &  2.42 &   0.25  &  8  & 0.033  & 127.38 &  -24.68  \\ % \\ %TIT=
      3  & 162,964  & 545.85   &  3.93 &   0.75  & 14  & 0.050  & 112.48 &  -12.10  \\ % \\ %TIT=
      4  & 259,372  & 272.94   &  6.29 &   2.29  & 26  & 0.076  & 104.99 &   -3.39  \\ % \\ %TIT=
      5  & 413,860  & 136.48   & 10.26 &   7.29  & 52  & 0.110  & 100.76 &   -0.38  \\ % \\ %TIT=
      \bottomrule
    \end{tabular}
  }
  \caption{\footnotesize Adaptive mesh refinement on {\tt UG}, with $p_r=20[\%]$ and $p_c=10[\%]$.
            %, total linear solution time after $7$ steps $=29.4[s]$
            Renumbering mesh cells on mesh level $L=5$ takes $0.16$ sec. 
            %Maximal $\mathcal{P}_h^t$ on the smallest mesh cell ranges from $4,676$ down to $32$.
            Linear solver used: BiCGSTAB + ILU(0) with reduction $10^{-8}$.
            $M=$ matrix assembly time,
            $IT=$ linear solver iterations,
            $T=$ linear solver time
          }
          \label{t:2D_transport_dg1A}
          %  \end{center}
\end{table}

\vspace{-5mm}
\begin{table}[H]
  %\footnotesize
  \centering
  \scalebox{0.69}{
    \ra{1.3}
    \hspace{-2mm}
    \begin{tabular}{@{}crrrrrr@{}}
      \toprule
      \multicolumn{7}{c}{Global $h$- or $p$-refinement on a structured mesh} \\
      Method &  $DOF$ &  $M[s]$ & $IT$  &  $T[s]$ &  $u_{min}$  & $u_{max}$ \\
      \midrule
      SDFEM ($L=0$) & 10,201 & 0.2 & 12  &  0.02  &  -13.79 & 116.79   \\ 
      SDFEM ($L=1$)  & 40,401 & 0.8 & 28   & 0.26  &  -15.29 & 113.43  \\
      \midrule
      \midrule
      DG($1$)  & 40,000 & 0.6 &  4 & 0.04  & -33.11 & 134.47 \\
      diffusive $L^2$-proj.   & 10,201 & +0.05 &  1  & +0.02 & -3.44  & 103.13 \\
      \midrule
      DG($2$)   & 90,000 & 1.7 &  4 & 0.13  & -34.04 & 135.88 \\
      diffusive $L^2$-proj.   & 10,201 & +0.09 &  1  & +0.02 & -3.27  & 104.82 \\
      \midrule
      DG($3$)  & 160,000 & 4.5 &  4 & 0.26  & -36.69 & 135.94 \\
      diffusive $L^2$-proj.   & 10,201 & +0.16 & 1 & +0.02 & -3.80  & 103.73 \\
      \bottomrule
    \end{tabular}
  }
  \caption
  {\footnotesize Computations on a structured mesh:
            SDFEM on mesh level $L=0$ and $L=1$ versus 
            post-processed DG methods on mesh level $L=0$ with different polynomial
            orders. Renumbering mesh cells on $L=0$ for the DG methods takes $0.004$ sec.
            Linear solver used for solving the transport equation:
            BiCGSTAB + ILU(0) with reduction $10^{-8}$.
            Linear solver used for the diffusive $L^2$-projection:
            BiCGSTAB + AMG with reduction $10^{-8}$.
            $M=$ matrix assembly time, $IT=$ linear solver iteration number, 
            $T=$ linear solver time.
          }
          \label{t:2D_transport_sdfem_dg_cg}
\end{table}

%\newgeometry{left=2.5cm,bottom=2cm,top=1.5cm,right=2.5cm}
\onecolumn
\begin{figure}[H]
  %\centering
  %\hspace{-1.3cm}
  \begin{tikzpicture}
  \begin{axis}[width=\linewidth,height=\linewidth,hide axis,clip=false]
      \addplot+[only marks,mark=.,mark size=1pt,draw=blue,fill=white] coordinates {
      (0.0,0.0)
      (0.0,100.0)
      (100.0,0.0)
      (100.0,100.0)
      };
      \node at (axis cs:50,79){\it 2-D solute transport: $u=m_0^c$};
      % 1. Reihe links:
      \node at (axis cs:5,74){(a)};
      \node[draw=black] at (axis cs:5,50) {\includegraphics[width=0.3\linewidth]{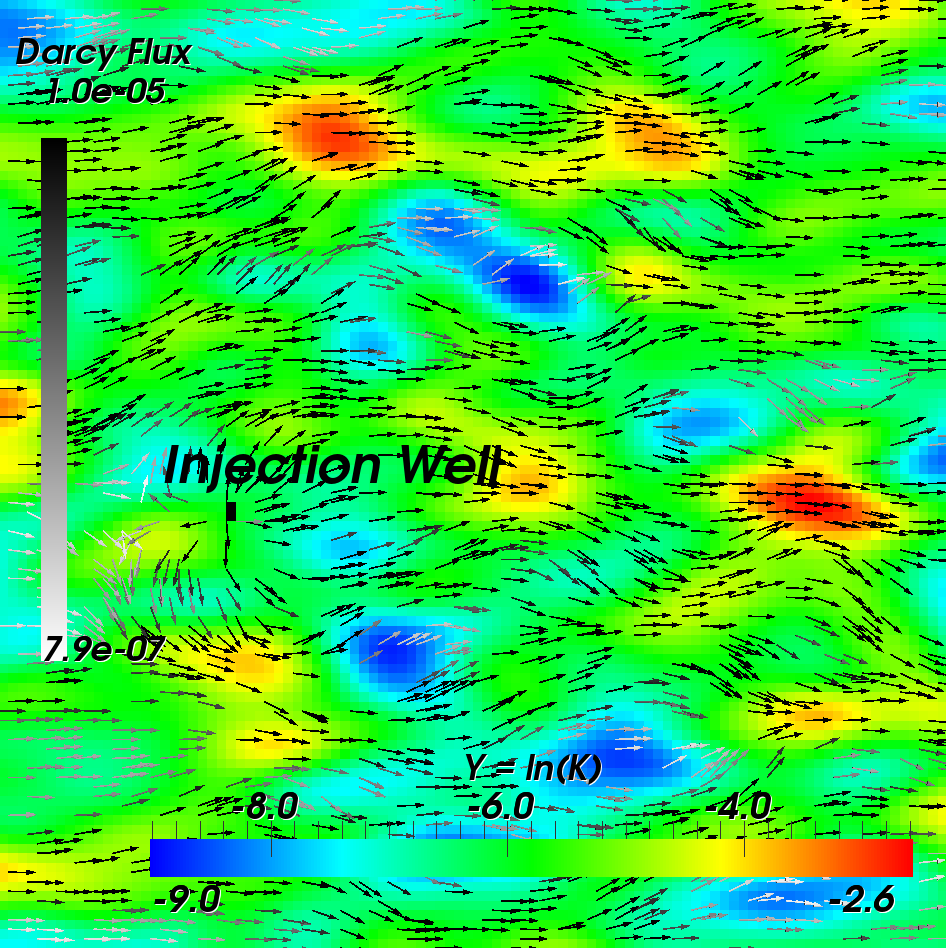}};
      % 1. Reihe mitte:
      \node at (axis cs:50,74){(b) adaptive DG($1$), level 7};
      \node[draw=black] at (axis cs:50,50) {\includegraphics[width=0.3\linewidth]{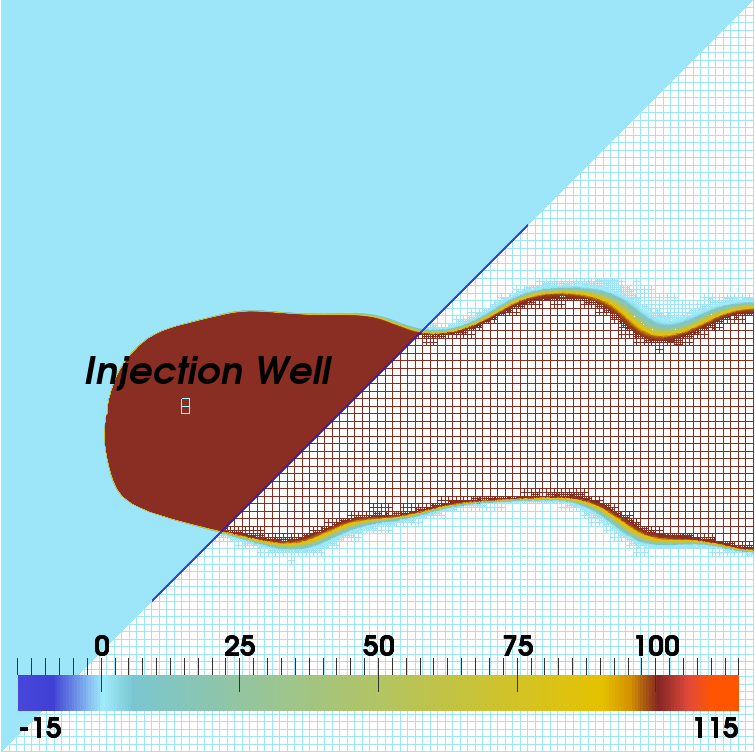}};
      \node at (axis cs:50,65){$u_{\max}=100.76$};
      \node at (axis cs:50,61){$u_{\min}=-0.38$};
      % 1. Reihe rechts:
      \node at (axis cs:95,74){(c) SDFEM};
      \node[draw=black] at (axis cs:95,50) {\includegraphics[width=0.3\linewidth]{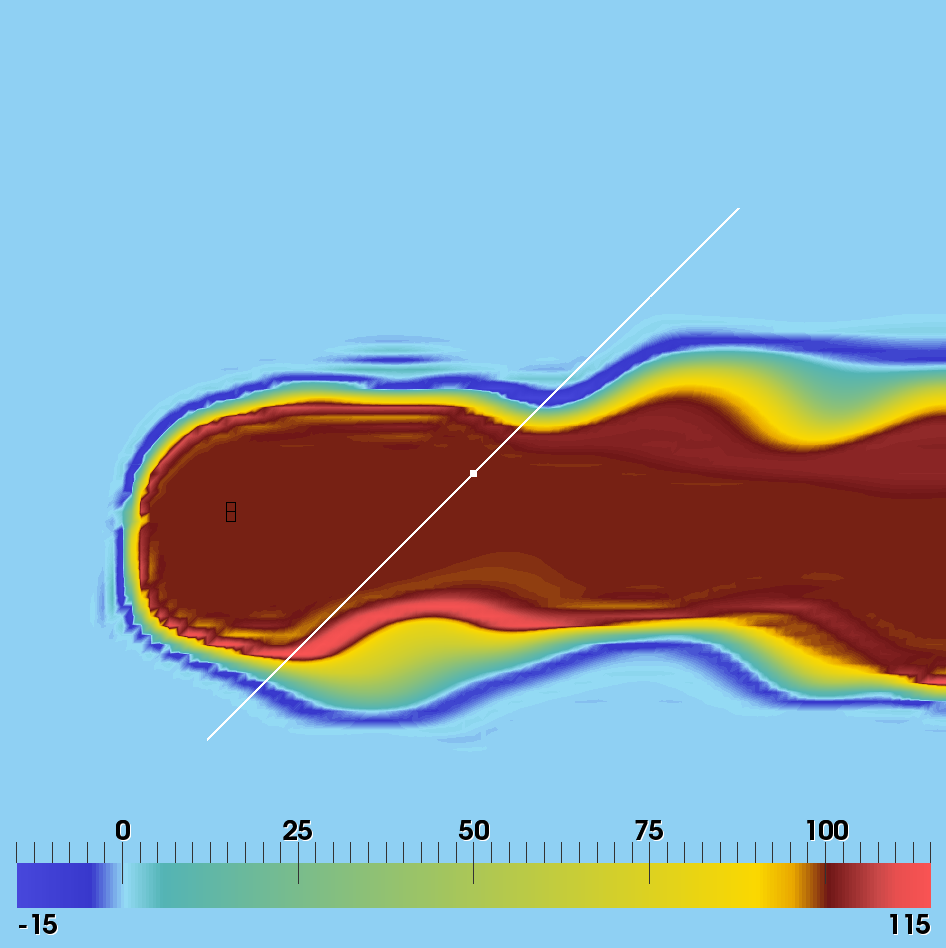}};
      \node at (axis cs:95,65){$u_{\max}=116.79$};
      \node at (axis cs:95,61){$u_{\min}=-13.79$};
      % 2. Reihe links:
      \node at (axis cs:5,24){(d) Post-processed DG($3$)};
      \node[draw=black] at (axis cs:5,00) {\includegraphics[width=0.3\linewidth]{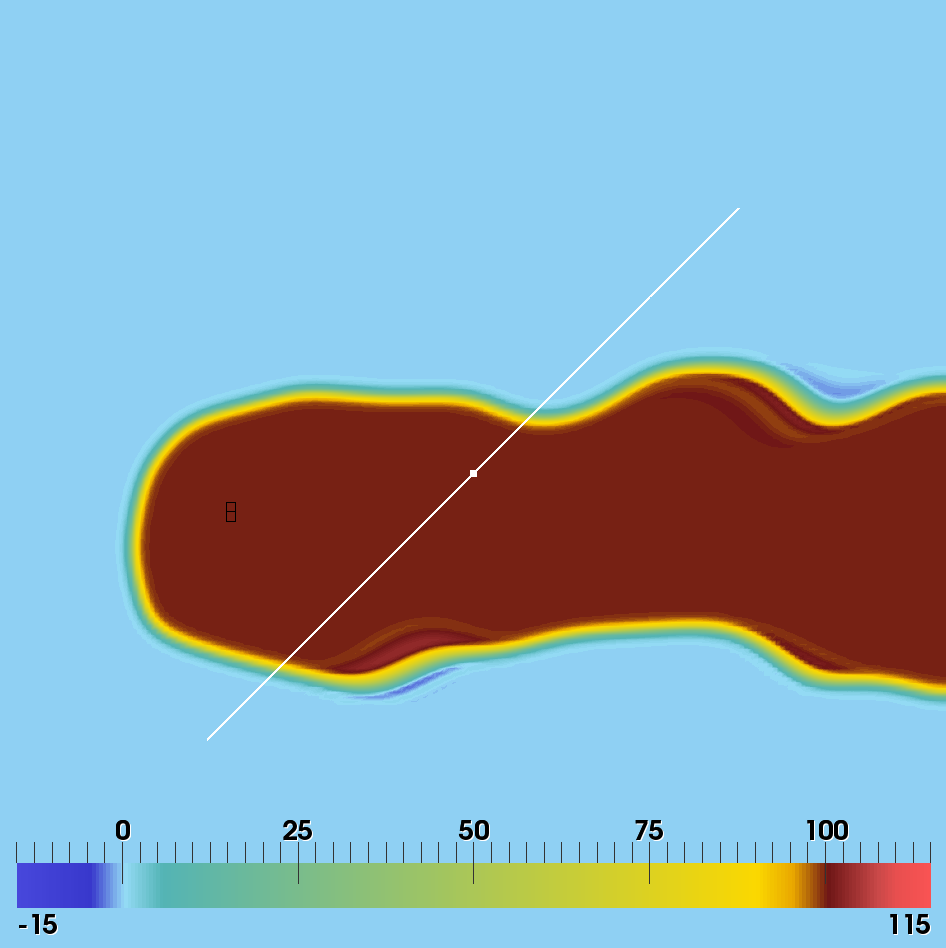}};
      \node at (axis cs:5,15){$u_{\max}=103.73$};
      \node at (axis cs:5,11){$u_{\min}=-3.80$};
      % 2. Reihe mitte:
      \node at (axis cs:50,24){(e) Post-processed DG($2$)};
      \node[draw=black] at (axis cs:50,00) {\includegraphics[width=0.3\linewidth]{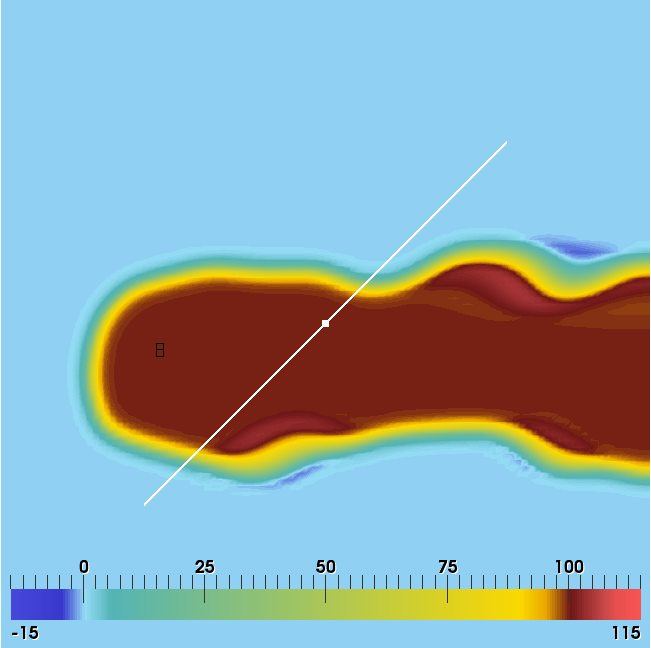}};
      \node at (axis cs:50,15){$u_{\max}=104.82$};
      \node at (axis cs:50,11){$u_{\min}=-3.27$};
      % 2. Reihe rechts:
      \node at (axis cs:95,24){(f) Post-processed DG($1$)};
      \node[draw=black] at (axis cs:95,00) {\includegraphics[width=0.3\linewidth]{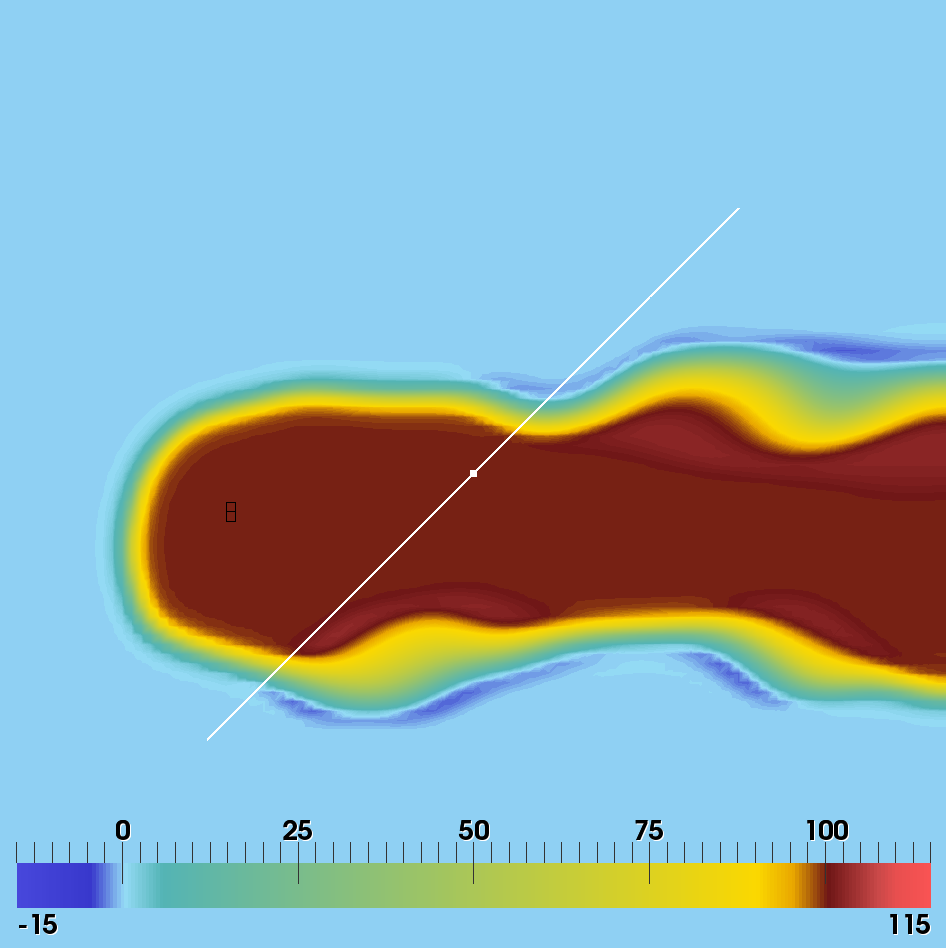}};
      \node at (axis cs:95,15){$u_{\max}=103.13$};
      \node at (axis cs:95,11){$u_{\min}=-3.44$};
  \end{axis}
  \end{tikzpicture}
  \caption{
    This Figure is continued in Figure \ref{fig:2D_transport_plotsB}.
    (a) Gaussian field $Y=\ln(K)$ and velocity field $\vec{q}_h$ on $\mathcal{T}_0$.
    (b) Adaptive DG($1$) solution with $p_{\text{osc}}\le 1\%$ is taken to be the reference solution for the
    stationary transport problem.
    (c) - (f) Comparing different solutions on the coarse structured mesh $\mathcal{T}_0$ with
    $100\times 100$ cells.
  }
  \label{fig:2D_transport_plotsA}
\end{figure}
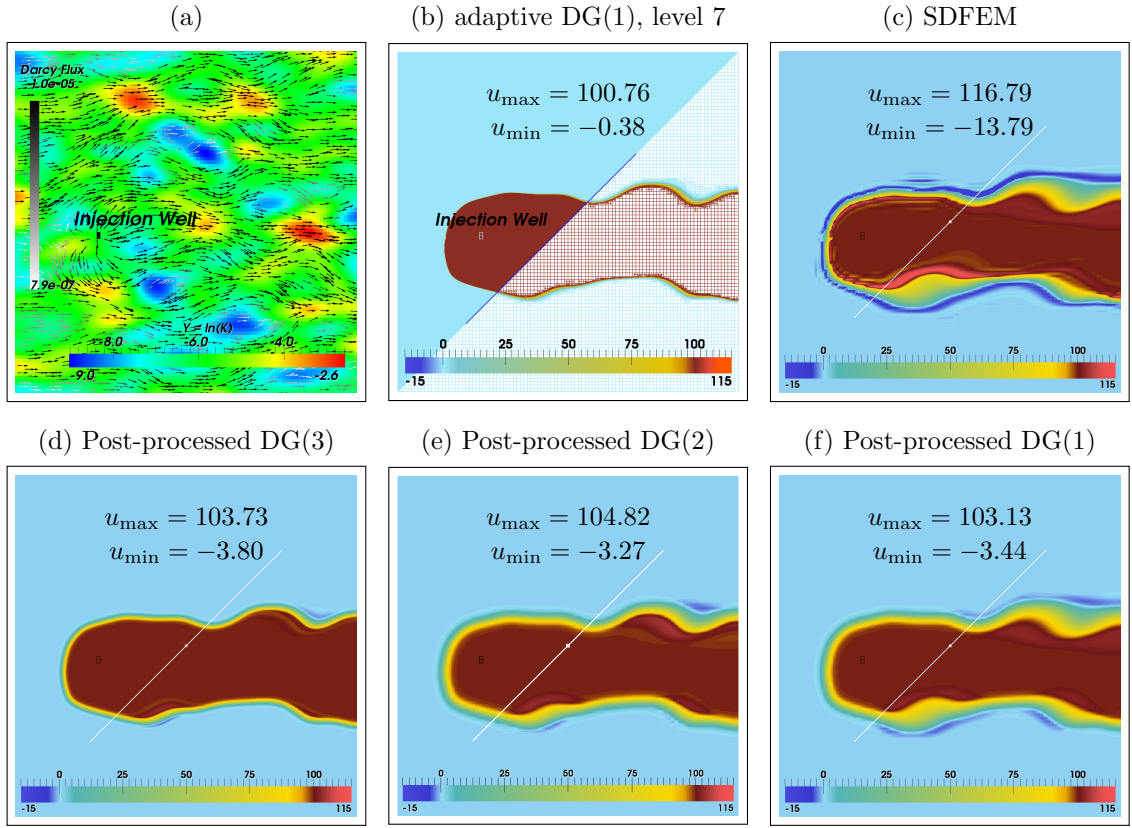
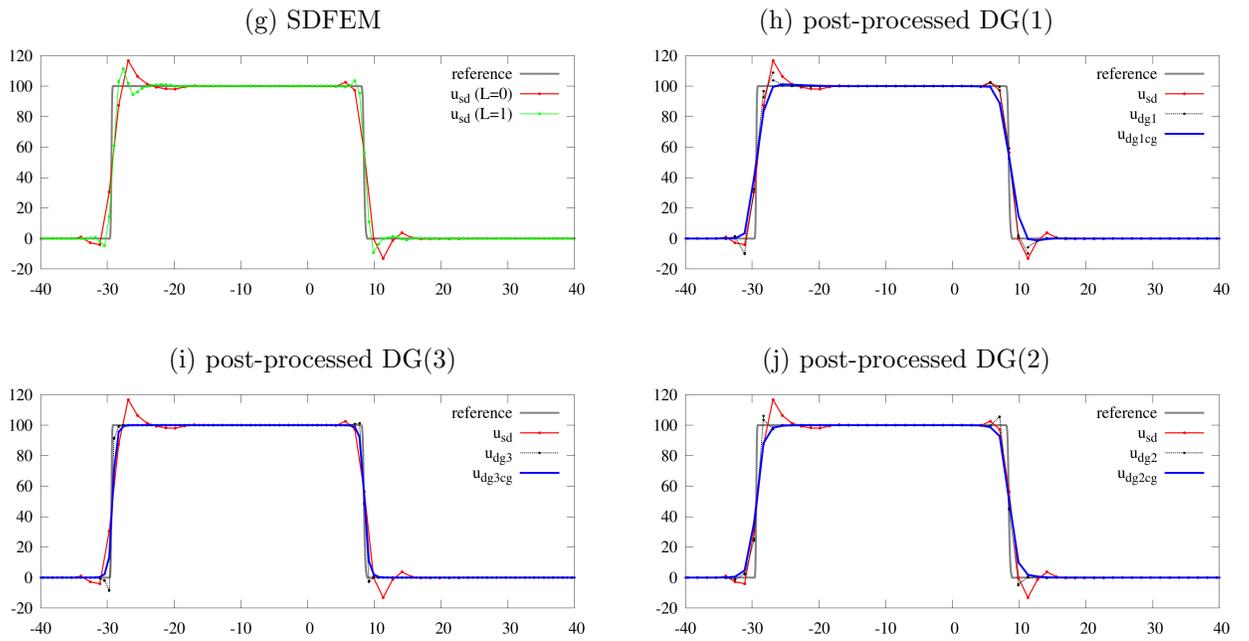
\begin{figure}[H]
  \begin{tikzpicture}
    \begin{axis}[width=\linewidth,height=\linewidth,hide axis,clip=false]
      \addplot+[only marks,mark=.,mark size=1pt,draw=blue,fill=white] coordinates {
      (0.0,0.0)
      (0.0,100.0)
      (100.0,0.0)
      (100.0,100.0)
      };
      \node at (axis cs:50,71){\it Cross sectional plots along the diagonal line indicated above:};
      % 3. Reihe links:
      \node at (axis cs:15,66){(g) SDFEM};
      \node at (axis cs:12,48) {\includegraphics[width=0.54\linewidth]{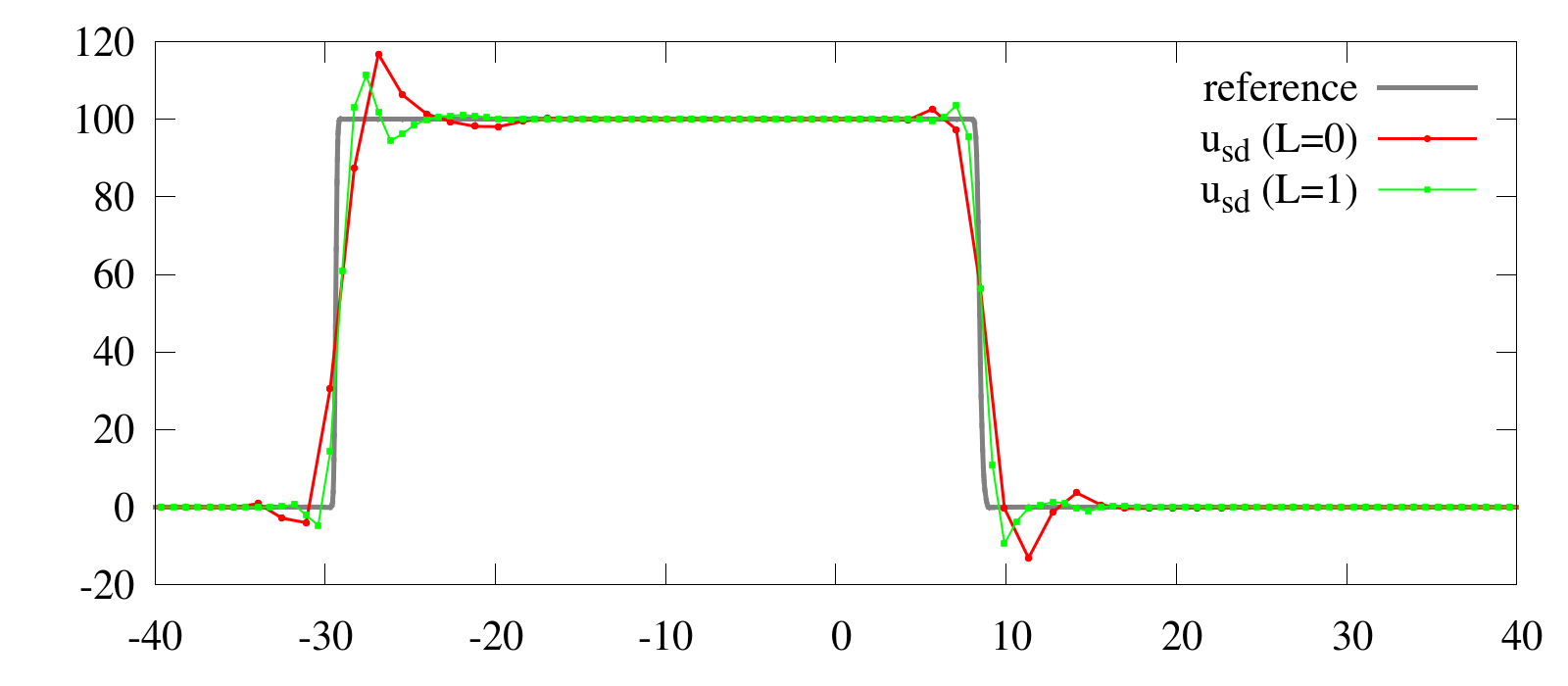}};
      % 3. Reihe rechts:
      \node at (axis cs:85,66){(h) post-processed DG($1$)};
      \node at (axis cs:88,48) {\includegraphics[width=0.54\linewidth]{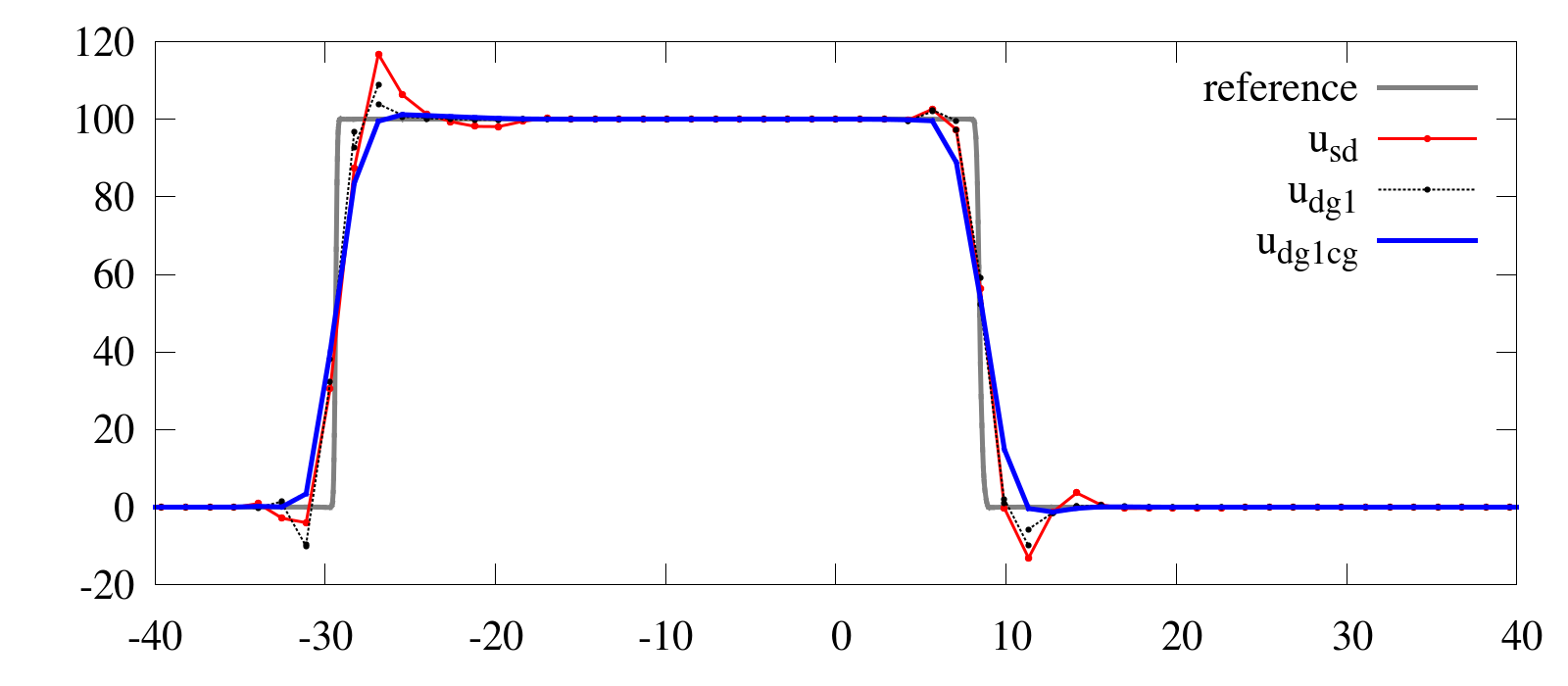}};
      % 4. Reihe links:
      \node at (axis cs:15,26){(i) post-processed DG($3$)};
      \node at (axis cs:12,8) {\includegraphics[width=0.54\linewidth]{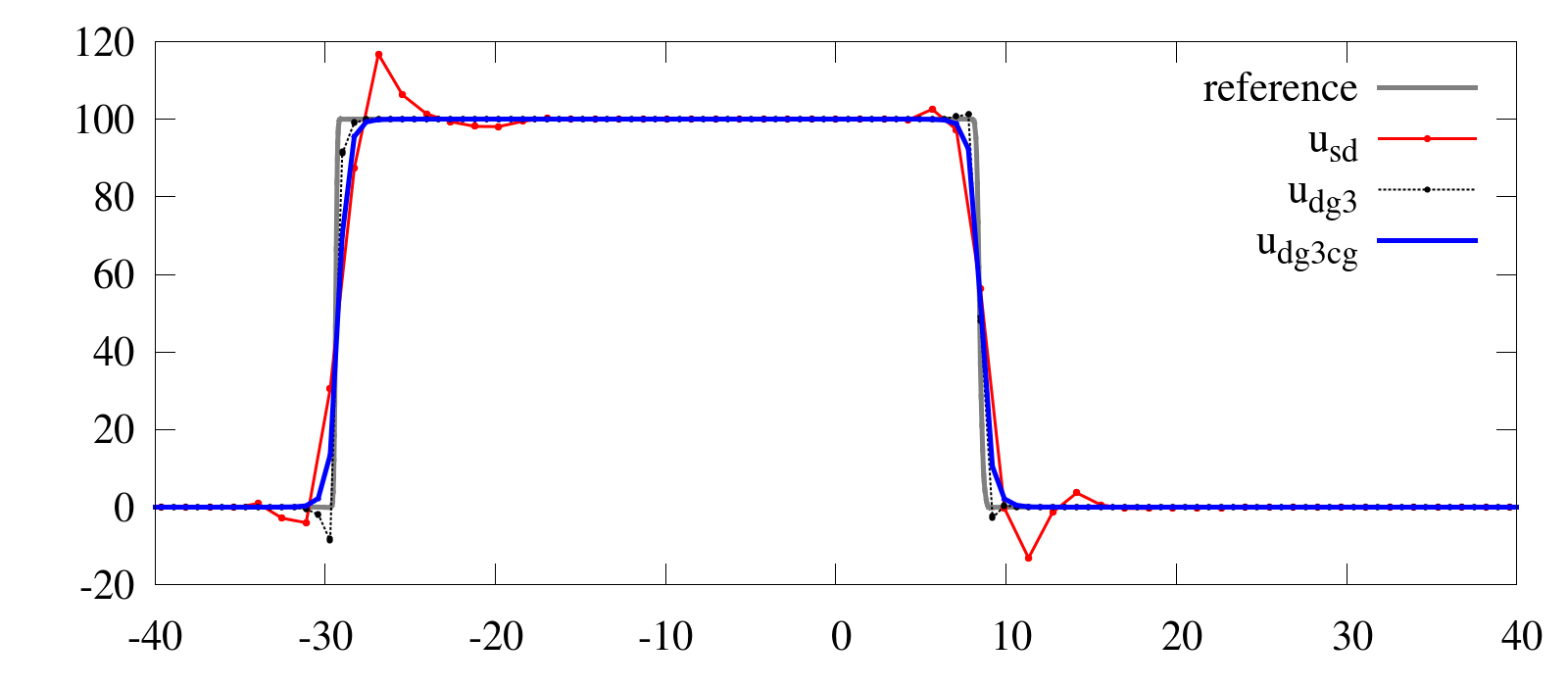}};
      % 4. Reihe rechts:
      \node at (axis cs:85,26){(j) post-processed DG($2$)};
      \node at (axis cs:88,8) {\includegraphics[width=0.54\linewidth]{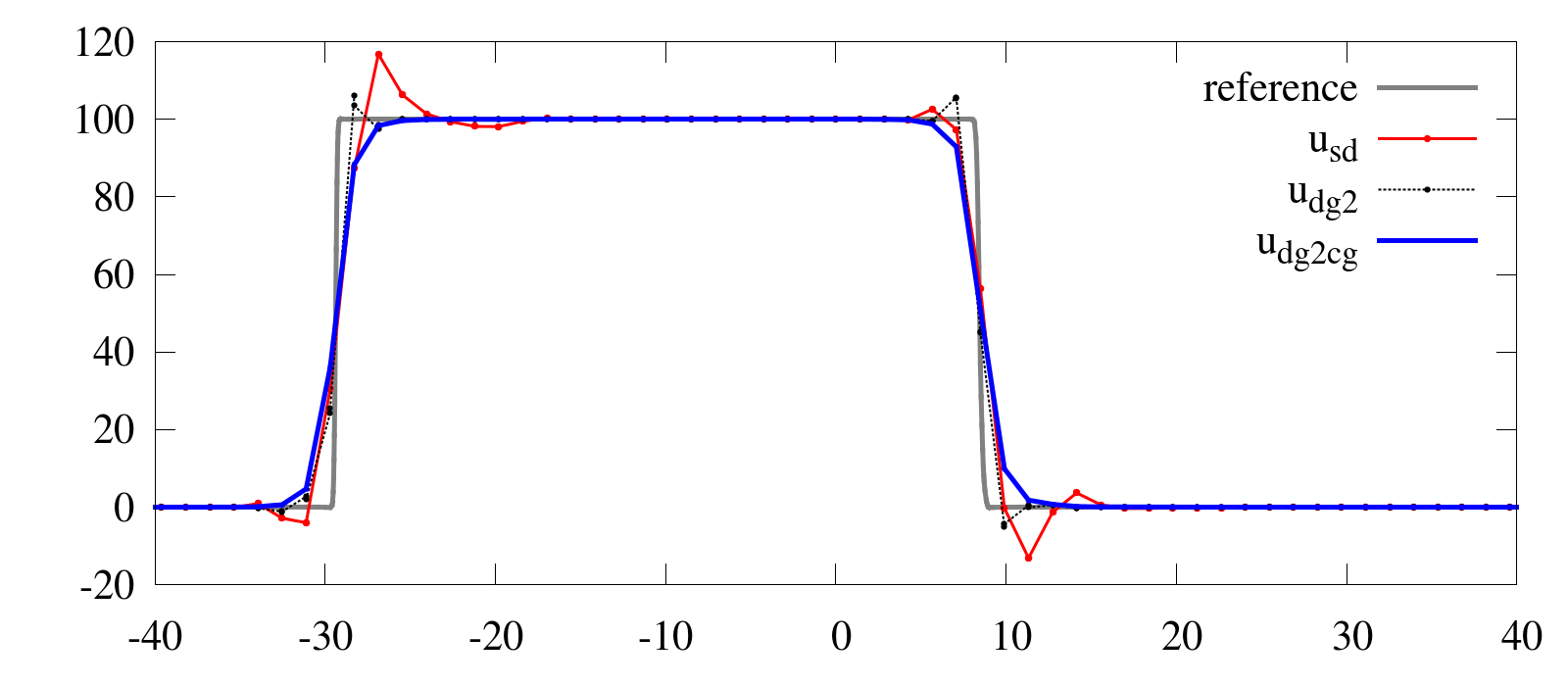}};
  \end{axis}
  \end{tikzpicture}
  \vspace{-8mm}
  \caption{
    (g) - (j):
    SDFEM solution $u_{sd}$ on $\mathcal{T}_0$ and on
    $\mathcal{T}_1^{\mbox{\scriptsize{global}}}$
    compared to DG($k$) solutions 
    ($u_{\text{dgk}}$) on $\mathcal{T}_0$ and post-processed DG($k$) solutions 
    ($u_{\text{dgkcg}}$) on $\mathcal{T}_0$
    along the indicated cut-line (for $k=1,2,3$).
    $u_{\max}$ and $u_{\min}$ are the maximal and minimal values of the
    numerical solution.
  }
  \label{fig:2D_transport_plotsB}
\end{figure}
%\restoregeometry
\twocolumn

\noindent {\it Observations:}
\begin{enumerate}
\item On the coarse mesh $\mathcal{T}_0$, DG($k$) solutions with $k=1,2,3$ may
  locally exhibit stronger over- and undershoots than the SDFEM solution, but these 
  can be reduced very efficiently with a diffusive $L^2$-projection
  (see the runtimes of $M$ and $T$ for the diffusive $L^2$-projection in Table
  \ref{t:2D_transport_sdfem_dg_cg} and 
  Figure \ref{fig:2D_transport_plotsB}(h)-(j)).
\item The over- and undershoots of the SDFEM solution may oscillate into the
  domain surrounding the steep front (see Figure \ref{fig:2D_transport_plotsA}(c)) 
  whereas over- and undershoots of the
  DG solutions do not show this behavior.
\item On the same mesh $\mathcal{T}_0$, higher order polynomials can be used
  to improve the qua\-lity of the DG solution with respect to the sharpness of the steep front. 
  The areas with over- and undershoots are shrunk. 
  A diffusive $L^2$-projection preserves this behavior
  and reduces over- and undershoots
  (see Figure \ref{fig:2D_transport_plotsA}(d)-(f)).
\item The SDFEM method applied to a globally refined mesh
  $\mathcal{T}_1^{\mbox{\tiny{global}}}$ requires a comparable number of
  unknowns as the DG($1$) method on $\mathcal{T}_0$.
  The steep front is resolved similarly well, but the matrix assembly and the
  linear solution takes longer 
  ($L=1$ in Table \ref{t:2D_transport_sdfem_dg_cg} and Figure \ref{fig:2D_transport_plotsB}).
\end{enumerate}

\subsection{Forward transport in \texorpdfstring{$\boldsymbol{3}$}{3}-D}\label{example:forward3D}

A $3$-D simulation with a similar setting is started on a coarse mesh $\mathcal{T}_{h_0}$
with $32\times 32\times 32$ cells. 
The domain extensions are $(L_x,L_y,L_z)=(10,10,10)[m]$.
The geostatistical field parameters for the Gaussian model are $\beta=-6.0$, $\sigma^2=1.0$ and 
$(\ell_x,\ell_y,\ell_z)=(2,2,1)[m]$.
The hydraulic head is prescribed on the left boundary by $\phi\big|_{x=0}=10[m]$ and on the
right boundary by $\phi\big|_{x=10}=9.8[m]$. 
An injection well is placed at the position $(x,y)=(2.1,5.1)[m]$ and its
$z$-range is $[0...-5][m]$.
% A tracer concentration of $1[g/m^3]$ is added through the injection well
% at a rate of $0.010 [m^3/s]$ over a period of $100[s]$. Integration over time
% yields a stationary solution with values ranging between $u_{\min}=0$ and
% $u_{\max}=100[gs/m^{3}]$.
The injection well parameters are $\tilde{w}=1\times 10^{-2} [m^3/s]$,
$\tilde{c}=1[g/m^3]$ and $T_{\hbox{\tiny{inj}}} = 100[s]$.
The values of the stationary solution range between $0$ and $100[gs/m^{3}]$.

Doing the same analysis as for the $2$-D case, we find that with global refinement,
we would end up with more than $10^{13}$ cells after $10$ refinement
steps in order to achieve $\mathcal{O}(\mathcal{P}^t_h)\sim1$.
The adaptive DG($1$) solution in Table \ref{t:3D_transport_dg1A} after $9$
steps shows a sharp resolution of the steep front (Figure
\ref{fig:3D_transport_plotsB}(c)). However, there are thin layers where the
over- and undershoots exceed $25\%$. This is still far away from being an acceptable reference
solution. 
% Nonetheless, this solution can be used to judge how well the steep
% fronts are resolved by the different solutions on the coarse structured mesh.
% We observe a similar behavior as in the $2$-D case.
To achieve our targeted reduction of under- and overshoots below $p_{\text{osc}}=5\%$, 
further refinement steps with increasing memory consumption and solution time
are necessary.
% A further refinement on the used hardware was not possible, 
% because in the following refinement step, the number of unknowns increases up to 
% $89,267,592$. 
% The memory consumption for the linear solver thereby exceeded 
% the limit.
%
% Regarding the quality of the diffusive $L^2$-projection applied to the
%DG($1$) and DG($3$) solutions, very similar observations as in the  $2$-D
%case can be made (Figure \ref{fig:3D_transport_plotsB}(e),(f) and 
%Table \ref{t:3D_transport_sdfem_dg_cg}).

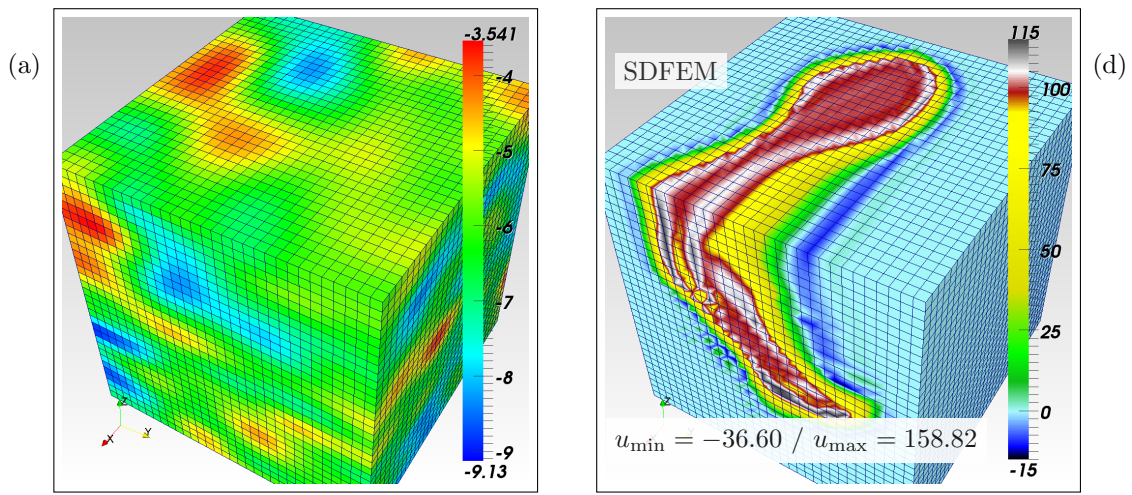
\begin{figure}[H]
  \centering
  \vspace{-3cm}
  \begin{tikzpicture}
    \begin{axis}[width=\linewidth,height=1.4\linewidth,hide axis,clip=false]
      \addplot+[only marks,mark=.,mark size=1pt,draw=white,fill=white] coordinates {
      (0.0,0.0)
      (0.0,100.0)
      (100.0,0.0)
      (100.0,100.0)
      };
      % links oben:
      \node[draw=black] at (axis cs:18,36) {\includegraphics[width=0.41\linewidth]{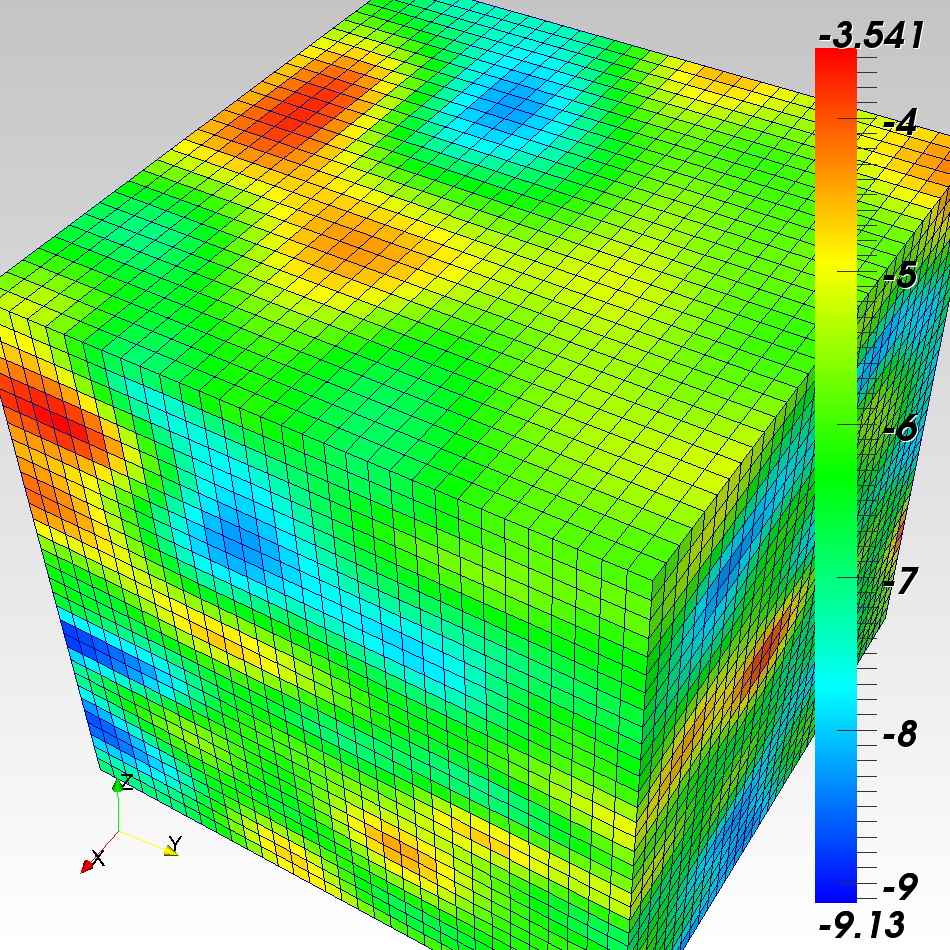}};
      \node at (axis cs:-14,51){\normalsize(a)};
      % rechts oben:
      \node[draw=black] at (axis cs:82,36) {\includegraphics[width=0.41\linewidth]{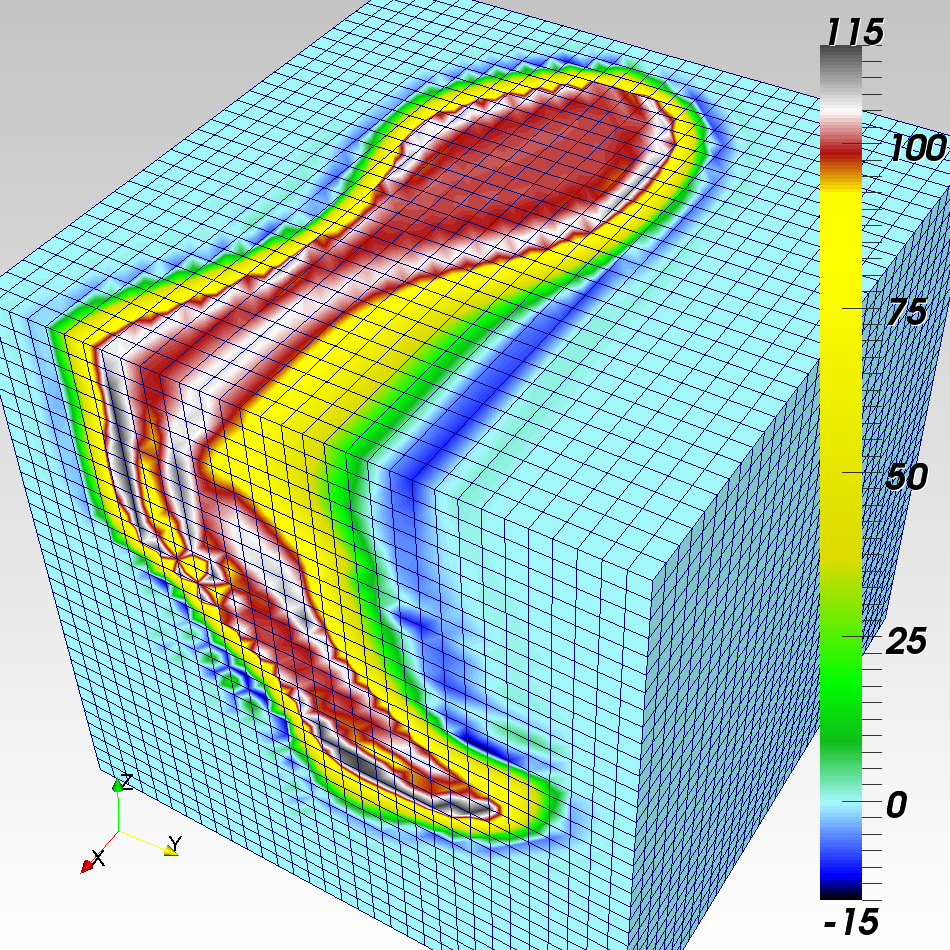}};
      \node at (axis cs:114,51){\normalsize(d)};
      \node[fill=white,fill opacity=0.8] at (axis cs:62.2,50.6){\normalsize{SDFEM}};
      \node[fill=white,fill opacity=0.9] at (axis cs:77,20.5){$u_{\min}=-36.60~/~u_{\max}=158.82$};
  \end{axis}
  \end{tikzpicture}
  \caption{
    This Figure is continued in Figure \ref{fig:3D_transport_plotsB}.
    (a): Gaussian field $Y=\ln(K)$ on $\mathcal{T}_0$ with $32\times
    32\times 32$ cells. 
  }
  \label{fig:3D_transport_plotsA}
\end{figure}

\begin{figure}[H]
  \vspace{-2cm}
  \begin{tikzpicture}
  \begin{axis}[width=\linewidth,height=1.4\linewidth,hide axis,clip=false]
      \addplot+[only marks,mark=.,mark size=1pt,draw=blue,fill=white] coordinates {
      (0.0,0.0)
      (0.0,100.0)
      (100.0,0.0)
      (100.0,100.0)
      };
      % links unten:
      \node[draw=black] at (axis cs:18,60) {\includegraphics[width=0.41\linewidth]{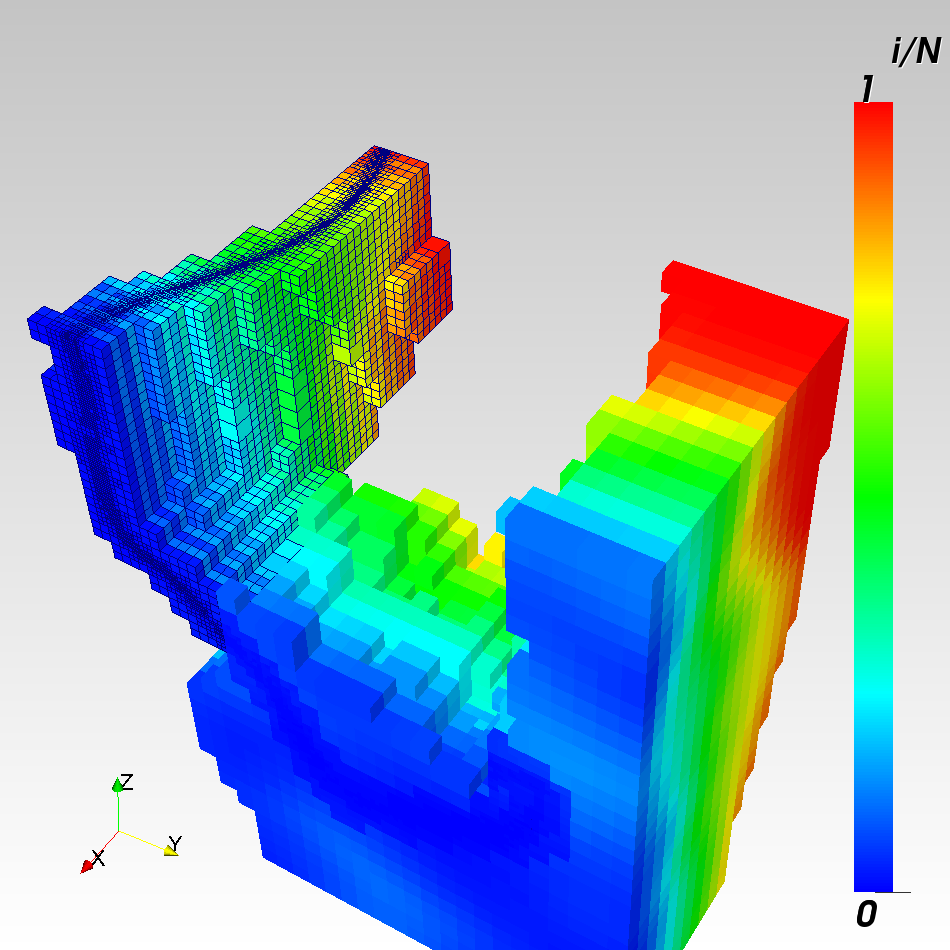}};
      \node at (axis cs:-14,75){\normalsize(b)};
      % rechts unten:
      \node[draw=black] at (axis cs:82,60) {\includegraphics[width=0.41\linewidth]{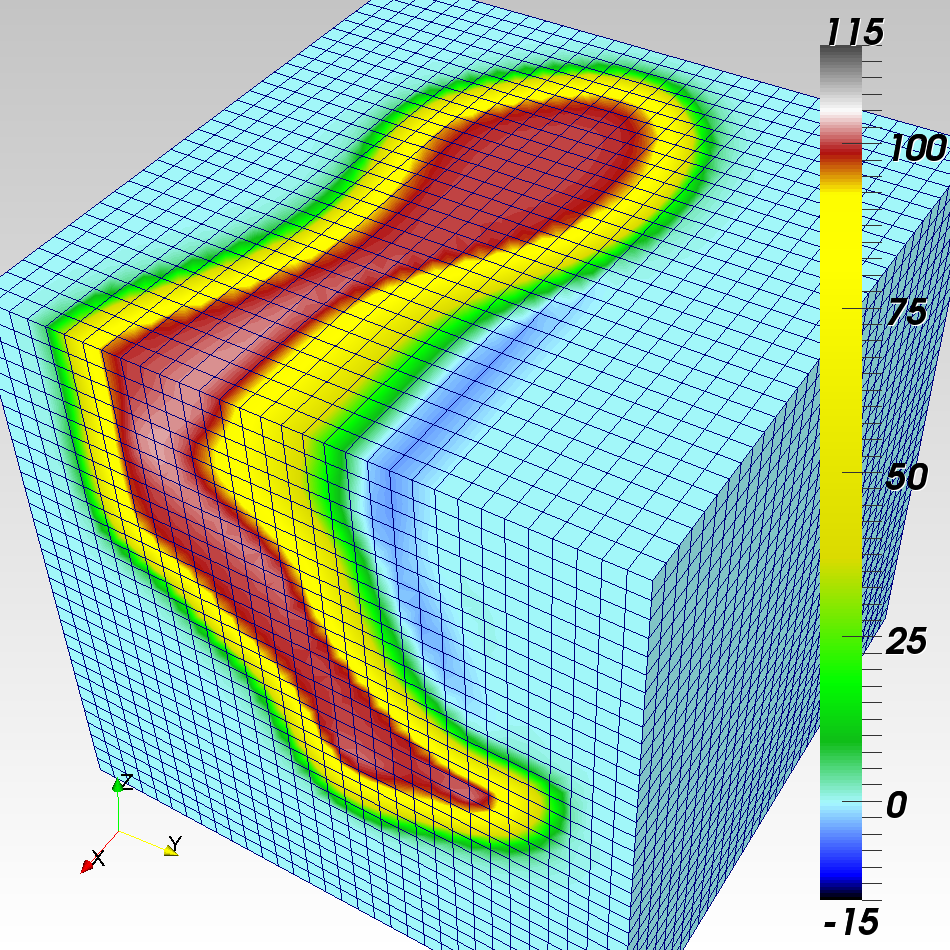}};
      \node at (axis cs:114,75){\normalsize(e)};
      \node[fill=white,fill opacity=0.8] at (axis cs:63.6,74.3){\normalsize{DG(1)+$L^2$}};
      \node[fill=white,fill opacity=0.8] at (axis cs:77,44.6){$u_{\min}=-4.28~/~u_{\max}=102.77$};
      % links unten:
      \node[draw=black] at (axis cs:18,14) {\includegraphics[width=0.41\linewidth]{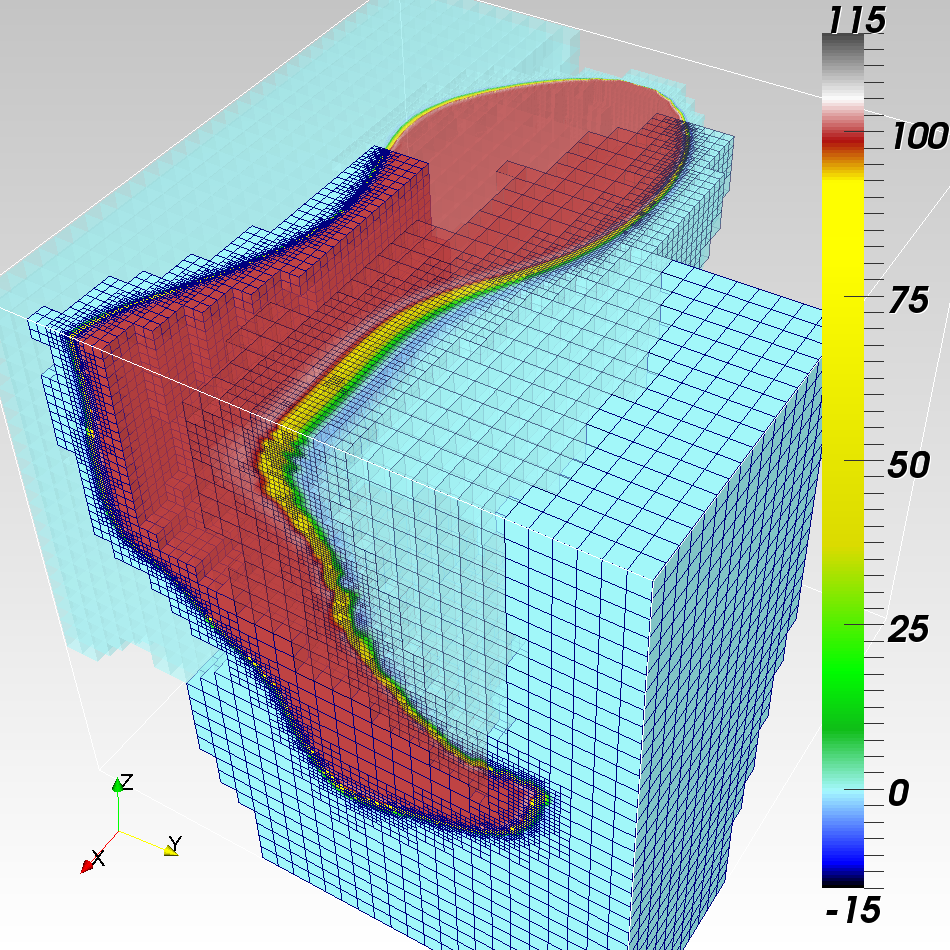}};
      \node at (axis cs:-14,29){\normalsize(c)};
      \node[fill=white,fill opacity=0.8] at (axis cs:4.0,28.3){\normalsize{Adaptive DG($1$)}};
      \node[fill=white,fill opacity=0.9] at (axis cs:11.7,-1.3){$u_{\min}=-26.69~/~u_{\max}=127.34$};
      % rechts unten:
      \node[draw=black] at (axis cs:82,14) {\includegraphics[width=0.41\linewidth]{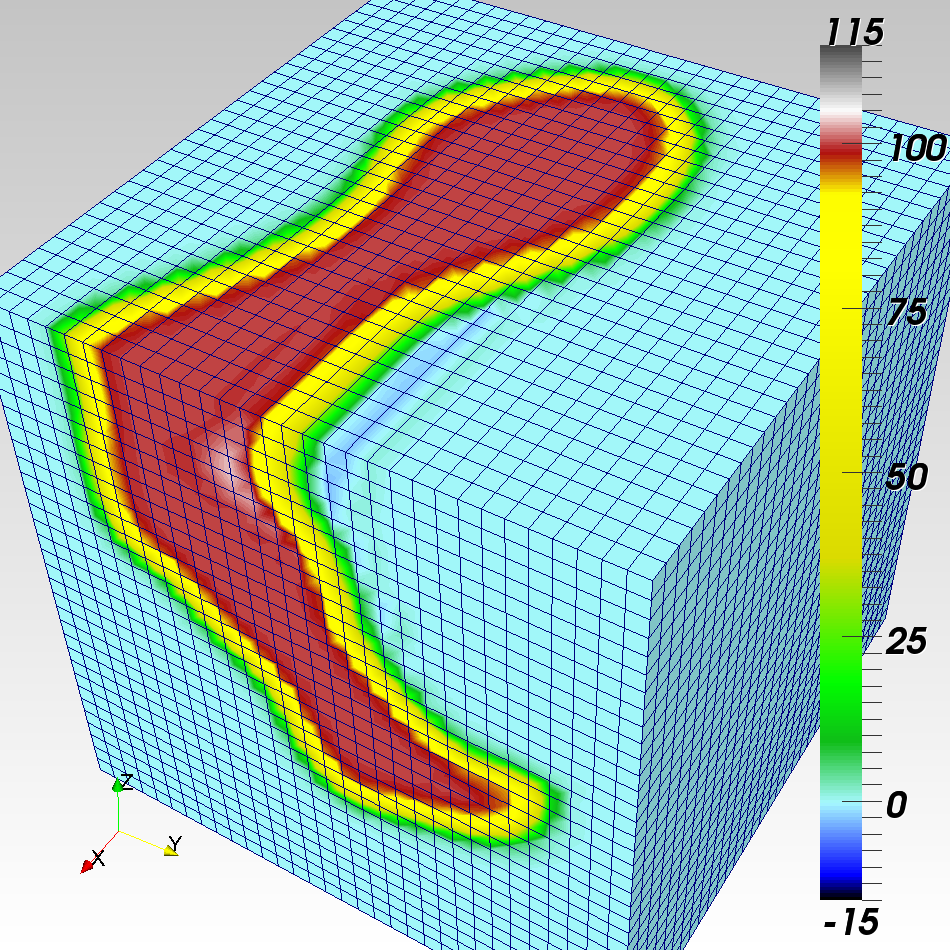}};
      \node at (axis cs:114,29){\normalsize(f)};
      \node[fill=white,fill opacity=0.8] at (axis cs:63.6,28.3){\normalsize{DG(3)+$L^2$}};
      \node[fill=white,fill opacity=0.9] at (axis cs:77,-1.3){$u_{\min}=-1.03~/~u_{\max}=102.14$};
  \end{axis}
  \end{tikzpicture}
  \vspace{-5mm}
  \caption{
    (b)-(c): Illustrating parallel adaptive refinement with dynamic
    load-balancing (step $9$). A reduction of under- and overshoots below $5\%$
    is possible, see Table \ref{t:3D_transport_dg1A}.
    Fig.\ref{fig:3D_transport_plotsA}(d),Fig.\ref{fig:3D_transport_plotsB}(e)+(f): Comparing different solutions on the coarse mesh $\mathcal{T}_0$.
    $u_{\max}$ and $u_{\min}$ are the maximal and minimal values of the
    displayed numerical solution.
  }
  \label{fig:3D_transport_plotsB}
\end{figure}
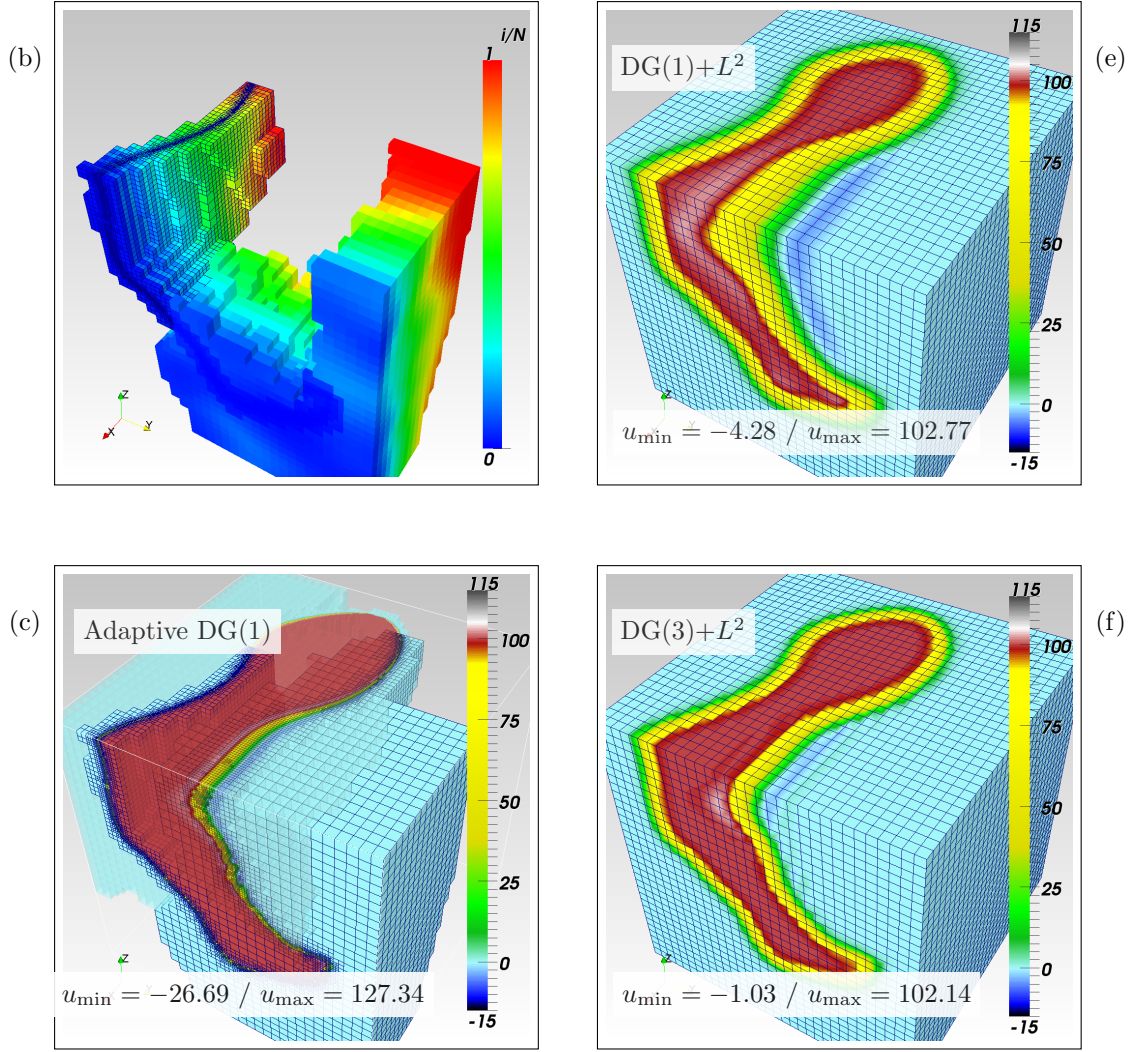

We make very similar observations as in the $2$-D case:
\begin{enumerate}
\item 
  The over- and undershoots generated by the DG($k$) solutions with $k=1,2,3$ 
  can be reduced very efficiently with a diffusive $L^2$-projection
  (see the runtimes of $M$ and $T$ for the diffusive $L^2$-projection in Table
  \ref{t:3D_transport_sdfem_dg_cg}).
\item The over- and undershoots of the SDFEM solution may oscillate into the
  domain surrounding the steep front (see Figure \ref{fig:3D_transport_plotsA}(d)) 
  whereas over- and undershoots of the DG solutions do not show this behavior.
\item On the same mesh $\mathcal{T}_0$, higher order polynomials can be used
  to improve the qua\-lity of the DG solution with respect to the sharpness of the steep front. 
  The areas with over- and undershoots are shrunk. 
  A diffusive $L^2$-projection preserves this behavior
  and reduces over- and undershoots
  (see Figure \ref{fig:3D_transport_plotsB}(e)+(f)).
\item The SDFEM method applied to a globally refined mesh
  $\mathcal{T}_1^{\mbox{\tiny{global}}}$ requires a comparable number of
  unknowns as the DG($1$) method on $\mathcal{T}_0$.
  The steep front is resolved similarly well (not shown here),
  but the matrix assembly and the
  linear solution takes longer 
  ($L=1$ in Table \ref{t:2D_transport_sdfem_dg_cg}).
\end{enumerate}

%\newgeometry{left=2.5cm,bottom=2cm,top=2cm,right=2.5cm}
%\thispagestyle{empty}
\onecolumn

%\restoregeometry
\twocolumn

%\newgeometry{left=2.5cm,bottom=2cm,top=2cm,right=2.5cm}
\begin{table}[H]
  %  \hspace{-4mm}
  %\footnotesize
  \centering
  \scalebox{0.73}{
    \ra{1.3}
    \begin{tabular}{@{}crrrrrr@{}}
      \toprule
      \multicolumn{7}{c}{Adaptive DG($1$)}\\
      $L$  &   \qquad $DOF$      & \quad  $IT$  & $T[s]$   &  \qquad $u_{\min}$ & \quad   $u_{\max}$ & \quad $\mathcal{P}^t_h$ \\
      \midrule
   0 & 346,560    &   8    & 0.7     & -38.8 & 136.59 &  1562 \\
   1 & 399,024    &   11   & 1.1     & -38.7 & 142.70 &  1560 \\ 
   2 & 537,240    &   12   & 1.6     & -44.8 & 144.25 &  1560 \\ 
   3 & 814,152    &   12   & 3.2     & -50.7 & 141.17 &  1560 \\ 
   4 & 1,410,136  &   18   & 8.2     & -50.2 & 142.59 &  389  \\ 
   5 & 2,740,704  &   17   & 15     & -56.3 & 142.78 &  194  \\ 
   6 & 5,275,496  &   18   & 31     & -43.1 & 147.91 &  194  \\ 
   7 & 9,273,552  &   24   & 88     & -38.8 & 151.62 &  97   \\ 
   8 & 15,368,112 &   28   & 152    & -28.9 & 132.08 &  97   \\ 
   9 & 23,806,944 &   35   & 349    & -26.7 & 127.34 &  48   \\ 
  10 & 39,897,040  &   52  & 743    & -24.0 & 121.37 &  48   \\ 
  11 & 59,903,672  &   64  & 1425   & -18.7 & 117.17 &  24   \\ 
  12 & 70,498,808  &   67  & 1935   & -15.2 & 114.24 &  24   \\ 
  13 & 95,837,712  &   87  & 3926   & -15.5 & 112.58 &  24   \\ 
  14 & 132,771,680 &   95  & 5724   & -11.7 & 113.12 &  24   \\ 
  15 & 188,298,456 &  110  & 8854   & -10.5 & 110.78 &  24   \\ 
  16 & 242,272,328 &  128  & 13489    &  -6.88 & 106.72 &  24   \\ 
  17 & 321,793,400 &  167  & 28075  &  -6.59 & 106.73 &  24   \\ 
  18 & 348,929,992 &  185  & 30403  &  -5.30 & 105.37 &  24   \\ 
  19 & 433,481,584 &  200  & 38251  &  -5.23 & 104.85 &  12   \\ 
  20 & 474,804,264 &  238  & 61351  &  -4.98 & 104.57 &  12   \\ 
    \bottomrule
  \end{tabular}
}
  \caption[$3$-D forward problem: Parallel adaptive DG]
          {$3$-D parallel adaptive refinement on {\tt ALUGrid}, 
            with $p_r=70[\%]$ and $p_c=5[\%]$.
            Renumbering mesh cells on level $20$ takes $24.5$ sec.
            Linear solver used: BiCGSTAB + SSOR with reduction $10^{-8}$.
            $L=$ refinement level, $IT=$ linear solver iterations, $T=$ linear solver
            time.
            The computation is performed on {\tt quadxeon4}. 
            %on a multi-processor
            %machine equipped with $80$ cores 
            %({\it Intel\textregistered  Xeon\textregistered Processor E7-4870 @ 2.40GHz})
            %and $1$ TB RAM.
            $P=16$ cores are used for the computation, more than $66\%$ of the RAM
            is required for the linear solver in step $20$ alone.
            % In the $13$-th refinement step, the total number of DOF is $88,695,432$.
            % The linear solver requires a memory consumption which exceeds $128$ GB main memory.
          }
  \label{t:3D_transport_dg1A}
\end{table}

\begin{table}[H]
  %  \hspace{-4mm}
%  \footnotesize
  \centering
  \scalebox{0.73}{
    \ra{1.3}
    \begin{tabular}{@{}crrrrrrr@{}}
      \toprule
      \multicolumn{7}{c}{Global $h$- or $p$-refinement on a structured mesh}\\
      Method         &  $DOF$      &  $M[s]$  & $IT$  &  $T[s]$ &  $u_{\min}$  & $u_{\max}$ \\
      \midrule
      SDFEM ($L=0$)   & 58,806   & 0.5 & 14  &  0.11  &  -36.6 & 158.8   \\ 
      SDFEM ($L=1$)   & 363,350  & 2.7 & 22  &  1.12  & -36.5  & 139.1   \\
      %    +$L^2$-proj.    & 327,080  &  2  & +0.57  & -4.91   & 110.32 \\
      \midrule
      \midrule
      DG($1$)           &   376,832  & 2.3 & 7 & 0.38    & -37.9 & 136.4 \\
      diffusive $L^2$-proj.    &   58,806   & +0.2 & 2 & +0.1   & -4.3  & 102.8 \\
      \midrule
      DG($2$)           & 1,271,808  & 22.5 & 8 & 3.8    & -50.8 & 146.7 \\
      diffusive $L^2$-proj.    &   58,806   & +0.5 & 2 & +0.1   & -3.07  & 103.4 \\
      \midrule
      DG($3$)           & 3,014,656 &  190.1 & 9 & 38.4  & -47.1 & 152.0 \\
      diffusive $L^2$-proj.   & 58,806    &    +1.2 & 2 & +0.1    & -1.0  & 102.1 \\
      \bottomrule
    \end{tabular}
  }
  \caption[$3$-D forward problem: DG($k$) vs. SDFEM on structured meshes]
          {
            $3$-D parallel computations with $P=8$ cores on {\tt fna} (see
            Table \ref{t:IWRmachines}) using
            a structured mesh with
            partitioning $(P_x,P_y,P_z)=(1,8,1)$ and overlap $=1$.
            SDFEM on mesh levels $L=0$ and $L=1$ compared to 
            DG methods on the coarse mesh level $L=0$ with different polynomial orders.
            Parallel renumbering of mesh cells for the DG methods on level $L=0$ takes
            $0.004$ sec.
            Linear solver used for solving the transport equation:
            BiCGSTAB + ILU(0) with reduction $10^{-8}$.
            Linear solver used for the diffusive $L^2$-projection:
            BiCGSTAB + AMG with reduction $10^{-8}$.
            $M=$ matrix assembly time, $IT=$ number of
            iterations, $T=$ linear solver time.
            %    All computations are performed on a multi-processor
            %    machine equipped with $48$ cores 
            %    ({\it AMD Opteron\texttrademark Processor 6172 @ 2.10 GHz})
            %    and $128$ GB RAM.
          }
          \label{t:3D_transport_sdfem_dg_cg}
\end{table}

%\restoregeometry
\newpage
\subsection*{Conclusion and Outlook}
For the solution of the steady-state convection-dominant transport equation,
we have compared Discontinuous Galerkin (DG) methods 
to the Streamline Diffusion (SDFEM) method.
Put\-ting special emphasis on a practical application, 
we have analyzed efficiency and accuracy.
Two main issues occurring in the solution of the 
convection-dominated transport equation were tackled: 
\begin{enumerate} 
  \item the efficient reduction of numerical under- and overshoots,
  \item the efficient solution of the arising linear systems.
\end{enumerate}

With respect to the efficiency (solution time)
and the quality of the solution (maximal amplitude of the over- and undershoots and smearing effects at the steep fronts) 
for the convective-dominant transport problem,  
the observations made in the sections 
\S \ref{paragraph:Lopez}--\S \ref{example:forward3D} 
favor the combination CCFV / DG($1$) + diffusive $L^2$-projection 
over the FEM / SDFEM approach. 
Considering computing time 
to be the ultimate measure of available hardware resources, we observed: 
\begin{itemize}
\item On the same mesh level, the DG solutions resolve the steep fronts
more sharply than the SDFEM solution. 
In order to obtain the same level of 
accuracy as DG($1$), SDFEM would have to work on a globally refined mesh. 
As a consequence, the SDFEM approach would take longer than DG($1$).
\item In heterogeneous fields, the layers of spurious
oscillations generated by SDFEM may spread into the surrounding domain whereas
they stay localized for the DG method.
\item The diffusive $L^2$-projection is able to reduce the over- and undershoots
without increasing smearing effects beyond the mesh size.
\end{itemize}
Therefore, DG($1$) post-processed by a diffusive $L^2$-projection offers an efficient and more accurate alternative to the well-known SDFEM method.

Without doubt, the best possible solution in terms of the $L^2$-error is achieved with adaptive mesh refinement.
However, numerical oscillations can be reduced to an acceptable level 
only if the mesh cells at the steep front 
become so small that their local mesh P\'eclet numbers approach $1$ (diffusion-dominant problem).
This comes at a very high price, especially in $3$-D.
A ``perfect'' solution in this sense 
may not be necessary for a stable 
inversion scheme that can cope with noisy data. 

Hence, regarding the integration of the forward solvers into an inversion framework,
we recommend the combination CCFV / DG($1$) post-processed by a diffusive $L^2$-projection 
for the solution of steady-state transport problems with high mesh P\'eclet numbers.
This combination works on the same structured mesh on which the hydraulic conductivity is
resolved and keeps the implementation of the inversion scheme simple.

For future developments, a natural extension of the presented methods 
is a combination of $h$- or $hp$-adaptive DG
with the diffusive $L^2$-projection on unstructured meshes (with hanging nodes refinement).

Further improvements regarding efficiency and parallel scalability of the
linear solver for the DG discretizations of the transport equation 
may be achieved 
by a multilevel preconditioner in which the block Gauss-Seidel method
with downwind numbering plays the role of a smoother \citep{Kanschat:2008}.

We have seen in Tables \ref{t:2D_transport_sdfem_dg_cg} and
\ref{t:3D_transport_sdfem_dg_cg} that 
the number of unknowns and therefore the matrix assembly and linear solver
times for DG($k$) grow rapidly with the order $k$ of the polynomial basis.
On quadrilateral/hexahedral meshes, where 
quadrature points and shape functions can be constructed from
a tensor product of 1-D objects, an excellent boost in performance can 
be achieved for the matrix assembly part with a technique called 
sum-factorization \citep{Melenk:2001}.

\section*{Appendix}

\begin{table}[H]
  \centering
  \scalebox{0.71}{
    % \footnotesize
    \ra{1.3}
    \hspace{-3mm}
    \begin{tabular}{@{}p{3.5cm}p{3cm}p{3cm}@{}}
      \toprule
      Machine name        & {\tt quadxeon4}            & {\tt fna} \\
      Number of nodes     &                   1        & 1  \\
      RAM per node        & 1024 GB (DDR-3/1066)       &  128GB (DDR-3/1333 MHz)  \\
      CPU-sockets per node  &                   4        &   4 \\
      Total $\#$cores     &                  40        & 48 \\
      OS                   & Debian GNU 7       &  Debian GNU 7  \\
      \midrule
      CPU socket     &
      {\it Intel\textregistered ~Xeon\textregistered ~E7-4870}
      &
      {\it AMD ~Opteron\texttrademark ~6172 }%~Interlagos}
      \\
      Clock speed  &        2.40 GHz     &   2.10 GHz   \\
      $\#$cores    &            10       &   12         \\
      $\#$threads  &            20       &   12         \\
      Launch date  &        Q2/2011      &    Q1/2010        \\
      % L1 Cache          &    32 KB per core    &   48 KB per core &   48 KB per core \\
      % L2 Cache          &   256 KB per core    &    1 MB per core &    1 MB per core \\
      L3 Cache     &               30 MB  &  12 MB \\
      $\#$memory channels &             4 &  4 \\
      % Max. memory bandwidth     &    59.7 GB/s & \\
      % QPI Speed         &      8.0 GT/s        &  6.4 GT/s  \\
      % \midrule
      % InfiniBand        &                      & & \\
      \bottomrule
    \end{tabular}
  }
  \caption{Single-node multi-core machines at the IWR Heidelberg}
  \label{t:IWRmachines}
\end{table}

{\small \noindent
\subsection*{Acknowledgements}
This study has been funded by the Baden-W\"urttemberg Stiftung in its
high-performance computing program, contract HPC-8.
}

{\small
\bibliographystyle{elsarticle-harv}
\bibliography{bib2}
}
\end{document}